\documentclass[reqno,a4paper,11pt]{amsart}
\usepackage{amsmath,amsthm,amssymb,mathrsfs}
\usepackage[
left = 2.25cm,
right = 2.25cm,
top = 2.2 cm,
bottom = 2.2 cm
]{geometry}
\usepackage{graphicx,xcolor,tikz,mathtools}
\usepackage{bm,stmaryrd}
\usepackage{amsaddr}
  \usepackage{enumerate}
  \usepackage[noadjust]{cite}
\usepackage[toc,page]{appendix}
\usepackage[
pdfencoding=auto,
colorlinks = true,
citecolor = blue,
linkcolor = blue,
anchorcolor = blue,
urlcolor = blue, 
filecolor=blue,
pdfkeywords={}
]{hyperref}
\usepackage{esint} 

\usepackage{todonotes}
\makeatletter
\define@key{todonotes}{Helmut}[]{%
	\setkeys{todonotes}{author=\textbf{Helmut},inline,color=red!10}}%
\define@key{todonotes}{Andrea}[]{%
	\setkeys{todonotes}{author=\textbf{Andrea},inline,color=green!20}}%
\define@key{todonotes}{Yadong}[]{%
	\setkeys{todonotes}{author=\textbf{Yadong},inline,color=cyan!20}}%
\makeatother

\theoremstyle{plain}
\newtheorem{theorem}{Theorem}[section]

\newtheorem{definition}[theorem]{Definition}
\newtheorem{lemma}[theorem]{Lemma}
\newtheorem{proposition}[theorem]{Proposition}
\newtheorem{remark}[theorem]{Remark}
\theoremstyle{definition}

\numberwithin{equation}{section}

\newcommand{\bv}{\mathbf{v}}

\newcommand{\bn}{\mathbf{n}}

\newcommand{\bH}{\mathbf{H}}

\newcommand{\bV}{\mathbf{V}}

\newcommand{\cH}{\mathcal{H}}
\newcommand{\cI}{\mathcal{I}}

\renewcommand{\d}{\mathrm{d}}
\newcommand{\dx}{\,\d x}
\newcommand{\dt}{\,\d t}

\newcommand{\ddt}{\frac{\d}{\d t}}

\newcommand{\ptial}[1]{ \partial_{#1} }
\newcommand{\pt}{\ptial{t}}

\newcommand{\abs}[1]{\vert #1 \vert}

\newcommand{\norm}[1]{\left\Vert #1 \right \Vert}

\newcommand{\normmm}[1]{\left\vert #1 \right\vert}

\def\Div{\mathrm{div}\,}
\DeclareMathOperator*{\esssup}{\mathrm{esssup}}
\newcommand{\dist}{\mathrm{dist}}

\allowdisplaybreaks[3]

\def \H{\mathcal{H}}
\def \L{\mathcal{L}}
\def \N{\mathbb{N}}
\def \R{\mathbb{R}}

\def \S{\mathbb{S}}
\def \Tstrong{T_{\mathrm{strong}}}

\newcommand {\supp} {\mathop \textup{supp}}

\renewcommand {\div} {\mathop \textup{div}}

\newcommand {\p} {\partial}
\newcommand {\eps} {\varepsilon}
\def\tto{\downarrow}
\def\lambdas{\lambda^*}
\def\Er{E[\mu,\chi|\mA](t)}
\def\Erel#1{E[\mu,\chi|\mA]({#1})}
\def\Eb{E_{\mathrm{bulk}}[\chi|\mA](t)}
\def\Ebulk#1{E_{\mathrm{bulk}}[\chi|\mA]({#1})}
\def\mA{\mathscr{A}}
\def\media{\langle\pi_{x,t}\rangle}
\def\mediak{\langle(\pi_n)_{x,t}\rangle}
\def\pA{\p^*A}
\def\bS{\mathbb S^{d-1}}
\def\cU{\mathcal{U}}
\def\Sigmas{\cI}
\def\si{s_{\cI}}
\def\vtheta{\vartheta}
\def\extra#1{{\color{black}#1}}
\def\RevX#1{{\color{black}#1}}
\def\RevY#1{{\color{black}#1}}
\newcommand{\lcorner}{%
	\,\raisebox{-.129ex}{\reflectbox{\rotatebox[origin=br]{-93}{$\lnot$}}}\,%
}

\def\xhi{\chi}
\allowdisplaybreaks[4]

\begin{document}
	
	\title[Varifold solutions to Volume-Preserving Mean Curvature Flow]{Varifold solutions\\ to Volume-Preserving Mean Curvature Flow:\\ existence and weak-strong uniqueness}
	
	\author[A. Poiatti]{
		\small
		Andrea Poiatti$^\ddagger$
	}

	\address{
		$^\ddagger$Faculty of Mathematics, University of Vienna, 
        1090 Vienna, Austria
	}
	\email{andrea.poiatti@univie.ac.at}

	

	\begin{abstract}
		In this contribution we introduce a novel weak solution concept for two-phase volume-preserving mean curvature flow, having both properties 
 of unconditional global-in-time existence and weak-strong uniqueness. These solutions extend the ones proposed by Hensel-Laux [J. Differential Geom. 130, 209-268 (2025)] for the standard mean curvature flow, and consist in evolving varifolds coupled with the phase volumes
by a transport equation. First, we show that, in the same setting as in Takasao [Arch. Ration.
Mech. Anal. 247, 52 (2023)], any  sharp interface limit of solutions to a slightly modified nonlocal Allen-Cahn equation is a varifold solution according to our new definition. Secondly, we crucially introduce a new  notion of
volume-preserving gradient-flow calibrations, allowing the extended velocity vector field to point in the normal direction on the interface. We show that any sufficiently regular strong solution is calibrated in this sense. Finally, we prove that any classical solution to volume-preserving mean curvature
flow, which is then automatically a calibrated flow, is unique in the class of our new varifold solutions. 
	\end{abstract}
	
	\maketitle
	
	
	\section{Introduction}
Volume-preserving mean curvature flow is the most simple geometric evolution equation for closed hypersurfaces with the main property of preserving the enclosed volume. Let $d \geq 2$ be an integer and let $\Omega$ be either $\R^d$ or the $d$-dimensional flat torus $\mathbb T^d$. Consider a family $\mA(t)$, $t\in[0,T)$, $T\in(0,\infty]$, of open sets with smooth boundary. The family of
hypersurfaces $\{\partial\mA(t)\}_{t\in[0,T)}$
is a volume-preserving mean curvature flow if the normal
velocity vector $\mathbf V$ satisfies
\begin{align}
\mathbf V=\bH-\left(\frac{1}{\H^{d-1}(\partial\mA(t))}\int_{\partial \mA(t)}\bH\cdot \bn_{\partial\mA}\d\H^{d-1}\right)\bn_{\partial\mA},
    \label{velox1}
\end{align}
for $t \in (0, T)$. Here, $\H^{d-1}$
is the $(d - 1)$-dimensional Hausdorff measure, whereas $\bH$  and $\bn_{\partial\mA}$ are the mean curvature vector and the inner unit normal vector of $\partial\mA(t)$, respectively. The fundamental property of this flow is that the enclosed volume is preserved, in the sense that
\begin{align}
\ddt \mathcal{L}^d(\mA(t))=-\int_{\partial\mA(t)}\bV\cdot \bn_{\partial\mA}\d\H^{d-1}=0,\quad \forall t\in[0,T),
    \label{preservation}
\end{align}
where $\mathcal L^d$
is the $d$-dimensional Lebesgue measure.
When $\mA(0)$ is convex, Gage \cite{13Taka} and Huisken \cite{16Taka} proved that there exists a solution to \eqref{velox1}
 converging to a sphere as time $t\to\infty$. Also, Escher and Simonett \cite{EscherSimonett} showed the short time
existence of the solution to the volume-preserving mean curvature flow for smooth initial data $\mA(0)$. They additionally proved that if
$\mA(0)$ is sufficiently close to a sphere in a suitable sense, then there
exists a global solution and it converges to some sphere as $t\to\infty$. Nevertheless, for more general smooth initial data singularities may also appear in finite time, even in the case of planar curves \cite{29Tim}. To be able to describe the evolution through these singularities, several authors have been interested in finding suitable notions of weak solutions. Mugnai, Seis and Spadaro \cite{MugnaiSeisSpadaro} constructed solutions using an almost volume-preserving
version of the scheme considered by Luckhaus and Sturzenhecker \cite{LukSturz}. The main drawback is that this result is only conditional to the energy convergence assumption
\begin{align}
\int_0^T\H^{d-1}(\partial^*\mA^k(t))\dt\to \int_0^T\H^{d-1}(\partial^* \mA(t)) \dt
    \label{enconvass}
\end{align}
as $k\to \infty$, where $\partial^*$ is the reduced boundary of the set, and $\{\mA^k(t)\}_{t\in[0,T)}$ is the time-discretized approximate solution to \eqref{velox1}. This is a reasonable assumption ruling out the piling up of layers of surfaces, and thus it guarantees  that the limit varifold preserves unit density. On the other hand, in general this property is hard to prove, since the framework of sets of finite perimeter intrinsically lacks compactness, and \eqref{enconvass} can so far only be shown for some specific two-phase cases, as for instance when the scalar mean curvature $H$ is nonnegative (see, e.g., \cite{DeFilLaux}). Convergence assumptions as in \eqref{enconvass} are quite common in the literature, and have been adopted, for instance, by Laux and Simon \cite{LauxSimon} to show the convergence of the nonlocal multi-phase Allen-Cahn
equation to the volume-preserving multi-phase mean curvature flow. 
As a different notion of weak solution, although volume-preserving mean curvature flow does not
satisfy a comparison principle, Kim and Kwon \cite{Tim19} proved the existence of a viscosity solution
to \eqref{velox1}, but only for the case when $\mA(0)$ satisfies a specific geometric condition called $\rho$-reflection. Nevertheless, this notion, albeit unique in the classical (nonpreserving) mean curvature flow (\cite{Tim12,Evans}), for the case of volume-preserving mean curvature flow is still not known to be unique (see \cite{Tim19}). 

Using yet another notion of weak solution, volume-preserving mean curvature flow can also be formulated for evolving varifolds as an extension of Brakke’s notion \cite{Brakke} of mean curvature flow
to the volume-preserving case. This notion has the important advantage that it is \textit{not} conditional to any energy convergence assumption. Takasao \cite{Tim34} first showed that solutions to a slightly modified
version of the nonlocal Allen-Cahn equation due to Golovaty \cite{Golovaty} converge to this varifold solution. More recently, the same author \cite{Takasao} refined his methods by relaxing the volume
constraint in the approximation and only recovering the precise volume preservation in the
sharp-interface limit. This allowed him to extend his result \cite{Tim34} to higher dimensions. Furthermore, in the recent paper \cite{CT} the authors extend the result in \cite{Takasao} by showing that such solution satisfies a suitable Brakke inequality in the case of periodic boundary conditions. Observe that Brakke flow is an appropriate
weak solution since it is also well suited for analyzing the multi-phase mean curvature flow.

In this contribution we propose a new weak solution concept for volume-preserving mean curvature flow, which is the first to enjoy both unconditional existence and weak-strong uniqueness. This notion exploits the gradient-flow structure of the flow, mirroring the concepts known for general gradient flows as introduced by De Giorgi \cite{TH15,TH4},
and Sandier and Serfaty \cite{SandierSeraty} (see also \cite{AmbrosioGigliSavare}). Indeed, note that the
volume-preserving mean curvature flow has a gradient-flow structure, as it is seen from the
energy dissipation relation
\begin{align*}
    \ddt \H^{d-1}(\partial\mA(t))=\int_{\partial\mA(t)}\bV\cdot\bH\d\H^{d-1}=-\int_{\partial\mA(t)}\normmm{\bV}^2\d\H^{d-1}.
\end{align*}
This weak notion of solution has been first introduced in the pioneering \cite{LH} for the standard (nonpreserving) mean curvature flow, and our aim is here to show that this notion is flexible enough to incorporate the volume-preserving constraint. 

As our first main result, in the present work we prove the existence of a global weak solution in this new sense, \textit{without} any assumption on the energy convergence as \eqref{enconvass}, when $\Omega=\mathbb T^d$. Namely, we show that solutions to the same modified nonlocal Allen–Cahn equation used by Takasao \cite{Takasao}, up to subsequences, converge to solutions of our new solution concept as the interface parameter vanishes. Our approach is based on \cite{Takasao}, which is itself rooted in the fundamental work of Ilmanen \cite{Ilmanen}. Since the existence
proof relies on the maximum principle, it is so far limited, as in \cite{LH}, to the two-phase volume-preserving mean curvature flow. Nevertheless, the notion of De Giorgi-type varifold solution can be still extended to multi-phase flows (see \cite{LH}), even though so far the existence proof is not yet available.

As a final comment, we point out that varifold solution concepts are nowadays well established in the literature of geometric flows. For instance, Abels proved the existence of varifold solutions to the two-phase Navier–Stokes equation with surface tension \cite{Abels}. More recently, a gradient flow approach à la De Giorgi, taking inspiration from \cite{LH}, has been extended to treat varifold solutions to Mullins-Sekerka flow \cite{SH}, as well as in the case of Navier-Stokes-Mullins-Sekerka equations \cite{AbelsPoiatti}.

In the second part of the paper, which deals either with $\Omega=\R^d$ or $\Omega=\mathbb T^d$, we show that the De Giorgi varifold solution we introduced does not allow for unphysical non-uniqueness before the appearance of singularities. In particular, we show that as long as a classical solution to volume-preserving mean curvature flow exists, any De Giorgi-type varifold solution in the sense of the present work coincides with this classical flow. The proof of this second result is based on a relative entropy approach. 

First, we need to extend the notion of gradient-flow calibrations introduced in the recent work by Fischer, Hensel, Laux and Simon \cite{FHLS} to the volume-preserving case and show that any sufficiently regular classical solution is calibrated in this sense (see Theorem \ref {thmcalibration}). The same task has already been performed in \cite{Timkroemer,Timvol} in different ways. The main technical issue is the extension $B$ of the velocity field in the definition of gradient-flow calibrations. In place of the trivial \textit{ad hoc} extension by nearest-point projection onto the classical solution, in \cite{Timvol} Laux solves a Neumann-Laplace equation to ensure that $B$ also satisfies the incompressibility condition $\div B=0$, at most with a linear error moving away from the (smooth) interface. This is enough to show a weak-strong uniqueness principle in the class of distributional solutions to volume-preserving mean curvature flow \cite{Timvol}. On the other hand, in \cite{Timkroemer} the authors need to relax the divergence constraint of $B$, in order to additionally fix the $\bn_{\partial\mA}\otimes \bn_{\partial\mA}$ component of the Jacobian $\nabla B$. To this aim they add a tangential part to the velocity field of the interface, so that $B=(\bV\cdot\bn_{\partial\mA})\bn_{\partial\mA}+X$ on the interface, where $X$ is a nontrivial tangential vector field. This modification of the calibration allows the authors to prove a quantitative convergence result of the nonlocal Allen-Cahn equation, in the formulation of Golovaty \cite{Golovaty}, to the volume-preserving mean curvature flow. 

Nevertheless, both the calibrations have the main drawback that the vector field $B$ in general does not point in the normal direction on the interface. On the other hand this property is necessary in the proof of weak-strong uniqueness of our new De Giorgi-type varifold solution concept, since we want to consider a truly general (oriented) varifold in our weak solution notion, without any assumption on its integrality. Indeed, to find weak-strong uniqueness without asking $B$ to be normal on the interface we would need to apply Brakke's orthogonality theorem \cite{Brakke} to the weak solution, which is only possible for an integral varifold. Therefore, we propose a third gradient-flow calibration of the volume-preserving mean curvature flow (see Definition \ref{calibration}). More precisely, we construct $B$ such that it solves an incompressible Stokes problem with nonhomogeneous Dirichlet boundary conditions in the domain enclosed by the interface, so that it automatically satisfies the incompressibility constraint $\div B=0$. Not only this, but the Dirichlet boundary conditions allow to directly impose that $B=(\bV\cdot\bn_{\partial\mA})\bn_{\partial\mA}$ on the interface, giving the desired direction of $B$ pointing in the normal direction on the interface. Surprisingly enough, in the context of gradient-flow calibrations, this is the first time that the vector field B is assumed to directly solve a vectorial elliptic problem rather than to be a (differential) transformation of a solution to a scalar elliptic one (like, e.g., the gradient of the solution to a scalar Neumann-Laplace equation).
Our second main result of the present contribution is thus the proof that any sufficiently smooth classical solution is calibrated in this new sense.

Making use of this new gradient-flow calibration, we can now pass to our third result, i.e., the weak-strong uniqueness of our new De Giorgi varifold solution concept, which is proven through a relative entropy technique. 
The basic idea of this approach is to construct a quantity which on one side captures the distance between a weak and a classical solution in a sufficiently strong sense, and on the other side allows for an estimate on its time evolution. This is the relative entropy (see the definition in \eqref{eq:relativeEntropyIntro}). This approach was recently widely adopted to study geometric evolution equations. For example, Fischer and Hensel showed in \cite{FischerHensel} a weak-strong uniqueness principle for Navier-Stokes two-phase flow with surface tension based on a suitable relative entropy functional. Other examples are the above mentioned proof of weak-strong uniqueness for distributional solutions to volume-preserving mean curvature flow by Laux \cite{Timvol}, as well as the recent result on the weak-strong uniqueness of two dimensional Mullins-Sekerka flow \cite{MullinsWS} in the case of the De Giorgi varifold solutions introduced in \cite{SH}. 

The method can also be used to
prove quantitative convergence of the Allen-Cahn equation to mean curvature flow as
shown by Fischer, Laux and Simon \cite{FischerLauxSimon}, \extra{see also \cite{FM}}. As already discussed, this method has also been extended in \cite{Timkroemer} to the convergence of a suitable nonlocal Allen-Cahn equation to a volume-preserving mean curvature flow. Also, the extension of the relative entropy technique to sufficiently
regular, strongly convex anisotropies in both surface tension and mobility function is possible, and we refer, for instance, to \cite{Timnotes,EP, LauxClemens}. \extra{We also refer to \cite{AMP} for the extension of the arguments to some non-isothermal models.}

The strength of this method is that one can also lower the assumptions on the weak solution, say as long as it satisfies a suitable optimal energy-dissipation relation. This is the case of the weak-strong uniqueness principle for De Giorgi varifold solutions to both 2D Mullins-Sekerka flow \cite{SH} and mean curvature flow \cite{LH}, and it highlights the importance of a gradient flow structure. In our third main result we thus follow the relative entropy approach to obtain a weak-strong uniqueness principle for our De Giorgi-type varifold solution to volume-preserving mean curvature flow, extending on one side the result of \cite{LH} to this more involved case, and on the other side lifting the weak-strong uniqueness result of \cite{Timvol} to the varifold setting. We point out that this is the first (weak-strong) uniqueness principle for volume-preserving mean curvature flow which relies on both a weak and strong solution whose existence is unconditionally proven. In the case of the standard (nonpreserving) mean curvature flow Stuvard and Tonegawa \cite{StuvardTonegawa} have shown the existence of a weak
solution, called the canonical Brakke flow, satisfying a transport equation in addition to Brakke’s family of inequalities, and thus enjoying further regularity properties, including weak-strong uniqueness. Nonetheless this is not the case for the volume-preserving mean curvature flow, for which to the best of our knowledge the notion of solution which is closer to a canonical Brakke flow is the one proposed in \cite{CT}, which so far is not enough to entail any weak-strong uniqueness principle.
Not only this, but our result is in the direction of finding a minimal solution concept that brings sufficient information to uniquely describe volume-preserving mean curvature flow, whereas the notion of solution proposed in \cite{StuvardTonegawa}, even if fully extended to the volume-preserving case, would be more regular than our new notion of solution à la De Giorgi.

   \subsection*{Structure of the paper}
The rest of the paper is organized as follows. In Section \ref{sec:notation} we give some preliminary notation which will be used in the sequel. Then Section \ref{sec:main} is devoted to the presentation and discussion of the three main results. Starting with the sections dedicated to the proofs, in Section \ref{sec:A} we prove Theorem \ref{thm:existence}, whereas in Section \ref{sec:B} we show the proof of Lemma \ref{lemma:further}. In Section \ref{sec:calibration} we prove the existence of the novel gradient-flow calibration for sufficiently smooth volume-preserving mean curvature flow, which leads to Section \ref{sec:C}, dedicated to the proof of Theorem \ref{theo:weakStrongUniqueness}.  In conclusion, in Appendix \ref{App:extra} we give an alternative simpler proof of Theorem \ref{theo:weakStrongUniqueness} under an additional assumption of integrality of the weak solution varifold.
	\section{Preliminaries}
  \label{sec:notation}
 We start by setting some basic notation.
 \begin{enumerate}
 \item \textbf{Lebesgue and Hausdorff measures.} Assume $\Omega$ to be  $\R^d$, $d=2,3$ or the $d$-dimensional flat torus $\mathbb T^d$. Given a Borel set $A\subset \Omega$, we denote by $\mathcal L^d(A)$ and $\H^{d-1}(A)$ the $d$-dimensional Lebesgue measure of $A$ and its $(d-1)$-dimensional Hausdorff measure, respectively. Given a Borel set $I\subset \R$, we also set $\mathcal L^1(I)$ to be the corresponding one-dimensional Lebesgue measure.
 	\item \textbf{Notation for general Banach spaces.} 
 	For any normed space $X$ of scalar-valued functions, we denote its norm by $\|\cdot\|_X$ and its dual space by $X'$. Given any vector space $X$, corresponding spaces of vector-valued or matrix-valued functions with each component belonging to $X$ are denoted by $\mathbf{X}$.
 	
 	\item \textbf{Lebesgue and Sobolev spaces.} 
 	Assume $\Omega$ to be  $\R^d$, $d=2,3$ or the $d$-dimensional flat torus $\mathbb T^d$.
 	For $1 \leq p \leq \infty$ and $k \in \N$, the canonical Lebesgue and Sobolev spaces defined on $\Omega$ are denoted by $L^p(\Omega)$ and $W^{k,p}(\Omega)$, and their standard norms are denoted by $\|\cdot\|_{L^p(\Omega)}$ and $\|\cdot\|_{W^{k,p}(\Omega)}$, respectively.
 	In case $p = 2$, we use $H^k(\Omega) = W^{k,2}(\Omega)$, and we write $L^2(\Omega)$ for $k=0$. Moreover, for any interval $I\subset\R$, any Banach space $X$, $1 \leq p \leq \infty$ and $k \in \N$, we write $L^p(I;X)$, $W^{k,p}(I;X)$ and $H^{k}(I;X) = W^{k,2}(I;X)$ ($L^p_{loc}(I;X)$, $W^{k,p}_{loc}(I;X)$ and $H^{k}_{loc}(I;X) = W^{k,2}_{loc}(I;X)$, respectively) to denote the (respectively, locally integrable) Bochner-Sobolev spaces of functions with values in $X$. The norms are indicated by $\|\cdot\|_{L^p(I;X)}$, $\|\cdot\|_{W^{k,p}(I;X)}$ and $\|\cdot\|_{H^k(I;X)}$, respectively. 

If $Y=X'$ is the dual space of $X$, then $L^\infty_{w*}(0,T;Y)$, $T>0$, is the space of all functions $\eta: (0,T)\to Y$ which are weakly* measurable and essentially bounded, meaning that the map
$$
t\mapsto \langle \eta(t),f\rangle_{Y,X}
$$ 
is measurable for any $f\in X$, and 
\begin{align}
\norm{\eta}_{L^\infty_{w*}(0,T;Y)}:=\esssup_{t\in(0,T)}\norm{\eta(t)}_Y.
	\label{norm1}
\end{align}
Recalling (see, e.g., \cite{AmbrosioFuscoPallara}) that $BV(\Omega)$ is the dual of a separable Banach space $X$, we obtain  $L^\infty_{w*}(0,T;BV(\Omega))=(L^1(0,T;X))'$ and that uniformly bounded sets in $L^\infty_{w*}(0,T;BV(\Omega))$ are weakly* compact. 
    
 	\item \textbf{Spaces of continuous functions.}
 	For  Banach spaces $X$ and $Y$, $C(X;Y)$ is the space of continuous functions mapping from $X$ to $Y$. Then, we denote by $C_c^k(X;Y)$ the spaces of $k$-times continuously differentiable functions with compact support with values in $Y$. If $Y=\R$ we omit its dependence. 
 	
 	\item \textbf{Radon measures, set of finite perimeters, and functions of bounded variation.}
 Assume $\Omega=\R^d$ or $\Omega=\mathbb T^d$, and consider any interval $I\subset\R$.
 We denote by $\mathcal{M}(\Omega \times I)$ the set of (signed) Radon measures and by $\langle \cdot, \cdot  \rangle$ the corresponding duality pairing. We denote by $BV(\Omega;\{0,1\})$ the set of functions of bounded variation in $\Omega$ with values in $\{0,1\}$.  
 For a set $A \subset \R^d$, $\chi_A$ is the associated characteristic function.
 We then say that $A$ is a set of finite perimeter in $\Omega$ if $\chi_A\in BV(\Omega;\{0,1\}),$ and  its distributional derivative in $\Omega$ is denoted as $\nabla \chi_A$. Its total variation, which is a positive Radon measure, is indicated by $|\nabla \chi_A|$. We use $\bn:=\frac{\nabla \chi_A}{|\nabla \chi_A|}$ for the measure-theoretic inner unit normal of $A$. Observe that we denoted by $\frac{\nabla \chi_A}{|\nabla \chi_A|}$ the Radon-Nikod\'ym derivative of $\nabla\chi_A$ with respect to its total variation. We also denote $\partial^\ast A$ as the reduced boundary of the set $A$ in $\Omega$. 
 For more details on functions of bounded variation and sets of finite perimeter, we refer to \cite{AmbrosioFuscoPallara,Maggi}.
 \item \textbf{Measures and varifolds.} For a locally compact and separable metric space $X$, we denote
	by $\mathcal M(X)$ the space of (finite) Radon measures on $X$. An oriented varifold $\mu$ on $\Omega$ is a positive Radon measure $\mu \in \mathcal M(\Omega\times\mathbb S^{d-1})$. The associated varifold mass measure is denoted by $\normmm{\mu}_{\mathbb S^{d-1}}\in \mathcal M(\Omega)$, with $\normmm{\mu}_{\mathbb S^{d-1}}(A):=\mu(A\times \mathbb S^{d-1})$, for any Borel-measurable $A\subset \Omega$. The varifold $\mu$ is
$(d-1)$-rectifiable if $\normmm{\mu}_{\mathbb S^{d-1}}$ is $(d-1)$-rectifiable. Also, we say that $\mu$ is $(d - 1)$-integer-rectifiable if, additionally, the $(d - 1)$-density $\theta^{d-1}(\normmm{\mu}_{\mathbb S^{d-1}},\cdot)$ of $\normmm{\mu}_{\mathbb S^{d-1}}$ is integer valued.  
An example of an oriented varifold is the measure $\normmm{\nabla\xhi}\otimes \delta_{\frac{\nabla\xhi}{\normmm{\nabla\xhi}}}$, associated to a set of finite perimeter represented by means of the characteristic function $\xhi\in BV(\Omega;\{0,1\})$. 
 \end{enumerate} 

  \section{Main results on De Giorgi-type varifold solutions \\for volume-preserving mean curvature flow}
  \label{sec:main}
  \subsection{Definition, existence results, and consistency}
We first extend the novel weak solution
  concept for two-phase mean curvature flow introduced in \cite{LH} to the case of volume-preserving mean curvature flow.
 Formally, given a smoothly evolving open set $A(t)$, $t\geq0$, this notion defines a solution only by means of the optimal energy-dissipation relation
  \begin{align}\label{eq:de giorgi inequality formal}
  	\frac{d}{dt}\H^{d-1}(\partial A(t)) \leq -\frac12 \int_{\partial A(t)} V^2 \d\H^{d-1} -\frac12 \int_{\partial A(t)} |\mathbf{H}-\lambda \bn|^2 \d\H^{d-1},
  \end{align}
  where $\bn$ denotes the inner unit normal, and the Lagrange multiplier is given by
  \begin{align*}
  \lambda:=\frac{\int_{\partial A(t)}\mathbf H\cdot  \bn\d\mathcal H^{d-1}}{\mathcal H^{d-1}(\partial A(t))}.
  \end{align*}
  Observe that here the curvature term appears in the energy inequality in the modified vectorial formulation $\bH-\lambda \bn$, due to the definition of the normal velocity in the volume-preserving mean curvature flow \eqref{velox1}. Notice that the scalar normal velocity $V$ and the mean curvature vector $\mathbf{H}$, as well as $\lambda$, can have very low regularity and they are suitable to be defined in a weak setting.  Since the evolution of $A(t)$ is assumed to be smooth, the volume rate of change satisfies
  \begin{align}
  	\frac{d}{dt} \mathcal{L}^d(A(t)) \underbrace{=}_{(1)} -\int_{\partial A(t)} V \d\H^{d-1}\underbrace{=}_{(2)}0,
  	\label{conservation}
  \end{align}
where here and in the sequel we use the sign convention $V>0$ for a shrinking set $A(t)$, say $\bV=V\bn$.
 First, from identity (1) in \eqref{conservation}, the resulting localized PDE  formulated using the characteristic function $\chi(x,t)=\chi_{A(t)}(x)$ then reads
  \begin{align*}
  	\partial_t \chi + V |\nabla \chi| =0,
  \end{align*}
 so that $\chi$ is transported by the vector field $\bV=V \bn$, recalling that $\bn=\frac{\nabla \chi}{|\nabla \chi|}$.
By substituting the term $|\nabla \chi|$, which would be lost in the limit in the varifold sense, the authors of \cite{LH} define
  \begin{align*}
  	\partial_t \chi + V \omega =0,
  \end{align*}
  where $\omega$ is the mass measure of the limit varifold $\mu$, and $V$ is measurable with respect to $\omega$, and not necessarily $|\nabla \chi|$. This gives rise to a new notion of normal velocity, which will be exploited also in the present contribution. 
  
  We come now to identity (2) in \eqref{conservation}. This has already been observed in \eqref{preservation}, and it is the main property of the volume-preserving mean curvature flow, since it entails volume conservation in the form
  $$
  \mathcal L^d({A(t)})=\mathcal L^d(A(0)),
  $$  
  for any $t>0$. Clearly, to give a meaning to our new notion of solution, this property has to be satisfied also by the weak solution we aim at defining.
  
   We give now our new definition of weak solution for (two-phase) volume-preserving mean curvature flow.

  
  \begin{definition}[De~Giorgi-type varifold solutions for volume-preserving mean curvature flow]
  	\label{def:twoPhaseDeGiorgiVarifoldSolutions}	Let $\Omega=\mathbb T^d$ or $\Omega=\R^d$, and let $\mu = \mathcal{L}^1\otimes (\mu_t)_{t\in (0,\infty)}$ be
  	a family of oriented varifolds $\mu_t \in \mathcal{M}(\Omega{\times}\S^{d-1})$, $t\in (0,\infty)$,
  	such that the map $(0,\infty) \ni t \mapsto \int_{\Omega{\times}\S^{d-1}} \eta(\cdot,\cdot,t) \d\mu_t$
  	is measurable for any $\eta \in L^1((0,\infty);C_c(\Omega{\times}\S^{d-1}))$. Also, consider  a family $A=(A(t))_{t \in (0,\infty)}$ of \RevX{Borel} subsets of $\Omega$
  	with finite perimeter such that the associated indicator function~$\chi(\cdot,t):=\chi_{A(t)}$,
  	$t \in (0,\infty)$, satisfies $\chi \in L^\infty_{w*}((0,\infty);\mathrm{BV}(\Omega;\{0,1\}))$.

  	Given an initial oriented varifold~$\mu_0 \in \mathcal{M}(\Omega{\times}\S^{d-1})$
  	and an initial phase indicator function $\chi_0 \in \mathrm{BV}(\Omega;\{0,1\}))$,
  	we call the pair~$(\mu,\chi)$ a \emph{De~Giorgi-type varifold solution for (two-phase)
  		volume-preserving mean curvature flow with initial data~$(\mu_0,\chi_0)$} if the following holds.
  	\begin{subequations}
  		\begin{itemize}
  			\item (Existence of a normal speed) Writing $\mu_t = \omega_t \otimes (\pi_{x,t})_{x\in\Omega}$
  			for the measure disintegration of $\mu_t$, $t \in (0,\infty)$, \RevX{where $\omega_t\in \mathcal M(\Omega)$ for a.a. $t>0$, and  $\pi_{x,t}\in \mathcal M(\mathbb S^{d-1})$ is a probability measure, i.e., such that  $\int_{\mathbb S^{d-1}}1\d \pi_{x,t}=1$, for $\omega_t$-a.a. $x\in\Omega$ and for a.a. $t>0$}, there exists $V \in L^2((0,\infty);L^2(\Omega,d\omega_t))$ encoding a normal velocity
  			in the sense of
  			\begin{align}
  				\label{eq:evolPhase}
  				\int_{\Omega} \chi(\cdot,T)\zeta(\cdot,T) \dx - 
  				\int_{\Omega} \chi_0\zeta(\cdot,0) \dx 
  				= \int_{0}^{T} \int_{\Omega} \chi \partial_t\zeta \dx dt
  				- \int_{0}^{T} \int_{\Omega} V \zeta \d\omega_t \dt
  			\end{align}
  			for almost every $T \in (0,\infty)$ and all $\zeta \in C^\infty_{c}(\Omega {\times} [0,\infty))$.
  			\item (Existence of a generalized mean curvature vector) There exists a generalized
  			mean curvature vector
  			$\mathbf{H} \in L^2_{loc}((0,\infty);L^2(\Omega,d\omega_t;\R^d))$  such that  
  			\begin{align}
  				\label{eq:weakCurvature}
  				\int_{0}^{\infty} \int_{\Omega} \mathbf{H} \cdot B \d\omega_t \dt
  				= - \int_{0}^{\infty} \int_{\Omega {\times} \S^{d-1}}
  				(I_d {-} p \otimes p) : \nabla B \d\mu_t dt
  			\end{align}
  			for all $B \in C^\infty_{c}(\Omega {\times} (0,\infty);\R^d)$.
  			\item(Square integrable Lagrange multiplier) There exists a measurable function $\lambda: (0, \infty) \to \R$ such that for any $T>0$ there exists $C_\lambda(T)>0$ satisfying
  			\begin{align}
  				\int_0^T\lambda(t)^2dt\leq C_\lambda(T).
  				\label{squareint}
  			\end{align}

  				\item (De~Giorgi-type inequality for volume-preserving mean curvature flow) A sharp energy dissipation
  			principle \`a la De~Giorgi holds  in the form 
  			\begin{align}
  				\label{eq:DeGiorgiInequality}
  				\int_{\Omega} 1 \d\omega_{T}
  				+ \frac{1}{2} \int_{0}^{T} \int_{\Omega} |V|^2 \d\omega_t \dt
  				+ \frac{1}{2} \int_{0}^{T} \int_{\Omega} |\mathbf{H}-\lambda \langle\pi_{x,t}\rangle|^2 \d\omega_t \dt
  				\leq \int_{\Omega} 1 \d\omega_{0}
  			\end{align}
  			for almost every $T \in (0,\infty)$, where $\langle\pi_{x,t}\rangle:=\int_{\mathbb S^{d-1}}p \d\pi_{x,t}$ for almost any $(x,t)\in\Omega\times(0,\infty) $.
  			\item (Volume conservation) The total volume is preserved, in the sense that 
  			\begin{align}
  				\int_{\Omega} \xhi(x,t)\dx=\int_{\Omega}\xhi_0(x)\,\d x,
  				\label{conservationvol}
  			\end{align}
  			for almost any $t>0$.
  			\item (Compatibility condition) For almost every $t \in (0,\infty)$
  			and for any $\zeta \in C^\infty_{c}(\Omega;\R^d)$ it holds
  			\begin{align}
  				\label{eq:compatibility}
  				\int_{\Omega} \zeta \cdot \, d\nabla\chi(\cdot,t) 
  				= \int_{\Omega {\times} \S^{d-1}} \zeta \cdot p \d\mu_t.
  			\end{align}
  		\end{itemize}
  	\end{subequations}
  \end{definition}
  \begin{remark}[Rectifiability of the varifold]\label{wsu}
  Since, by the definition of varifold solution given above, the mean curvature vector $\mathbf{H} \in L^2_{loc}((0,\infty);L^2(\Omega,d\omega_t;\R^d))$, it also holds that $\int_\Omega \normmm{\bH}^2\d\omega_t<+\infty$ for almost any $t>0$. Then, Allard's rectifiability criterion \cite{Allard} allows to conclude that the unoriented varifold $\widehat\mu_t$, corresponding to $\mu_t$  (identifying antipodal points $\pm p$ on $\S^{d-1}$), is $(d-1)$-rectifiable for almost any $t\in(0,\infty)$. Notice that in general we do not know anything about the integer rectifiability of the varifold. Nevertheless, we will see that Definition \ref{def:twoPhaseDeGiorgiVarifoldSolutions} is enough to guarantee the existence of a solution and a weak-strong uniqueness principle without resorting to the integrality of the varifold (to apply Brakke's orthogonality theorem \cite{Brakke}), in perfect analogy with the standard (nonpreserving) mean curvature flow studied in \cite{LH}.
  \end{remark}
  \begin{remark}[Local time integrability of the mean curvature vector $\bH$]
  	 Observe that, differently from the weak solution notion in \cite{LH}, here we can only require the generalized mean curvature to be locally integrable in time, i.e., $\mathbf{H} \in L^2_{loc}((0,\infty);L^2(\Omega,d\omega_t;\R^d))$, due to the presence of the locally integrable Lagarange multiplier $\lambda\in L^2_{loc}(0,\infty)$ in the De Giorgi inequality \eqref{eq:DeGiorgiInequality}.
  \end{remark}
  \begin{remark}
  	\label{compatib}
    \RevX{Note that $\omega_t$ coincides with the mass measure $\normmm{\mu_t}_{\mathbb S^{d-1}}$ associated with the varifold $\mu_t=\omega_t\otimes (\pi_{x,t})_{x\in \Omega}$, since, for a Borel-measurable set $A\subset \Omega$, 
\begin{align*}
   \vert \mu_t\vert_{\mathbb S^{d-1}}(A):=\mu_t(A\times \mathbb S^{d-1})=\int_A\left(\int_{\mathbb S^{d-1}}1 \mathrm{d} \pi_{x,t}\right)\mathrm{d}\omega_t=\omega_t(A),\quad \text{ for a.a. }t>0,
\end{align*}
and thus $\vert \mu_t\vert_{\mathbb S^{d-1}}=\omega_t$. This is possible as $\pi_{x,t}$ is indeed a probability measure. }

  	 \RevX{Moreover, }as observed in \cite{LH}, the compatibility condition~\eqref{eq:compatibility} entails $\, |\nabla\chi(\cdot,t)| \leq \omega_t$
  	in the sense of measures, for almost any $t>0$. As an immediate consequence, for almost every~$t \in (0,\infty)$
  	the Radon--Nikod\'ym derivative $\rho(\cdot,t) := \smash{\frac{ d|\nabla\chi(\cdot,t)|}{d\omega_t}}$
  	exists and additionally satisfies 
  	\begin{align}
  		\label{eq:inverseMultiplicity}
  		\rho(x,t) &\in [0,1] 
  		&&\text{ for } \omega_t\text{-almost every } x \in \Omega,
  		\\ \label{eq:RadonNikodymProperty}
  		\int_{\p^*A(t)} f \d\H^{d-1} &= \int_{\Omega{\times}\S^{d-1}}
  		f\rho(\cdot,t) \d\omega_t
  		&&\text{ for any } f \in C^\infty_{c}(\Omega).
  	\end{align}
  	As a consequence, we immediately infer from~\eqref{eq:compatibility} that  
  	\begin{align}
  	&\nonumber	\int_{\Omega} \zeta\cdot \, \rho(x,t)\bn(x,t)\, d\omega_t =\int_{\Omega} \zeta \cdot \bn(x,t)\, d\abs{\nabla\chi} 	=\int_{\Omega} \zeta \cdot \, d\nabla\chi(\cdot,t) 
  		\\&= \int_{\Omega {\times} \S^{d-1}} \zeta \cdot p \d\mu_t=\int_{\Omega}\zeta\cdot \langle \pi_{x,t}\rangle\, d\omega_t,
  		\label{equivalence}
  	\end{align} 
  	for almost every $t \in (0,\infty)$
  	and all $\zeta\in C^\infty_{c}(\Omega)$, so that we deduce 
  	\begin{align}
  		\rho(x,t)\bn(x,t)=\langle \pi_{x,t}\rangle
  		\label{equiv2}
  	\end{align}
  	for $\omega_t$-almost any $x$ and for almost any $t>0$. This allows to substitute the curvature term in the De Giorgi inequality with
  	\begin{align*}
  	\frac{1}{2} \int_{0}^{T} \int_{\Omega} |\mathbf{H}-\lambda \langle\pi_{x,t}\rangle|^2 \d\omega_t \dt=\frac{1}{2} \int_{0}^{T} \int_{\Omega} |\mathbf{H}-\lambda \rho \bn|^2 \d\omega_t \dt,
  	\end{align*}
  	which is in accordance, for instance, with the term appearing in the Brakke inequality proposed in \cite[Definition 2.2]{CT}.
  \end{remark}

  We now show that any limit as $\eps\downarrow0$ of the following modified nonlocal Allen--Cahn equation, with relaxed volume constraint and periodic boundary conditions,
  \begin{align}\label{eq:allen cahn}
  &	\eps\partial_t u_\eps = \eps\Delta u_\eps - \frac1{\eps} W'(u_\eps)+\lambda_\varepsilon\sqrt{2W(u_\varepsilon)}
  	\quad\text{in } \Omega \times(0,\infty),\\&
  	\lambda_\varepsilon(t)=\frac1{\varepsilon^\alpha}\int_\Omega (\phi(u_\varepsilon(x,0))-\phi(u_\varepsilon(x,t)))\dx,\quad \forall t>0,
  \end{align}
  for $\alpha\in(0,1)$, is a weak solution in the above sense, where we have introduced the well-known quantity
  \begin{align}
  	\phi(s):=\int_0^s\sqrt{2W(\tau)}\d\tau.\label{Modica}
  \end{align}
  Here $W\colon \R\to [0,\infty)$ is a classical double-well potential, which for simplicity, we can assume to be of the canonical form $W(u)=18 u^2 (u-1)^2$, where the prefactor is used to ensure, for the ease of readability, that $\phi(1)=1$. Note that the standard PDE theory implies the global existence and uniqueness of the solution to \eqref{eq:allen cahn} (see also \cite[Remark 8]{Takasao}). The standard Allen--Cahn equation is the $L^2$-gradient flow of the Modica-Mortola energy
  \begin{align}
  	E_\eps^S(u) := \int_{\Omega} \left( \frac\eps2|\nabla u|^2 + \frac1\eps W(u)\right) dx.
  \end{align}
  On the other hand, in this modified Allen-Cahn equation, proposed in \cite{Golovaty}, the total energy is a modified version of the one above, taking into account the relaxed volume constraint. Namely, we need to define the additional energy
  \begin{align}
   E^P_\eps(u):=\frac1{2\varepsilon^\alpha}\left(\int_\Omega (\phi(u(x,0))-\phi(u(x,t)))\dx\right)^2,
  \end{align}
  and, setting $E_\eps:=E_\eps^S+E_\eps^P$, observe that \eqref{eq:allen cahn} is indeed an $L^2$ gradient flow, thanks to the identity
  \begin{align*}
  	\frac{d}{dt} E_\eps(u_\eps(\,\cdot\,,t)) = - \int_{\Omega} \eps(\partial_tu_\eps(\,\cdot\,,t))^2 \dx\leq 0.
  \end{align*}
 Note that $E_\eps^S$ corresponds to the classical surface energy, whereas $E_\eps^
  P$ penalizes the
  difference between the initial volume and the one at time $t$ and it corresponds to a relaxed volume constraint which becomes stronger and stronger as $\eps\downarrow0$. 
  
   According to the approach looking for an energy dissipation inequality à la De Giorgi (see~\eqref{eq:de giorgi inequality formal}), here we need to measure the  energy dissipation in the balanced way
  \begin{align}
  	&\nonumber\frac{d}{dt} E_\eps(u_\eps(\,\cdot\,,t)) \\&= -\frac12 \int_{\Omega}  \eps(\partial_tu_\eps(\,\cdot\,,t))^2 \dx - \frac12\int_{\Omega}  \frac1\eps \left( \eps \Delta u_\eps(\,\cdot\,,t) - \frac1\eps W'(u_\eps(\,\cdot\,,t)+\lambda_\varepsilon\sqrt{2W(u_\eps(\cdot,t))}) \right)^2 \dx,
  	\label{eq:AC}
  \end{align}
 which comes from \eqref{eq:allen cahn} by straightforward computations. 
 
 Our first main result is the following, where we need to restrict ourselves to the case of $\Omega=\mathbb T^d$, since the proof is strongly based on the seminal work of \cite{Takasao}, which is developed only for the $d-$dimensional flat torus with periodic boundary conditions to \eqref{eq:allen cahn}, although we expect a similar proof to hold for the case $\Omega=\R^d$, as it is for the standard Allen-Cahn equation in Ilmanen's work \cite{Ilmanen}. To be precise, the results in \cite{Takasao} have been obtained for a double well potential $W(u)=\frac12(1-u^2)$, with wells in $\pm1$, but it is immediate to see that the same results hold for the choice of the same (up to a rescaling) potential with the two wells shifted in $0,1$. We prefer to consider the wells in $0,1$ to be coherent with the notation in \cite{LH}, which we extend to the volume-preserving case. In the following theorem we consider a specific (rather technical) class of well-prepared initial data as in \cite{Takasao}, whose precise introduction will be in Definition \ref{well-prepared} in Section \ref{sec:A}, for the sake of readability.
  
  \begin{theorem}[Convergence of the modified nonlocal Allen--Cahn equation to a De Giorgi solution]
  	 \label{thm:existence}\ 
  	Let $\Omega=\mathbb T^d$ and let $u_\eps$ denote the solution to the modified Allen--Cahn equation with well-prepared initial conditions  $u_{\eps,0}$ in the sense of Definition \ref{well-prepared}.
  	Then there exist a measurable function $\chi \colon \Omega \times (0,\infty ) \to \{0,1\} $, a measurable function $\lambda\colon (0,\infty)\to \R$, and a family of oriented varifolds $\mu = \mathcal{L}^1\otimes (\mu_t)_{t\in(0,\infty)} \in \mathcal{M}(  \Omega \times (0,\infty) \times \S^{d-1})$ with associated mass measure $\omega \in \mathcal{M}(\Omega\times (0,\infty))$, such that 
  	\begin{align}
  		&\lim_{\eps \downarrow0} u_\eps = \chi&& \text{strongly in } L^1_{loc}(\Omega\times(0,\infty)),\\&
  		\lim_{\eps\downarrow0} \lambda_\eps= \lambda,&& \text{weakly in } L^2_{loc}(0,\infty),
  		\\
  		&\lim_{\eps\downarrow0} \left(\frac\eps2|\nabla u_\eps|^2 +\frac1\eps W(u_\eps)\right) \mathcal{L}^d\otimes \mathcal{L}^1 = \omega && \text{weakly* as Radon measures.}
  	\end{align}
  	Furthermore, the pair $(\chi,\mu)$ is a varifold solution to the volume-preserving mean curvature flow in the sense of Definition \ref{def:twoPhaseDeGiorgiVarifoldSolutions}, where the initial conditions $(\mu_0,\xhi_0)$ are given by
  	\begin{align}
  		& \lim_{\eps\downarrow0} u_{\eps,0}=	\chi_0& & \text{strongly in }L^1(\Omega),\\
  		\label{energ}
  		&\lim_{\eps \downarrow 0} \left(\frac\eps2|\nabla u_{\eps,0}|^2 +\frac1\eps W(u_{\eps,0})\right) \mathcal{L}^d =\omega_0 &&  \text{weakly* as Radon measures},
  	\end{align}
  	with $\omega_0=\normmm{\nabla \xhi_0}$, and $\mu_0=\omega_0\otimes \delta_{\frac{\nabla \xhi_0}{\normmm{\nabla \xhi_0}}}.$ 
  	In addition, the unoriented varifold $\widehat\mu_t$, corresponding to $\mu_t$  (identifying antipodal points $\pm p$ on $\S^{d-1}$), is integer $(d-1)$-rectifiable for almost any $t\in(0,\infty)$. Furthermore, the energy-dissipation inequality \eqref{eq:DeGiorgiInequality} holds also with $(0,T)$ replaced by $(s,t)$ for {all} $t\in (0,\infty)$ and almost every $s<t$, with $s=0$ included, as well as the total mass is non-increasing, in the sense that  $\omega_t(\Omega) \leq \omega_s(\Omega)$ for all $0<s<t$.
  \end{theorem}
	\begin{remark}
		The validity of the De Giorgi-type inequality on $(s,t)$ for \emph{all} $t\in (0,\infty)$, is here not conditional  to some additional assumption on the second moments of the initial datum, as in \cite[Theorem 1]{LH}, since in this case we have $\Omega=\mathbb T^d$ and not $\R^d$.\label{remark1}
	\end{remark}
	\begin{remark}(Brakke's orthogonality theorem)\label{brakke}
	Notice that, since the varifold $\mu_t$ constructed in Theorem \ref{thm:existence} additionally enjoys integrality, we can apply the well-known Brakke's orthogonality of the mean curvature vector (see \cite[Chapter 5]{Brakke}) and deduce, as a straightforward consequence, that $\mathbf H$ points in the same direction of the expected value $\langle \pi_{x,t}\rangle$ (if this does not vanish), namely
		\begin{align}
			\abs{\langle \pi_{x,t}\rangle}^2\mathbf H(x,t)=(\langle \pi_{x,t}\rangle\cdot \mathbf H(x,t))\langle\pi_{x,t}\rangle,
			\label{basciWSuniq}
		\end{align} 
		for $\omega_t$-almost every $x\in \Omega$ and almost every $t>0$. As already pointed out in Remark \ref{wsu}, rectifiability of the varifold cannot in general be deduced for De Giorgi-type varifold solutions according to Definition \ref{def:twoPhaseDeGiorgiVarifoldSolutions}.
	\end{remark}
	
		In addition to Theorem \ref{thm:existence}, as in \cite{LH} we can also show the consistency of the varifold solution given in Definition \ref{def:twoPhaseDeGiorgiVarifoldSolutions} with respect to the classical solution to the volume-preserving mean curvature flow, which we now define. Let $\Tstrong \in (0,\infty)$ be a finite time horizon,
		and let $\mA=(\mA(t))_{t \in [0,\Tstrong]}$
		be a family of bounded domains in~$\Omega$, smoothly evolving,
		where the associated family of interfaces $\mathcal{I}=(\mathcal{I}(t))_{t \in [0,\Tstrong]}$
		is assumed to evolve smoothly following a volume-preserving mean curvature flow, in the sense of the following 
		\begin{definition}[Strong solution to volume-preserving mean curvature flow]
			Let $T_{strong}>0$ be a finite time horizon, and let $\mA=(\mA(t))_{t \in [0,\Tstrong]}$
			be a family of bounded domains in~$\Omega$, say $\mA(t)\subset B_{R_*}(0)$ for all $t\in[0,\Tstrong]$, for some $R_*>0$. Let also $\mathcal I(t)=\partial\mathscr A(t)$ be bounded and smoothly evolving in time. Then we say that $(\mathcal I(t))_t$ is a strong solution of volume-preserving mean curvature flow if $\mathcal I(t)$ is of class $C^{2,\alpha}$, for some $\alpha\in(0,1)$, for all $t\in[0,\Tstrong]$, with normal velocity $V_{\mathcal I}$ (with respect to the inner unit normal) of class $C^{2,\alpha}$ in space, and for all $t\in[0,\Tstrong]$ it holds 
			\begin{align}
				V_{\mathcal I}\bn_{\cI}=\mathbf H_{\mathcal I}-\lambda_{\cI}\bn_{\cI},\quad \text{ on }\cI(t),\label{Vi}
			\end{align}
			where $\bH_\cI$ is the associated mean curvature vector, $\bn_{\cI}$ is the inner unit normal, and $\lambda_\cI=\lambda_\cI(t)$ is the Lagrange multiplier corresponding to the volume constraint $\mathcal L^{d}(\mA(t))=\mathcal L^{d}(\mA(0))$, which is explicitely given by
			\begin{align}
				\lambda_\cI(t)=\frac{1}{\mathcal H^{d-1}(\cI(t))}\int_{\cI(t)}\bH_\cI\cdot\bn_\cI\d\mathcal H^{d-1}.
			\end{align}
			\label{strongsolmcf}
		\end{definition}

		 We then have the following consistency result, which holds for both $\Omega=\R^d$ and $\Omega=\mathbb T^d$, together with a sequential compactness property on the torus.
	
	\begin{lemma}\label{lemma:further}
		\begin{enumerate}[(i)]
			\item 	(Classical solutions are weak solutions)
			\label{item:classical is weak}
			If $(\partial \mA(t))_{t\in[0,T]} $ 
			is a classical solution to volume-preserving mean curvature flow for some $T\in (0,\infty]$, i.e., it satisfies Definition \ref{strongsolmcf} with $\Tstrong=T$ and it is smooth, 
			then the pair given by  $\chi(x,t):= \chi_{\mA(t)}(x)$ 
			and $\mu_t:= |\nabla \chi(\cdot,t)| \otimes 
			\delta_{\frac{\nabla \chi(\cdot,t)}{|\nabla \chi(\cdot,t)|}} $ 
			is a De~Giorgi-type varifold solution according to  
			Definition~\ref{def:twoPhaseDeGiorgiVarifoldSolutions} up to suitably restricting it to $[0,T]$.
			
			\item (Smooth weak solutions are classical solutions) 
			\label{item:weakclassic}
			If $(\chi,\mu)$ is a De~Giorgi-type varifold solution 
			such that there exists $T\in (0,\infty]$ so that $\chi(x,t)=\chi_{\mA(t)}(x)$
			and $\mu_t = |\nabla \chi(\cdot,t)| \otimes 
			\delta_{\frac{\nabla \chi(\cdot,t)}{|\nabla \chi(\cdot,t)|}}$ 
			for all $t\in [0,T]$ and some smoothly evolving family $(\mA(t))_{t\in [0,T]}$, 
			then $(\partial \mA(t))_{t\in [0,T]}$
			is a classical (smooth) solution to volume-preserving mean curvature flow satisfying Definition \ref{strongsolmcf} with $\Tstrong=T$.
\item (Sequential compactness on the torus) Let $\Omega=\mathbb T^d$, and $(\mu_n,\xhi_n)$ be a sequence of De Giorgi-type varifold solutions with initial data $(\mu_{n,0},\xhi_{n,0})$, $n\in\N$, such that the the energy $t\mapsto (\omega_{n})_t$ is nonincreasing in time for any $n\in\N$, as well as the initial energy $(\omega_n)_0$ is uniformly bounded and $\int_0^T\lambda_n(t)^2\dt$ is uniformly bounded in $n\in \N$ for any $T>0$. Then there exists a subsequence $\{n_l\}_{l\in \N}$ and a De Giorgi-type varifold solution $(\mu,\xhi)$ with initial conditions $\xhi_0=\lim_{l\to\infty}\xhi_{n_l,0}$ (in $L^1(\Omega))$ and $\mu_0=\lim_{l\to\infty}\mu_{n_l,0}$ (as Radon measures), such that $\mu_{n_l}\to \mu$ as Radon measures and $\xhi_{n_l}\to \xhi$ strongly in $L^1(\Omega\times(0,T))$ for any $T>0$.
		\end{enumerate}
	\end{lemma}	

\subsection{Gradient-flow calibration and weak-strong uniqueness}
We now address the weak-strong uniqueness of
	  De~Giorgi-type varifold solutions for (two-phase) volume-preserving
	mean curvature flow. Let $(\mu,\chi)$ be a De~Giorgi-type varifold solution
	according to Definition~\ref{def:twoPhaseDeGiorgiVarifoldSolutions}.
	
	To measure the distance
	between~$\mA$ and~$(\mu,\chi)$ we exploit the relative entropy
	\begin{align}
		\label{eq:relativeEntropyIntro}
		\Er := \int_{\Omega {\times} \S^{d-1}} 
		1 - p \cdot \xi(x,t) \d\mu_t(x,p) \geq 0, 
		\quad t \in [0,\Tstrong],
	\end{align}
	together with the bulk error
	\begin{align}
		\label{eq:bulkErrorIntro}
		\Eb := \int_{\Omega} 
		|\chi_{A(t)}(x) {-} \chi_{\mA(t)}(x)| |\vartheta(x,t)| \dx \geq 0, 
		\quad t \in [0,\Tstrong].
	\end{align}
    \RevY{We point out that in the last definition there is a slight inconsistency in the notation, since to be precise one should write the left-hand side as $E_{\mathrm{bulk}}[\xhi_A|\xhi_{\mA}](t)$. This is done for consistency with the relative entropy functional used in \eqref{eq:relativeEntropyIntro}.}
	The vector field~$\xi$ appearing in the definition \eqref{eq:relativeEntropyIntro} is a suitable
	extension of the normal vector field associated
	with the classical solution~$\mA$ (pointing inward), whereas the weight~$|\vartheta|$
	is a smooth version of $(x,t) \mapsto \min\{1,\dist(x,\mathcal{I}(t))\}$, which disappears on $\mathcal I(t)$.
	For a rigorous construction of $(\xi,\vartheta)$ here we introduce a new notion of gradient-flow calibrations in the context of volume-preserving mean curvature flow, which differs from the one first introduced in \cite{Timvol}. Indeed, to complete the weak-strong uniqueness argument we further need that the vector field $B$, approximating the velocity, is normal to the interface $\cI$. This is a completely nontrivial property which is new of this contribution.
 In particular, we will show that the flow admits a gradient-flow calibration $(\xi,B,\vartheta, \lambdas)$ in the following sense:
	\begin{definition}(Novel gradient-flow calibration for volume-preserving mean curvature flow)\label{calibration}\ \ \ \ \  \
	Let $T_{strong}>0$ and let $\mA=(\mA(t))_{t \in [0,\Tstrong]}$
	be a family of bounded domains in~$\Omega$, say $\mA(t)\subset B_{R_*}(0)$ for all $t\in[0,T_{strong}]$, for some $R_*>0$. Let also $\mathcal I(t)=\partial\mathscr A(t)$. Let $\xi,B:\Omega\times[0,\Tstrong]$, $\vartheta:\Omega\times[0,\Tstrong]\to\R$, and $\lambdas:[0,\Tstrong]\to\R$. We call the tuple $(\xi,B,\vartheta,\lambdas)$ a gradient-flow calibration for the volume-preserving mean curvature flow if the following statements are satisfied.
	\begin{itemize}
		\item [(i)]The vector field $\xi$, the scalar function $\vartheta$, and the function $\lambdas$ are such that 
		\begin{align}
		\nonumber&	\xi\in C_c^{0,1}(\Omega\times[0,\Tstrong];\R^d),\\& \vartheta\in C^{0,1}(\Omega\times[0,\Tstrong])\cap L^\infty(\Omega\times[0,\Tstrong]),\quad \lambdas\in L^2(0,\Tstrong).\label{controltheta}
		\end{align}
		Moreover it holds
		\begin{align}
			B(\cdot,t)\in C^{1,1}_c(\Omega;\R^d),\\&
			\label{regB}
		\end{align}
		as well as, for any $t\in[0,\Tstrong]$,
		\begin{align}
			\supp B(\cdot,t)\subset \supp \xi(\cdot,t)=\RevX{\{x\in \Omega:\ \dist(x,\cI(t))\leq c\}},\label{supports}
		\end{align}
        for some $c>0$,
		and 
		\begin{align}
			\sup_{t\in[0,\Tstrong]}\norm{B(\cdot,t)}_{C^{1,1}(\Omega;\R^d)}\leq C(\mA,\Tstrong),
			\label{unifcontrol1}
		\end{align}
      for some $C=C(\mA,\Tstrong)>0$.
		\item[(ii)](Vanishing divergence of $B$). The vector field $B$ satisfies, \RevX{for any $(x,t)\in\Omega\times [0,\Tstrong]$ such that $\dist(x,\cI(t))\leq c$},
		\begin{align}
			\label{div}\Div B(\cdot,t)=O(\dist(\cdot,\cI(t))),
		\end{align}
		where the constants are independent of $t\in[0,\Tstrong]$.
		\item[(iii)] (Normal direction of $B$ on the interface).  The vector field $B$ satisfies, \RevX{for any $(x,t)\in\Omega\times [0,\Tstrong]$ such that $\dist(x,\cI(t))\leq c$},
		\begin{align}
			\label{normal}(Id-\nabla s_\cI\otimes \nabla s_\cI) B(\cdot,t)=O(\dist(\cdot,\cI(t))),
		\end{align}
		where $s_\cI(\cdot,t)$ is the signed distance function with respect to $\cI(t)$, which is well defined over $\supp\xi(\cdot,t)$, and where the constants are independent of $t\in[0,\Tstrong]$.
		\item[(iv)] (Normal extension and shortness). The vector field $\xi$ extends the (interior) unit normal vector field of $\cI$, i.e., 
		\begin{align}
		\xi(\cdot,t)=\bn_\cI(\cdot,t)\quad\text{ on }\cI(t),
		\label{identityn}
		\end{align}
        \begin{align}
        \xi(x,t)=-\normmm{\xi(x,t)}\nabla s_\cI(x,t), \quad \forall (x,t)\in\Omega\times[0,\Tstrong],
            \label{xiprop}
        \end{align}
		and there exists a constant $c>0$ such that
		\begin{align}
			\label{short}
			\min\{1,\frac1{c^2}\ \dist(x,\cI(t))^2\}\leq 1-\normmm{\xi(x,t)},\quad \forall (x,t)\in\Omega\times[0,\Tstrong].
		\end{align}
		\item[(v)] (Approximate transport equations). The weight $\vartheta$ is approximately transported by $B$:
		\begin{align}
			\label{thetabasic}
			(\pt\vartheta+(B\cdot\nabla)\vartheta)(\cdot,t)=O\left(\min\{1,\frac1 c\dist(\cdot,\cI(t))\}\right),\quad \forall t\in[0,\Tstrong],
		\end{align}
		and also the length of $\xi$ is approximately transported by the flow:
		\begin{align}
			\label{transport}
			(\pt\normmm{\xi}^2+(B\cdot\nabla)\normmm{\xi}^2)(\cdot,t)=2\xi(\cdot,t)\cdot(\pt{\xi}+(B\cdot\nabla){\xi})(\cdot,t)=O(\dist(\cdot,\cI(t))^2),
		\end{align}
		\RevX{for any $(x,t)\in\Omega\times [0,\Tstrong]$ such that $\dist(x,\cI(t))\leq c$}, where the constants are independent of $t\in[0,\Tstrong]$.
		Additionally, there exists a constant $C>0$ and a function $f:\Omega\times[0,\Tstrong]\to \R$ with $\norm{f}_{L^\infty(\Omega\times(0,\Tstrong))}\leq C$, such that the vector field $\xi$ itself is almost transported by the flow $B$, say
		\begin{align}
			\label{almosttranported}
			(\pt\xi+(B\cdot \nabla)\xi+(\nabla B)^T\xi)(\cdot,t)=f(\cdot,t)\xi(\cdot,t)+O(\dist(\cdot,\cI(t))),
		\end{align} 
			\RevX{for any $(x,t)\in\Omega\times [0,\Tstrong]$ such that $\dist(x,\cI(t))\leq c$}, where the constants are independent of $t\in[0,\Tstrong]$.
		\item[(vi)] (Geometric evolution equation). It holds the approximated geometric evolution
		\begin{align}
			\label{geometric}
			B(\cdot,t)\cdot\xi(\cdot,t)+\Div\xi(\cdot,t)-\lambdas(t)=O(\dist(\cdot,\cI(t))),
		\end{align}
		 \RevX{for any $(x,t)\in\Omega\times [0,\Tstrong]$ such that $\dist(x,\cI(t))\leq c$}, where $\lambdas$ is defined as 
		\begin{align}
			\lambdas(t):=\frac{1}{\mathcal H^{d-1}(\cI(t))}\int_{\cI(t)}\Div \xi(\cdot,t)\d\mathcal H^{d-1},
			\label{lambdas}
		\end{align}
		and the constants are independent of $t\in[0,\Tstrong]$.
		\item[(vii)] (Sign and coercivity of the transported weight). It holds
		\begin{align*}
	     &\vartheta(\cdot,t)<0,\quad\text{ in }\mA(t),\\&
	     \vartheta(\cdot,t)>0,\quad\text{ in }\Omega\setminus\mA(t).
		\end{align*}
		Moreover, there exists a constant $C_\vartheta>0$, such that
		\begin{align}
			\min\{1,\frac1c\ \dist(x,\cI(t))\}\leq \normmm{\vartheta(x,t)}\leq C_\vartheta \min\{1,\frac1c\ \dist(x,\cI(t))\}.
			\label{varthetacontrol}
		\end{align} 
	\end{itemize}
	\end{definition}
	All the quantities $(\xi,B,\vartheta,\lambdas)$	in the definition have immediate interpretation. Indeed $\xi$ is	an extension of the normal vector field $\bn_\cI$,
	whereas, $B$ is an extension of the velocity vector field $V_\cI\bn_\cI$, with the fundamental additional properties that it is
	solenoidal (which is compatible with the volume-preservation of the corresponding PDE in Definition \ref{strongsolmcf}). As already observed, the novelty in this gradient-flow calibration with respect to the one in \cite{Timvol} lies in the fact that here $B(\cdot,t)$ is also normal to the interface $\cI(t)$, i.e., it coincides with the velocity vector field $V_\cI\bn_\cI$ on $\cI(t)$.
	 Then $\normmm{\vartheta}$ is simply a smooth version of $(x,t)\mapsto \min\{1,\RevX{\frac1c}\ \dist(x,\cI(t))\}$, whereas $\lambdas$ is the analog of the Lagrange multiplier $\lambda_\cI$. Note that this calibration has the peculiar property that we do not need any time regularity of the flow $B$, although, as it is somehow hidden in the proofs in \cite{Timvol, Timkroemer}, one needs that some norms of $B$ are uniformly bounded in time. 
	
	Our second main result is that such a gradient-flow calibration of Definition \ref{calibration} exists for any smooth  volume-preserving mean curvature flow. 
	
	\begin{theorem}[Existence of a gradient-flow calibration to volume-preserving mean curvature flow]\label{thmcalibration}
	 \ Let $\Tstrong\in(0,\infty)$, and let $(\mA(t))_{t\in [0,T_{strong}]}$ be a smoothly evolving family of bounded domains \RevY{(open and connected)} in $\Omega$ whose associated family of interfaces $(\mathcal I(t))_{t\in[0,T_{strong}]}$, $\mathcal I(t)=\partial \mA(t)$  evolves smoothly by a volume-preserving mean curvature flow as in Definition \ref{strongsolmcf}. Then there exists a tuple $(\xi,B,\vartheta,\lambdas)$ which is a gradient-flow calibration for volume-preserving mean curvature flow as in Definition \ref{calibration}.
	\end{theorem}
    Surprisingly enough, as in \cite{Timvol, Timkroemer}, also with this new calibration we do not need any estimate on the closeness of the respective Lagrange multipliers
	to deduce the weak-strong uniqueness result. As a byproduct this also gives a quantitative stability for De~Giorgi-type varifold solutions as in Definition \ref{def:twoPhaseDeGiorgiVarifoldSolutions}. We thus have our third main result as follows
	
	\begin{theorem}[Weak-strong uniqueness  for De~Giorgi-type varifold solutions]
		\label{theo:weakStrongUniqueness}
		Let $\Tstrong\in (0,\infty)$ be a finite time horizon, 
		and let $\mA=(\mA_t)_{t\in [0,\Tstrong]}$ be a smoothly evolving 
		family of bounded domains in~$\Omega=\R^d$ or $\Omega=\mathbb T^d$,
		whose associated family of interfaces $(\mathcal{I}(t))_{t \in [0,\Tstrong]}$
		evolves smoothly following a volume-preserving mean curvature flow as in Definition \ref{strongsolmcf}. 
		Let also $(\mu,\chi)$ be a De~Giorgi-type varifold solution
		for volume-preserving mean curvature flow with initial data~$(\mu_0,\chi_0)$ according to
		Definition~\ref{def:twoPhaseDeGiorgiVarifoldSolutions}.
		The following stability estimate holds
		\begin{align}
		 &	E[\mu,\chi|\mA](T)+E_{\mathrm{bulk}}[\chi|\mA](T)
			\leq C(T) (E[\mu_0,\chi_0|\mA(0)]+E_{\mathrm{bulk}}[\chi_0|\mA(0)]
			)
			 \label{eq:stabilityglobal}
		\end{align}
		for some constant $C=C(T)>0$ and for almost every $T \in [0,\Tstrong]$.
		
		More precisely, if it holds that $\chi_0=\chi_{\mA(0)}$
		and $\mu_0 = |\nabla\chi_0|\otimes{\delta_{\frac{\nabla\chi_0}{|\nabla\chi_0|}}}$, then for almost any $t \in [0,\Tstrong]$ we have
		\begin{align}
			\label{eq:weakStrong1}
			\chi(\cdot,t) &= \chi_{\mA(t)}
			\text{ almost everywhere in } \Omega,
			\\ \label{eq:weakStrong2}
			\mu_t &= |\nabla\chi(\cdot,t)|\otimes(\delta_{\frac{\nabla\chi(\cdot,t)}{|\nabla\chi(\cdot,t)|}(x)})_{x\in\Omega}.
		\end{align}
	\end{theorem}
	\begin{remark}
		We point out that, due to the presence of the Lagrange multiplier $\lambda$, one cannot retrieve a stability estimate for the relative entropy $E[\mu,\chi|\mA]$ decoupled from the stability estimate of $E_{\mathrm{bulk}}[\chi|\mA]$, differently from what happens with the standard mean curvature flow (see \cite[Theorem 2]{LH}).  
	\end{remark}
We will prove this Theorem in Section \ref{sec:C} by making use of the gradient-flow calibration given in Theorem \ref{thmcalibration}. On the other hand, there are cases in which the essential property \eqref{normal}  of the extended velocity vector field $B$ to be normal to the interface cannot be satisfied, as for instance in the gradient-flow calibration proposed in \cite{Timkroemer}. Also, in general we cannot expect a similar property for the gradient-flow calibration of a multi-phase mean curvature flow (due to the presence of moving triple junctions, see, e.g. \cite{LH, FHLS}). To deal with these cases, the price to pay for two-phase flows is simply a further assumption in the definition of the weak (varifold) solution, namely that the unoriented varifold $\widehat\mu_t$, corresponding to $\mu_t$  (identifying antipodal points $\pm p$ on $\S^{d-1}$), has integer multiplicity for almost any $t\in(0,\infty)$. This is not a limiting assumption, in the sense that the varifold solution constructed as a sharp interface limit in Theorem \ref{thm:existence} is even integer $(d-1)$-rectifiable. In Appendix \ref{App:extra} we will thus prove the weak-strong uniqueness principle under this extra assumption on the integrality of the varifold in the definition of weak solution, and without making use of property \eqref{normal} of the gradient-flow calibration.	
	
	\section{Proof of Theorem \ref{thm:existence}}
	\label{sec:A}
	Let $u_\eps$ be the solution to \eqref{eq:allen cahn} with given initial condition $u_\eps(0)=u_{\eps,0}$. We also set once and for all the notation
	\begin{align}
	E_0:=\sup_{\eps>0}E_\eps(u_{\eps,0})=\sup_{\eps>0}E_\eps^S(u_{\eps,0})<\infty.\label{energybound}
	\end{align}	
As anticipated when stating the main results, the essential ingredient is the energy identity \eqref{eq:AC}, which we rewrite here after a time integration in $(0,T)$, $T\in(0,\infty)$, as an inequality, as it will be used in the sequel:
  \begin{align}
	&E_\eps(u_\eps(\,\cdot\,,T))\nonumber\\& +\frac12 \int_0^T\int_{\Omega}  \eps(\partial_tu_\eps(\,\cdot\,,t))^2 \dx\dt+ \frac12\int_0^T\int_{\Omega}  \frac1\eps \left( \eps \Delta u_\eps(\,\cdot\,,t) - \frac1\eps W'(u_\eps(\,\cdot\,,t))+\lambda_\varepsilon\sqrt{2W(u_\eps(\cdot,t))} \right)^2 \dx\dt\nonumber\\&\leq  E_\eps(u_\eps,0).
	\label{eq:AC2}
\end{align}
Note that, this inequality is equivalent to the standard energy inequality 
  \begin{align}
	&E_\eps(u_\eps(\,\cdot\,,T)) + \int_0^T\int_{\Omega}  \eps(\partial_tu_\eps(\,\cdot\,,t))^2 \dx\dt\leq  E_\eps(u_{\eps,0}),
	\label{eq:AC3}
\end{align}
and, analogously, the same estimate holds for any $0<s<t$:
  \begin{align}
	&E_\eps(u_\eps(\,\cdot\,,t))\nonumber\\& + \int_s^t\int_{\Omega}  \eps(\partial_tu_\eps(\,\cdot\,,\tau))^2 \dx\d\tau+\frac12\int_s^t\int_{\Omega}  \frac1\eps \left( \eps \Delta u_\eps(\,\cdot\,,\tau) - \frac1\eps W'(u_\eps(\,\cdot\,,\tau))+\lambda_\varepsilon\sqrt{2W(u_\eps(\cdot,\tau))} \right)^2 \dx\d\tau\nonumber\\&\leq  E_\eps(u_{\eps}(\cdot,s)).
	\label{eq:AC4}
\end{align}
Let us now recall the definition of $\phi$ in \eqref{Modica}, and set
$$
\psi_\eps:=\phi(u_\eps).
$$
Then, by the standard Modica-Mortola/Bogomol'nyi trick \cite{Modica, Bogo}, we have by \eqref{eq:AC2} that
\begin{align}
\sup_{t\in(0,\infty)}\int_\Omega\abs{\nabla\psi_\eps(\cdot,t)}\dx\leq \sup_{t\in(0,\infty)}E^S_\eps(u_\eps(\cdot,t))\leq \sup_{t\in(0,\infty)} E_\eps(u_{\eps}(\cdot,t))\leq E_\eps(u_{\eps,0}). \label{basicbound}
\end{align}
We now introduce the corresponding oriented space-time varifold
\begin{align}\label{eq:def mu_eps}
	\mu_\eps := \left( | \nabla \psi_\eps| \L^d\llcorner \Omega\right) \otimes 
	\left(\L^1 \llcorner (0,\infty)\right) \otimes 
	\left(\delta_\frac{\nabla \psi_\eps}{|\nabla \psi_\eps|}\llcorner\S^{d-1}\right).
\end{align}
Additionally, we denote the associated mass measure as
\begin{align}\label{eq:def omega_eps}
	\omega_\eps :=  \left( | \nabla \psi_\eps| \L^d\llcorner \Omega\right) \otimes \left(\L^1\llcorner (0,\infty)\right),
\end{align}
i.e., as the $(x,t)$-marginal of $\mu_\eps$.
Then, we recall one of the main novelties in \cite{LH}, namely we define the approximate normal velocity $V_\eps$ and the approximate 
mean curvature vector~$\mathbf{H}_\eps$ to be the $\omega_\eps$-measurable functions
\begin{align}
	\label{eq:def Veps} V_\eps &:= -\frac{\eps \partial_t u_\eps}{\sqrt{2W(u_\eps)}},
	\\ \label{eq:Heps} \mathbf{H}_\eps &:=  -\left( \eps \Delta u_\eps - \frac1\eps W'(u_\eps)\right) \frac{\nabla \psi_\eps}{|\nabla \psi_\eps|},
\end{align}
and we recall that the Lagrange multiplier $\lambda_\eps$ is defined, for almost any $t>0$, as 
\begin{align}
	\lambda_\eps(t):=\frac1{\varepsilon^\alpha}\int_\Omega (\psi_\eps(x,0)-\psi_\varepsilon(x,t))\dx,
	\label{lambdaeps}
\end{align}
for \RevY{some fixed} $\alpha\in (0,1)$. In simple words, the idea is here to decouple the roles of $E_\eps^S$ and $E_\eps^P$, since the first variation with respect to $u_\eps$ \RevY{of $E_\eps^S$ gives rise to the definition of $\bH_\eps$, whereas the first variation of $E^P_\eps$ gives $\lambda_\eps$}.
We point out that we define $\frac{\nabla \psi_\eps}{|\nabla \psi_\eps|}:=e_1$ on the set $\{|\nabla \psi_\eps|=0\}$ (where $e_1$ is the first element of the canonical basis of $\R^d$) and $ \frac{\eps \partial_t u_\eps}{\sqrt{2W(u_\eps)}} := 0$ on the set $\{\sqrt{2W(u_\eps)} =0\} \subset \{|\nabla \psi_\eps|=0\}$, which are clearly outside the support of the measure $\omega_\eps$. 

The following compactness statement can be proved similarly to \cite[Proposition 1]{LH}, recalling \cite[Lemma 1]{Takasao}.
\begin{proposition}[Compactness properties]\label{prop:compactness}
	Let the initial conditions $u_\eps(\,\cdot\,,0)$ satisfy the energy bound~\eqref{energybound} and assume $u_{\eps,0}\to \xhi_{A_0}$ in $L^1(\Omega)$ as $\eps\downarrow0$, for some finite perimeter set $A_0\subset \Omega$.
	Then it holds
	\begin{align}
		\label{eq:nabla}\sup_{t\in(0,\infty)} \int_{\Omega}|\nabla \psi_\eps| \dx \leq& E_0,
		\\ 	\label{eq:dt}\int_0^T\int_{\Omega}|\partial_t \psi_\eps| \dx dt \leq &(1+2T)E_0 \quad \forall T<\infty.
	\end{align}
	Furthermore, there exists $\eps_1>0$ and, for any $T>0$, $C(T)>0$ depending on $T>0$ and $E_0$, but not on $\eps\in(0,\eps_1)$, such that 
		\begin{align}
		\label{lambda_bound}
		\int_0^T\lambda_\varepsilon(t)^2\dt\leq& C(T),\quad \forall \eps\in(0,\eps_1),\quad  \forall T<\infty.
	\end{align}
	 As a consequence, after possibly passing to a subsequence $\eps\downarrow0$, there exists a family of 
 finite perimeter sets $(A(t))_{t\in (0,\infty)}$, with corresponding characteristic function $\xhi(x,t):=\xhi_{A(t)}(x)$, such that 
	\begin{align}
	&	\label{eq:compactness psi} \psi_\eps \to \psi  \quad \text{in }L^1_{loc}(\Omega\times (0,\infty)),\\&
			\label{compactness_lambda} \lambda_\eps \rightharpoonup \lambda  \quad \text{weakly in }L^2_{loc} (0,\infty),
	\end{align}
	where $\lambda\in L^2_{loc} (0,\infty)$, and $\psi(x,t) :=  \chi(x,t)$ with $\phi(1) =\int_0^1\sqrt{2W(s)}\d s=1$,
	\begin{align}
	\label{classicalbound}	\operatorname{ess\,sup}_{t\in(0,\infty)} 
		\limits  \mathcal{H}^{d-1}(\partial^\ast A(t)) \leq E_0,
	\end{align}
	and $\chi\in L^\infty_{w*}(0,\infty;BV(\Omega;\{0,1\}))$.
	Additionally, the volume is conserved, i.e.,
	\begin{align}\label{eq:compactness cont vol}
	\mathcal{L}^d(A(t))=\mathcal{L}^d(A_0).
	\end{align}	
Moreover, there exist Radon measures $\mu$ and $\omega$ on $\Omega\times(0,\infty)\times \S^{d-1}$ and $\Omega\times(0,\infty)$, respectively, such that 
	\begin{align}
		\label{eq:compactness mu} \mu_\eps \stackrel{\ast}{\rightharpoonup} \mu 
		&\quad 
		\text{weakly* as Radon measures in } \Omega\times (0,\infty)\times \S^{d-1},
		\\\label{eq:compactness omega} \omega_\eps \stackrel{\ast}{\rightharpoonup} \omega& 
		\quad 
		\text{weakly* as Radon measures in } \Omega\times (0,\infty).
	\end{align}
\end{proposition}
\begin{remark}
	\label{additional}
	Takasao in \cite[Proposition 9]{Takasao} actually proves more. Indeed, if we define
	\begin{align}
			\widetilde \omega_\eps:= \left( \left(\frac\eps2\normmm{\nabla u_\eps}^2+\frac1\eps W(u_\eps)\right) \L^d\llcorner \Omega\right) \otimes \left(\L^1\llcorner (0,\infty)\right),\label{tildeomega}
	\end{align}
	 it holds 
	\begin{align}\label{more}
		(\widetilde\omega_{\eps})_t\stackrel{\ast}{\rightharpoonup}\widetilde{\omega}_t\quad \text{weakly* as Radon measures in } \Omega,
	\end{align}
	as $\varepsilon\downarrow0$, for almost any $t\geq 0$, and for some Radon measure $\widetilde \omega_t$. This convergence is not strictly necessary for the proof of Proposition \ref{prop:compactness}, but it is fundamental to prove the last part of Theorem \ref{thm:existence}. Indeed, as already observed in Remark \ref{remark1}, since $\Omega$ is the $d$-dimensional flat torus, it immediately holds
		\begin{align}
		(\widetilde\omega_{\eps})_t(\Omega)\to \widetilde\omega_t(\Omega), \label{convomega}
	\end{align}
as $\varepsilon\downarrow0$, for almost any $t\geq0$, and this is necessary to pass to the limit in the energy inequality for generic time instants $0\leq s\leq t$. We refer to Lemma \ref{lemmaslicing} for more details about this convergence.
\end{remark}
\begin{proof}
	We only give a sketch of the immediate results, since the only nontrivial one, i.e., \eqref{lambda_bound}, can be found in \cite[Proposition 1]{Takasao}, to which we refer the interested reader. Now, \eqref{eq:nabla} directly comes from \eqref{basicbound}, and, by the usual Modica-Mortola/Bogomol'nyi trick, this time applied to the time derivative of $\psi_\eps$, we also have, exploiting \eqref{eq:AC3},
	 \begin{align}
	 \nonumber
	 	&\int_0^T\int_\Omega\abs{\partial_t\psi_\eps}\dx\dt\leq \sqrt 2\left(\int_0^T\int_\Omega \frac1\eps W(u_\eps)\dx\dt\right)^\frac12\left(\int_0^T\int_\Omega \eps\abs{\partial_t u_\eps}^2\dx\dt \right)^{\frac12} \\&\leq \int_0^T\int_\Omega \frac2\eps W(u_\eps)\dx\dt+\int_0^T\int_\Omega \eps\abs{\partial_t u_\eps}^2\dx\dt\leq \int_0^T2E^S_\eps(u_\eps)+\int_0^T\int_\Omega \eps\abs{\partial_t u_\eps}^2\dx\dt\nonumber\\&\leq  (2T+1)E_\eps(u_{\eps,0})\leq (2T+1)E_0, \label{basicbound2}
	 \end{align}
	 for any $T>0$, immediately giving \eqref{eq:dt}. As \eqref{lambda_bound} is already obtained in \cite[Proposition 1]{Takasao}, we pass to show the convergences. Since $\psi_\eps$ is uniformly bounded in $BV(\Omega\times(0,T))$ for any $T>0$, by compactness we infer \eqref{eq:compactness psi}. To identify $\psi$ one can exploit the fact that 
	 $$
	 \sup_{t\in(0,\infty)}\int_\Omega W(u_\eps)\dx\leq \eps E_0, 
	 $$ 
	 and thus, thanks to the double well shape of $W$, by standard arguments one deduces that $u_\eps(t)\to \xhi_{A(t)}$ in $L^1(\Omega)$, for some finite perimeter set $A(t)$, for almost any $t>0$, and thus $\psi= \chi_A$, where we recall $\phi(1)=1$. Also \eqref{classicalbound} is retrieved in a classical way. Then, since $\psi_\eps\in L^\infty(0,\infty; W^{1,1}(\Omega))\hookrightarrow L^\infty(0,\infty; BV(\Omega))$, and, as seen in Section \ref{sec:notation}, uniformly bounded sets in $L^\infty_{w*}(0,\infty; BV(\Omega))$ are
	 weakly* compact, we have, by uniqueness of the weak* limit, 
	 $$
	 \psi_\eps\stackrel{\ast}{\rightharpoonup} \psi,\quad\text{ weakly* in }L^\infty_{w*}(0,\infty; BV(\Omega)),
	 $$ 
	 entailing that $\xhi_A\in L^\infty_{w*}(0,\infty; BV(\Omega))$. Moreover, convergence \eqref{compactness_lambda} directly comes from \eqref{lambda_bound} by weak compactness in $L^2(0,T)$ for any $T>0$. Then, to show property \eqref{eq:compactness cont vol}, observe that, from \eqref{eq:AC3}, 
	 \begin{align*}
	 	\sup_{t\in(0,\infty)}E_\eps^P(u_\eps(\cdot,t))\leq E_0,
	 \end{align*}
	 and thus, by the assumption $u_{\eps,0}\to \xhi_{A_0}$ in $L^1(\Omega)$, 
	 \begin{align}
	 	\sup_{t\in(0,\infty)}\left\vert{\int_\Omega \psi_\eps(x,t)\dx-\int_\Omega \psi_\eps(x,0)\dx}\right\vert\leq \sqrt2\eps^{\frac\alpha2} \sqrt{E_0}.\label{controlbasic1}
	 \end{align}
	 Since the left-hand side, for almost any $t>0$, converges to $(\mathcal L^d(A(t))-\mathcal L^d(A_0))$ as $\eps\downarrow0$, we infer \eqref{eq:compactness cont vol}. In conclusion, convergences \eqref{eq:compactness mu}-\eqref{eq:compactness cont vol} are immediate from the definitions of $\omega_\eps$ and $\mu_\eps$, by weak* compactness recalling \eqref{eq:nabla}. This concludes the proof.  
\end{proof}

		The following proposition can be found in \cite{Takasao}, and is rooted in the result of Ilmanen \cite{Ilmanen} based on the maximum principle. Namely, if one chooses in a suitable way the initial data, then the discrepancy  $\frac\eps2|\nabla u_\eps(x,t)|^2  - \frac1\eps W(u_\eps(x,t))$ between the two terms of the energy $E^S_\varepsilon$ will remain nonpositive for any $t\geq0$. This is well known from \cite{Ilmanen} in the standard Allen-Cahn equation, but it is not in general so trivial for other problems, as in our modified nonlocal Allen-Cahn equation. Additionally, the discrepancy vanishes in the sense of measures: 
	\begin{align}
		\left(\frac\eps2|\nabla u_\eps|^2  - \frac1\eps W(u_\eps)\right) \mathcal{L}^d\otimes \mathcal{L}^1  \stackrel{\ast}{\rightharpoonup} 0 \quad \text{weakly* as Radon measures},\label{discrepancy}
	\end{align} 
	meaning that there is equipartition of energy as $\eps\downarrow0$. Before stating the theorem, we need to introdcuce a class of well-prepared initial data as in \cite{Takasao}. First, we say that a bounded open set $U_0\subset\subset \Omega$ is an admissible initial set if 
	\begin{itemize}
		\item[(1)] There exists $D_0>0$ such that 
		\begin{align*}
			\sup_{x\in \Omega,\ 0<r<1}\frac{\mathcal H^{d-1}(\partial U_0\cap B_r(x))}{\omega_{d-1}r^{d-1}}\leq D_0,
		\end{align*}
		where $\omega_{d-1}$ is the measure of the $(d-1)$- dimensional ball.
		\item[(2)] There exists a family of open sets $\{U_0^i\}_{i=1}^\infty$ such that $U_0^i$ has $C^3$ boundary for any $i\in \N$ and it holds
		\begin{align*}
			\lim_{i\to\infty}\mathcal L^d(U_0\Delta U_0^i)=0,\quad\text{ and }\lim_{i\to\infty}\normmm{\nabla \xhi_{U_0^i}}=\normmm{\nabla \xhi_{U_0}}\ \text{weakly* as Radon measures.}
		\end{align*}
	\end{itemize}
	Observe that the second assumption is always satisfied for a finite perimeter set $U_0$, whereas both
	conditions are fulfilled when $\partial U_0$ is of class $C^1$ (see \cite{Giusti}).
	Let us now introduce the class of well-prepared initial data:
	\begin{definition}[Well-prepared initial data]\label{well-prepared}
		Given an admissible initial set $U_0$, we say that the sequence of initial data $\{u_{\eps,0}\}_\eps$ is well-prepared if
		\begin{enumerate}
			\item[(a)] For any $\varepsilon>0$ and any $x\in \Omega$, we have 
				\begin{align}
				\label{eq:Ilmanen initial 1}\frac\eps2|\nabla u_{\eps,0}(x)|^2  
				\leq \frac1\eps W(u_{\eps,0}(x)).
			\end{align}
			\item[(b)] Setting the measure 
			$$
			\widetilde{\omega}_{\eps,0}:=\left(\frac\eps2\normmm{\nabla u_{	\eps,0}}^2+\frac1\eps W(u_{\eps,0})\right)\mathcal L^d\lcorner\Omega,
			$$
			there exists $D_1(D_0)>0$ such that 
			\begin{align}
					\label{energyb1}
				\sup_{\eps>0}E_\eps^S(u_{\eps,0})=\sup_{\eps>0}\widetilde\omega_{\eps,0}(\Omega)\leq D_1,
				\end{align}
				and
			\begin{align*}
		 \sup_{\eps>0,\ x\in \Omega, r\in(0,1)}\frac{\widetilde\omega_{\eps,0}(B_r(x))}{\omega_{d-1}r^{d-1}} \leq D_1.
			\end{align*} 
			\item[(c)] It holds 
			\begin{align}
				\widetilde\omega_{\eps,0}\stackrel{\ast}{\rightharpoonup} \omega_0:=\mathcal H^{d-1}\lcorner {\partial U_0}=\normmm{\nabla \xhi_{U_0}},\quad \text{weakly* as Radon measures},
				\label{eps0conv}
			\end{align}
			as $\eps\downarrow0$.
		 
			\item[(d)] It holds \begin{align}
				&u_{\eps,0}\to \xhi_{U_0},\quad\text{ strongly in }L^1(\Omega),\label{wellprep}\\&
				\normmm{\nabla u_{\eps,0}}\stackrel{\ast}{\rightharpoonup} \omega_0,\quad \quad \text{weakly* as Radon measures},
				\label{conv2}
			\end{align}	
			as $\eps\downarrow 0$.
		\end{enumerate}
		
	\end{definition}
As shown in \cite[Proposition 2]{Takasao}, such well prepared initial data exist. We then state a slight reformulation of the results in \cite{Takasao}.
	
	\begin{proposition}[Takasao \cite{Takasao}]\label{prop:takasao}
		If the initial conditions $u_{\eps,0}$ are well-prepared in the sense of Definition \ref{well-prepared}, then it holds
		\begin{align}
			&\psi_\eps(\cdot,0)=\phi(u_{\eps,0})\to \psi(\cdot,0)= \chi_{U_0},\quad \text{strongly in }L^1(\Omega),\text{ as }\eps\tto0\label{L1conv},\\&
			\label{eq:Ilmanen initial 3} \lim_{\eps\downarrow0} E_\eps^S(u_{\eps,0}) = \omega_0(\Omega).
		\end{align}
		Also, for all $t>0$ it holds:
		\begin{align}
			\label{eq:takasao2} \frac\eps2|\nabla u_\eps(x,t)|^2  \leq \frac1\eps W(u_\eps(x,t)) &\quad \text{ for any $x\in \Omega$ and any $t\geq 0$, and} 
			\\ \label{eq:takasao 3}\sup_{\eps>0}E_\eps(u_\eps(\,\cdot\,,t))\leq D_1<+\infty& \quad \text{for all $t\geq 0$}.
		\end{align}
		Furthermore,	for any test function $\xi \in C_c(\Omega\times (0,\infty))$, it holds
		\begin{align}\label{eq:discrepancy 1}
			\lim_{\eps\downarrow 0} \int_{\Omega\times(0,\infty)} \xi \left(\frac\eps2|\nabla u_\eps(x,t)|^2  - \frac1\eps W(u_\eps(x,t))\right) dx dt =0,
		\end{align}
		as well as
		\begin{align}
			\notag&\lim_{\eps\downarrow 0} \int_{\Omega\times(0,\infty)} \frac12 \xi \left(\sqrt{\eps} |\nabla u_\eps| 
			- \frac1{\sqrt{\eps}}\sqrt{2W(u_\eps)} \right)^2 dx dt 
			\\\label{eq:discrepancy 2}&=\lim_{\eps\downarrow 0} \int_{\Omega\times(0,\infty)} \xi \left( \frac\eps2 |\nabla u_\eps|^2 
			+ \frac1\eps W( u_\eps) - |\nabla \psi_\eps|\right)dx dt =0.
		\end{align}
	\end{proposition}
	\begin{proof}
		First, notice that, from \eqref{energyb1} and since $W(u)=18 u^2(u-1)^2\geq \frac{72}5u^4$, for $u>M$, $M$ sufficiently large, it holds
		\begin{align}
			\int_\Omega u_{\eps,0}^4\dx\leq \int_{\normmm{u_{\eps,0}}\leq M}M^3u_{\eps,0}\dx+\int_{\normmm{u_{\eps,0}}> M} \frac{5}{72}W(u_{\eps,0})\dx\leq C,\quad \forall \eps>0,\label{L4}
				\end{align}
				i.e., $u_{\eps,0}$ is uniformly bounded in $L^4(\Omega)$. Also, from \eqref{wellprep} we deduce that, up to a nonrelabeled subsequence, 
				$$
				u_{\eps,0}\to \xhi_{U_0}, \quad \text{almost everywhere in }\Omega,
				$$
				as $\eps\tto0$, and, since $\phi:s\mapsto \int_0^s\sqrt{2W(\tau)}\d\tau$ is (Lipschitz) continuous, we also have $\psi_\eps(\cdot,0)=\phi(u_{\eps,0}(\cdot))\to\xhi_{U_0}=\psi(x,0)$ almost everyhere in $\Omega$. Then we have, for some $M>0$ sufficiently large, recalling \eqref{L4}, 
				\begin{align*}
					\int_\Omega \normmm{\psi_\eps(x,0)}\dx \leq \int_{\{\normmm{u_{\eps,0}>M}\}}C\normmm{u_{\eps,0}(x)}^3\dx+Lip({\phi_{|[0,M]}})\int_\Omega \normmm{u_{\eps,0}(x)}\dx\leq C<+\infty,
				\end{align*}
				so that, also applying the Modica-Mortola/Bogomol'nyi trick by means of \eqref{energyb1}, we infer $\psi_\eps(\cdot,0)\in W^{1,1}(\Omega)\hookrightarrow L^{\frac d{d-1}}(\Omega)$ uniformly in $\eps$. Since $\frac d{d-1}>1$, this allows to apply Vitali's Convergence Theorem and deduce that, up to a subsequence, $\psi_\eps(\cdot,0)\to\psi(x,0)$ strongly in $L^1(\Omega)$ as $\eps\downarrow0$. Notice that, by uniqueness of the limit, this also holds for the entire sequence $\eps$, not only up to subsequences, entailing \eqref{L1conv}. 
				
		Proceeding in the proof, convergence \eqref{eq:Ilmanen initial 1} can be easily deduced from property \eqref{eps0conv} of the initial data, since $E^S_\varepsilon(u_{\eps,0})=\widetilde\omega_{\eps,0}(\Omega)$.
		 Then property \eqref{eq:takasao2} is shown in \cite[Proposition 6]{Takasao} based on maximum principle, whereas inequality \eqref{eq:takasao 3} can be immediately inferred from the energy inequality \eqref{eq:AC3}. In conclusion, \eqref{eq:discrepancy 1} corresponds to the vanishing discrepancy measure \eqref{discrepancy}, whose proof can be found in \cite[Theorem 13]{Takasao}, and it is crucially based on the fundamental result in \cite{Ilmanen}.  Then \eqref{eq:discrepancy 2} can easily be obtained exactly as in \cite[Proof of Proposition 2]{LH}, by slightly manipulating \eqref{eq:discrepancy 1}.
	\end{proof}
	We now aim at showing that the time slices of the measure $\omega$ defined in Proposition \ref{prop:compactness} are (sharply) bounded above by the total energy. This is \RevY{analogous to} \cite[Lemma 2]{LH}. In this case it is slight more delicate, since $E_\eps^S(u_\eps(\cdot,t))$ is \textit{not} monotone nonincreasing in time, since this property is enjoyed by the total energy $E_\eps(u_\eps)$. This is actually not an issue since $E^P_\eps(u_\eps)\geq0$. One could wonder if the inequality is actually sharp, due to the presence of the energy $E_\eps^P$, and the answer is indeed affirmative, since the inequality is even an equality. Additionally, in this Lemma we also prove the compatibility condition \eqref{eq:compatibility}. In particular, we have 
	\begin{lemma}
	Under the same assumptions of Proposition \ref{prop:takasao},	the measures $\omega$ and $\mu$ constructed in Proposition \ref{prop:compactness} can be disintegrated as $\omega=\left(\mathcal{L}^1\llcorner (0,\infty)\right) \otimes (\omega_t)_t $ and $\mu=\left(\mathcal{L}^1\llcorner (0,\infty)\right) \otimes (\omega_t)_t \otimes (\pi_{x,t})_{x,t} $, where $(\omega_t)_t$ is a weakly* $\mathcal{L}^1$-measurable family of
	finite Radon measures
	and $(\pi_{x,t})_{x,t} $ is a weakly* $\mathcal{L}^1 \otimes (\omega_t)_t$-measurable family of
probability measures. Moreover, for a.e.  $t\in(0,\infty)$ it holds
	\begin{align}\label{eq:comp}
		\liminf_{\eps\downarrow0}E_\eps(u_\eps(\,\cdot\,,t)) \geq \omega_t(\Omega),
	\end{align}
	and this inequality is actually sharp, namely 
\begin{align}
	\lim_{\eps\downarrow0}E_\eps(u_\eps(\,\cdot\,,t))=\omega_t(\Omega),\label{sharp}
\end{align}
	for almost any $t>0$. 
	
	\noindent Furthermore, $\mu_t$ and $\nabla \xhi$ enjoy the following compatibility condition: for almost every $t \in (0,\infty)$
	and for any $\xi \in C^\infty_{c}(\Omega;\R^d)$, it holds
	\begin{align}
		\label{eq:compatibility2}
		\int_{\Omega} \xi \cdot \, d\nabla\chi(\cdot,t) 
		= \int_{\Omega {\times} \S^{d-1}} \xi \cdot p \d\mu_t.
	\end{align}
	\label{lemmaslicing}
	\end{lemma}
	\begin{proof}
		The proof follows the lines of \cite[Lemma 2]{LH} (see also \cite[Proof of Proposition 10]{Takasao} for a similar argument).
		First, notice that, due to the energy inequality \eqref{eq:AC4}, we infer that, for any $\eps>0$, the function $t\mapsto E_\eps(u_\eps(\cdot,t))$ is monotonically nonincreasing, and also, recalling \eqref{eq:takasao 3}, uniformly bounded for all $t\geq0$. As a consequence, by Helly's Selection Theorem there exists a subsequence $\eps\downarrow0$ and a nonincreasing function $t\mapsto E(t)$ such that 
		\begin{align}
			E_\eps(u_\eps(\cdot,t))\to E(t),\quad \forall t\geq 0,\label{energia}
		\end{align}
		as $\eps\tto0$. We now follow \cite{LH} and consider a cylindrical test function, namely, given $\eta\in C_c(0,\infty)$, with $\eta\geq 0$, and $\zeta\in C_c(\Omega)$, $\zeta\in[0,1]$, it is easy to infer, recalling the definition of $\omega_\eps$ in \eqref{eq:def omega_eps} and the usual Modica-Mortola/Bogomol'nyi trick,  
		\begin{align*}
		\int_{\Omega\times (0,\infty)}\zeta(x)\eta(t)\d\omega_\eps(x,t)\leq \int_{0}^\infty\eta(t)E_\eps^S(u_\eps(\cdot,t))\dt\leq \int_{0}^\infty\eta(t)E_\eps(u_\eps(\cdot,t))\dt ,
		\end{align*}
		recalling $\zeta\in[0,1]$ and $E^P_\eps(u_\eps)\geq0$. By means of the compactness \eqref{eq:compactness omega} and dominated convergence first, then by monotone convergence to eliminate $\zeta$, we are led to  
		$$
		\int_{\Omega\times (0,\infty)}\eta(t)\d\omega(x,t)\leq \int_0^\infty \eta(t)E(t)\dt,
		$$
		which allows to conclude (see the proof of \cite[Lemma 2]{LH} for the details) that $\omega$ can be disintegrated as in the statement of the theorem, and also that \eqref{eq:comp} is satisfied. For $\mu$ a similar result can be shown, up to another disintegration in the $p$-variable, $p\in \mathbb S^{d-1}$. We are left to show \eqref{sharp}, where we make use of Remark \ref{additional}. To this aim, we first notice that, considering the same subsequence for which \eqref{energia} holds, we have, for any $T>0$, by Fatou's Lemma, 
		\begin{align}
			&\int_0^T\liminf_{\eps\tto0}E_\eps^P(u_\eps(\cdot,t))\dt\leq \liminf_{\eps\tto0}\int_0^TE_\eps^P(u_\eps(\cdot,t))\dt\nonumber
            \\&=\liminf_{\eps\tto0}\frac{\eps^\alpha}2\int_0^T\lambda_\eps(t)^2\dt\leq \liminf_{\eps\tto0}C(T)\eps^\alpha= 0,\label{ff}
		\end{align} 
		where in the last inequality we used the uniform bound \eqref{lambda_bound}, so that $\liminf_{\eps\tto0}E_\eps^P(u_\eps(\cdot,t))=0$ for almost any $t>0$. Now, recalling the second identity in \eqref{eq:discrepancy 2}, together with the definition of $\omega_\eps$ and its weak* convergence \eqref{eq:compactness omega}, as well as the definition of $\widetilde\omega_\eps$ in \eqref{tildeomega}, we deduce that, up to further subsequences, 
		\begin{align}
		\widetilde \omega_\eps\stackrel{\ast}{\rightharpoonup} \omega,\quad\text{weakly* as Radon measures in }\Omega\times(0,\infty),
	\label{conven1}
		\end{align} 
		and thus, recalling \eqref{more} and the disintegration of $\omega$, we infer $\widetilde \omega_t=\omega_t$ for almost any $t>0$, and
		\begin{align*}
		(\widetilde \omega_\eps)_t\stackrel{\ast}{\rightharpoonup}  \omega_t,\quad\text{weakly* as Radon measures in }\Omega,
		\end{align*}
		for almost any $t>0$, as $\eps\tto0$. As a consequence, we also infer, since $\Omega=\mathbb T^d$,
		\begin{align}
			E^S_\eps(u_\eps(\cdot,t))=(\widetilde\omega_\eps)_t(\Omega)\to \omega_t(\Omega),
			\label{cv1}
		\end{align}
		for almost any $t>0$, and thus, recalling \eqref{energia} and \eqref{ff},
		\begin{align*}
		&	E(t)=\lim_{\eps\tto0} E(u_\eps(\cdot,t))=\liminf_{\eps\tto0}(E^S_\eps(u_\eps(\cdot,t))+E^P_\eps(u_\eps(\cdot,t)))\\&=\liminf_{\eps\tto0}E^P_\eps(u_\eps(\cdot,t))+\lim_{\eps\tto0}E^S_\eps(u_\eps(\cdot,t))=\omega_t(\Omega),
		\end{align*} 
		for almost any $t>0$, giving \eqref{sharp}. In conclusion, the validity of \eqref{eq:compatibility2} is immediate, since for any $\xi\in C_c^\infty(\Omega;\R^d)$ we have the relation
		\begin{align*}	-\int_0^\infty\int_\Omega \psi_\eps\Div\xi\dx\dt	=\int_0^\infty\int_\Omega \xi\cdot \nabla \psi_\eps\dx\dt=\int_{\Omega\times(0,\infty)\times \mathbb S^{d-1}}\xi\cdot p\d\mu_\eps,
		\end{align*}
		and we can pass to the limit as $\eps\tto0$, by means of \eqref{eq:compactness psi} and \eqref{eq:compactness mu}, to obtain 
			\begin{align*}	-\int_0^\infty\int_\Omega \xhi\Div\xi\dx\dt	=\int_0^\infty\int_\Omega \xi\cdot \nabla \xhi\dt=\int_{\Omega\times(0,\infty)\times \mathbb S^{d-1}}\xi\cdot p\d\mu.
		\end{align*}  
		A further localization in time, recalling the slicing in time of $\mu$, gives the result \eqref{eq:compatibility2} and concludes the proof of the lemma.
	\end{proof}
We need now to pass to the limit in the velocity equation. To this aim, as first introduced in \cite{LH}, we aim at showing that the couple $(\psi_\eps,V_\eps)$ satisfies an approximate transport equation of the form 
\begin{align}
	\partial_t\psi_\eps+V_\eps\omega_\eps=\partial_t\psi_\eps+V_\eps\bn_\eps\cdot\nabla \psi_\eps	\approx 0,\label{transport1}
\end{align}	
where $\bn_\eps:=\frac{\nabla \psi_\eps}{\normmm{\nabla\psi_\eps}}=\frac{\nabla u_\eps}{\normmm{\nabla u_\eps}}$ is the approximate normal vector which transports the quantity $\psi_\eps$. In particular, as $\eps\downarrow0$, we rigorously verify \eqref{transport1} and we also obtain that $\partial_t\psi\ll\omega$, so that we can define a final normal velocity $V:=-\frac{\d\partial_t\psi}{\d\omega}$, which is $\omega$-measurable and such that, by construction,  
$$
\partial_t\psi+V\omega=0,
$$
in the sense of measures, which is a kind of transport equation. The only difference lies in the fact that the information on the support of $\omega$ is not enough to guarantee that it is the same as $\normmm{\nabla\psi}$, but we have only (see Remark \ref{compatib}) $\normmm{\nabla \psi}\leq \omega$ in the sense of measures. In conclusion, in the proposition we provide a sharp energy dissipation inequality between the dissipative term $\int \eps(\partial_t u_\eps)^2$ and the velocity $V$, which is essential to prove the validity of an energy inequality à la De Giorgi.    
\begin{proposition}
	\label{propvel}
Under the same assumpions of Propositions \ref{prop:compactness}-\ref{prop:takasao}, then, for any $\eps>0$,
\begin{align}
	\int_{\Omega\times(0,\infty)}V_\eps^2\d\omega_\eps\leq	\int_{\Omega\times(0,\infty)}\eps(\partial_tu_\eps)^2\dx\dt\leq  E_\eps(u_{\eps,0}),\label{firstbound}
\end{align}
and, for any test function $\xi\in C_c(\Omega\times (0,\infty))$,
\begin{align*}
	\lim_{\eps\tto0} \left(\int_{\Omega\times(0,\infty)}\xi V_\eps\d\omega_\eps+\int_{\Omega\times(0,\infty)}\xi\partial_t\psi_\eps\dx\dt\right)=0.	
\end{align*}
Furthermore, it holds $\partial_t\psi	\ll \omega$, and thus there exists an $\omega$-measurable function \RevY{$V=-\frac{\d\partial_t\psi}{\d\omega}:\Omega\times(0,\infty)\to \R$}, corresponding to the normal velocity of $\xhi$ in the sense of \eqref{eq:evolPhase}. Also, for almost every $T>0$ the following estimate holds:
\begin{align}
	&\liminf_{\eps\tto0} \frac12\int_0^T\int_\Omega V_\eps^2\d\omega_\eps\geq \frac12 \int_0^T\int_{\Omega\times(0,T)} V^2\d\omega,
\end{align}
which, together with \eqref{firstbound}, gives the sharp dissipation lower bound
\begin{align}	
\liminf_{\eps\tto0} \frac12\int_0^T\int_\Omega\eps(\partial_t u_\eps)^2\dx\dt\geq \frac12 \int_0^T\int_\Omega V^2\d\omega= \frac12 \int_0^T\int_\Omega V^2\d\omega_t\dt.
	\label{Degiorgi}
\end{align} 
\end{proposition}
\begin{proof}
	The proof can be carried out following word by word the proofs of Lemma 3 and Proposition 3 of \cite{LH}. The essential ingredient is the pointwise nonpositivity of the discrepancy measure, which corresponds to \eqref{eq:takasao2}, together with \eqref{eq:discrepancy 1}-\eqref{eq:discrepancy 2}.
\end{proof}

In the following last proposition, we give a sharp inequality between the gradient terms in the energy dissipation inequalities. Namely it (formally) gives a sharp lower bound for the term $\int\frac1\eps(\delta E_\eps(u_\eps))$ in terms of the corresponding mean curvature terms in the sharp interface limit. Most of the proof is already present in \cite{Ilmanen} (see also \cite[Proof of Proposition 4]{LH}), but here we need to crucially adapt it due to the presence of the Lagrange multipliers. \RevY{Similar} arguments can be found in \cite{Takasao, CT}. We have
\begin{proposition}
	\label{propcurvature}
Under the same assumpions of Propositions \ref{prop:compactness}-\ref{prop:takasao}, then, there exists an $\omega$-measurable vector field $\mathbf H:\Omega\times(0,\infty)\to \R^d$, $\mathbf H\in L^2_{loc}((0,\infty);L^2(\Omega;\d\omega_t;\mathbb R^d))$, which is the mean curvature vector of the oriented space-time varifold $\mu$ in the sense of \eqref{eq:weakCurvature}. Additionally, for any $T>0$, the following (sharp) lower bound holds
\begin{align}
&\nonumber\liminf_{\eps\tto0}\frac12\int_{0}^T\int_\Omega \frac1\eps\normmm{\bH_\eps-\lambda_\eps\sqrt{2W(u_\eps(\cdot,t))}\frac{\nabla \psi_\eps}{\normmm{\nabla \psi_\eps}}}^2\dx\dt\\&
\geq \frac 12\int_0^T\int_\Omega \normmm{\bH-\lambda\langle\pi_{x,t}\rangle}^2\d\omega=\frac 12\int_0^T\int_\Omega \normmm{\bH-\lambda\langle\pi_{x,t}\rangle}^2\d\omega_t\dt,
	\label{lowerbound}
\end{align}
where, for any $\eps>0$, the quantity $\mathbf H_\eps$ is defined in \eqref{eq:Heps}, and $\lambda$ is the limit of the Lagrange multiplier $\lambda_\eps$, in the sense of \eqref{compactness_lambda}.
	\end{proposition}
\begin{proof}
	Let us first notice that the quantity in the left-hand side of inequality \eqref{lowerbound} coincides with the second dissipative term in \eqref{eq:AC4}, since, for any $\eps>0$ and for almost any $(x,t)\in \Omega\times(0,\infty)$, 
	\begin{align}
		& \nonumber\frac1\eps\normmm{\bH_\eps-\lambda_\eps\sqrt{2W(u_\eps(\cdot,t))}\frac{\nabla \psi_\eps}{\normmm{\nabla \psi_\eps}}}^2\\&\nonumber= \frac1\eps \left( \eps \Delta u_\eps(\,\cdot\,,t) - \frac1\eps W'(u_\eps(\,\cdot\,,t))+\lambda_\varepsilon\sqrt{2W(u_\eps(\cdot,t))} \right)^2\frac{\nabla \psi_\eps(\cdot,t)}{\normmm{\nabla \psi_\eps(\cdot,t)}}\cdot \frac{\nabla \psi_\eps(\cdot,t)}{\normmm{\nabla \psi_\eps(\cdot,t)}} \\&= \frac1\eps \left( \eps \Delta u_\eps(\,\cdot\,,t) - \frac1\eps W'(u_\eps(\,\cdot\,,t))+\lambda_\varepsilon\sqrt{2W(u_\eps(\cdot,t))} \right)^2.\label{identity}
	\end{align}
	Then, let us fix $T\in(0,\infty]$ and $\zeta\in C_c^1(\Omega\times(0,T);\R^d)$. From the inequality $\frac12\normmm{a}^2\geq a\cdot b-\frac12\normmm{b}^2$, with $a=\eps \Delta u_\eps(\,\cdot\,,t) - \frac1\eps W'(u_\eps(\,\cdot\,,t))+\lambda_\varepsilon\sqrt{2W(u_\eps(\cdot,t))}$ and $b=\eps\xi\cdot\nabla u_\eps$, we infer from the identity \eqref{identity} that
		\begin{align}
		\nonumber&\frac12\int_{0}^\infty\int_\Omega \frac1\eps\normmm{\bH_\eps-\lambda_\eps\sqrt{2W(u_\eps(\cdot,t))}\frac{\nabla \psi_\eps}{\normmm{\nabla \psi_\eps}}}^2\dx\dt\\&
		=\frac1{2\eps}\int_{0}^\infty\int_\Omega  \left( \eps \Delta u_\eps - \frac1\eps W'(u_\eps)+\lambda_\varepsilon\sqrt{2W(u_\eps)} \right)^2\dx\dt\nonumber\\&
		\geq \underbrace{\frac1\eps\int_0^\infty\int_\Omega \left(\eps\Delta u_\eps-\frac1\eps W'(u_\eps)\right)(\eps\zeta\cdot \nabla u_\eps)\dx\dt}_{I_1}+\underbrace{\frac1\eps\int_0^\infty\int_\Omega \left(\lambda_\varepsilon\sqrt{2W(u_\eps)}\right)(\eps\zeta\cdot \nabla u_\eps)\dx\dt}_{I_2}\nonumber\\&
		\underbrace{-\frac1{2\eps}\int_0^\infty\int_\Omega\normmm{\eps\zeta\cdot\nabla u_\eps}^2\dx\dt}_{I_3}.
		\end{align}
		Now the last term $I_3$ can be easily handled as
		\begin{align*}
			-\frac1{2\eps}\int_0^\infty\int_\Omega\normmm{\eps\zeta\cdot\nabla u_\eps}^2\dx\dt\geq 	-\frac12\int_0^\infty\int_\Omega\normmm{\zeta}^2\eps\normmm{\nabla u_\eps}^2\dx\dt,
		\end{align*}
		and the right-hand side converges to $-\frac12 \int_{\Omega\times(0,\infty)}\normmm{\zeta}^2\d\omega$,
		recalling the (asymptotic) equipartition of the energy as well as \eqref{eq:compactness omega}. Then, concerning $I_1$, we can argue as in the proof of \cite[Proposition 4]{LH}, but we give here a sketch of the argument for the sake of the reader. In particular, we first notice that, after an integration by parts and nowadays classical computations,
		\begin{align}
			I_1=\int_0^\infty\int_\Omega \mathbf T_\eps:\nabla \zeta\dx\dt,
		\end{align} 
		where
		\begin{align*}
			\mathbf T_\eps:=\left(\frac\eps2\normmm{\nabla u_\eps}^2+\frac1\eps W(u_\eps)\right)Id-\eps\nabla u_\eps\otimes \nabla u_\eps.
		\end{align*}
		Then we can write 
		\begin{align*}
			\int_0^\infty\int_\Omega \mathbf T_\eps:\nabla \zeta\dx\dt&=\int_0^\infty\int_\Omega \left(\frac\eps2\normmm{\nabla u_\eps}^2+\frac1\eps W(u_\eps)\right)\Div\zeta\dx\dt \\&\quad-\int_0^\infty\int_\Omega \left(\frac{\nabla u_\eps}{\normmm{\nabla u_\eps}}\cdot \nabla \zeta\right)\frac{\nabla u_\eps}{\normmm{\nabla u_\eps}}\eps\normmm{\nabla u_\eps}^2\dx\dt.
		\end{align*}
	Now the right-hand side is easily handled, since, from again from the equipartition of the energy \eqref{eq:discrepancy 2} and \eqref{eq:compactness omega}, 
	\begin{align*}
		\int_0^\infty\int_\Omega \left(\frac\eps2\normmm{\nabla u_\eps}^2+\frac1\eps W(u_\eps)\right)\Div\zeta\dx\dt\to \int_{\Omega\times (0,\infty)}\Div \zeta d\omega,
	\end{align*}
	as $\eps\tto0$.	Concerning the last term, recalling $\frac{\nabla \psi_\eps}{\normmm{\nabla \psi_\eps}}=\frac{\nabla u_\eps}{\normmm{\nabla u_\eps}}$ on $\{\nabla \psi_\eps\not=0\}$, we can write 
	\begin{align*}
	&	-\int_0^\infty\int_\Omega \left(\frac{\nabla u_\eps}{\normmm{\nabla u_\eps}}\cdot \nabla \zeta\right)\frac{\nabla u_\eps}{\normmm{\nabla u_\eps}}\eps\normmm{\nabla u_\eps}^2\dx\dt=\\&-\int_0^\infty\int_{\Omega\times \mathbb S^{d-1}} \left(p\cdot \nabla \zeta\right)p\d\mu_\eps-\int_0^\infty\int_\Omega \left(\frac{\nabla u_\eps}{\normmm{\nabla u_\eps}}\cdot \nabla \zeta\right)\frac{\nabla u_\eps}{\normmm{\nabla u_\eps}}\left(\eps\normmm{\nabla u_\eps}^2-\normmm{\nabla\psi_\eps}\right)\dx\dt.
	\end{align*}
Then the first term in the right-hand side converges to $-\int_{\Omega\times(0,\infty)\times  \mathbb S^{d-1}}\left(p\cdot \nabla \zeta\right)p\d\mu$, by means of \eqref{eq:compactness mu}. On the other hand, the other term converges to zero by Cauchy-Schwarz inequality and \eqref{eq:discrepancy 2}, as 
\begin{align*}
&\normmm{\int_0^\infty\int_\Omega \left(\frac{\nabla u_\eps}{\normmm{\nabla u_\eps}}\cdot \nabla \zeta\right)\frac{\nabla u_\eps}{\normmm{\nabla u_\eps}}\left(\eps\normmm{\nabla u_\eps}^2-\normmm{\nabla\psi_\eps}\right)\dx\dt}\\&
\leq \int_0^\infty\int_\Omega \normmm{\nabla \zeta}\sqrt\eps\normmm{\nabla u_\eps}\left\vert\sqrt\eps\normmm{\nabla u_\eps}-\frac1{\sqrt \eps}\sqrt{2W(u_\eps)}\right\vert\dx\dt\\&
\leq \left(\int_0^\infty\int_\Omega \normmm{\nabla \zeta}\eps\normmm{\nabla u_\eps}^2\dx\dt\right)^\frac12\left(\int_0^\infty\int_\Omega\normmm{\nabla \zeta}\left(\sqrt \eps\normmm{\nabla u_\eps}-\frac1{\sqrt \eps}\sqrt{2W(u_\eps)}\right)^2\dx\dt\right)^\frac12\\&\leq C\left(\int_0^\infty\int_\Omega\normmm{\nabla \zeta}\left(\sqrt \eps\normmm{\nabla u_\eps}-\frac1{\sqrt \eps}\sqrt{2W(u_\eps)}\right)^2\dx\dt\right)^\frac12\to 0,
\end{align*}
	as $\eps\tto0$. We are left with the new term $I_2$. This term is nontrivial, since in the integral there is a highly nonlinear quantity, and the convergence of $\lambda_\eps$ is only weak.  To handle this, let us first observe that, after an integration by parts, 
	\begin{align*}
		&I_2=\frac1\eps\int_0^\infty\int_\Omega \left(\lambda_\varepsilon\sqrt{2W(u_\eps)}\right)(\eps\zeta\cdot \nabla u_\eps)\dx\dt\\&=\int_0^\infty\int_\Omega \lambda_\varepsilon\zeta \cdot\nabla \psi_\eps\dx\dt=-\int_0^\infty\int_\Omega \lambda_\varepsilon\psi_\eps\Div\zeta  \dx\dt.
	\end{align*}
	Now, since $\zeta$ has compact support, there exists $T_\zeta\in(0,\infty)$ ($T_\zeta=T$ in case $T\in(0,\infty)$), such that 
	$$
	\supp\zeta\subset (0,T_\zeta)\times K,
	$$ 
	with $K\subset \Omega$ compact. We can thus estimate, by Cauchy-Schwarz inequality, 
	\begin{align*}
		&\normmm{\int_0^\infty\int_\Omega \lambda_\varepsilon \psi_\eps\Div\zeta  \dx\dt-\int_0^\infty\int_\Omega \lambda\psi\Div\zeta  \dx\dt}\\&\leq \sup_{(x,t)\in \Omega\times(0,\infty)}\normmm{\Div \zeta(x,t)}\int_0^\infty\normmm{\lambda_\eps}\normmm{\int_\Omega({\psi_\eps-\psi})\dx}\dt+\normmm{\int_0^\infty\int_\Omega (\lambda-\lambda_\eps)\psi\Div\zeta\dx\dt}\\&
		\leq 
		C_\zeta\left(\int_0^{T_\zeta}\lambda_\eps(t)^2\dt\right)^\frac12\left(\int_{0}^{T_\zeta}\normmm{\int_\Omega{(\psi_\eps-\psi)}\dx}^2\dt\right)^\frac12+\normmm{\int_0^{T_\zeta}\int_\Omega (\lambda-\lambda_\eps)\psi\Div\zeta\dx\dt}\\&
		\leq \underbrace{C_\zeta(T_\zeta)\left(\int_{0}^{T_\zeta}\normmm{\int_\Omega{(\psi_\eps-\psi)}\dx}^2\dt\right)^\frac12}_{I_4} +\underbrace{\normmm{\int_0^{T_\zeta}(\lambda-\lambda_\eps)\int_\Omega \psi\Div\zeta\dx\dt}}_{I_5},
	\end{align*}
	where in the last inequality we have used \eqref{lambda_bound}.
	We now observe that, due to mass conservation \eqref{eq:compactness cont vol}, since $\psi= \xhi$, it holds 
	$$
	\int_\Omega\psi(x,t)\dx=\int_\Omega\psi(x,0)\dx,\quad \forall t\geq 0,
	$$
	and thus, recalling \eqref{controlbasic1},
	\begin{align*}
		I_4&\leq C_\zeta(T_\zeta)\left(\left(\int_0^{T_\zeta}\normmm{\int_\Omega{(\psi_\eps(x,t)-\psi_\eps(x,0))}\dx}^2\dt\right)^{\frac12}+\left(\int_0^{T_\zeta}\normmm{\int_\Omega({\psi_\eps(x,0)-\psi(x,0)})\dx}^2\dt\right)^{\frac12}\right)\\&
		\leq C_\zeta(T_\zeta)T_\zeta^\frac12\left(\sup_{t\in(0,T_\zeta)}\normmm{\int_\Omega{(\psi_\eps(x,t)-\psi_\eps(x,0))}\dx}+\normmm{\int_\Omega({\psi_\eps(x,0)-\psi(x,0)})\dx}\right)\\&
		\leq C_{\zeta}\eps^{\frac\alpha2}+C_\zeta\normmm{\int_\Omega({\psi_\eps(x,0)-\psi(x,0)})\dx},
	\end{align*}
	so that, resorting to \eqref{L1conv}, we can conclude that $I_4\to0$ as $\eps\tto0$. In conclusion, since $\int_\Omega\psi\Div \zeta\dx \in L^2(0,T_\zeta)$, we can exploit the weak convergence of $\lambda_\eps$ in $L^2(0,T_\zeta)$ given by \eqref{compactness_lambda} to infer that $I_5\to0$ as $\eps\tto0$. This allows to conclude that 
	\begin{align}
	-	\int_0^\infty\int_\Omega \lambda_\varepsilon \psi_\eps\Div\zeta  \dx\dt\to -\int_0^\infty\int_\Omega \lambda\psi\Div\zeta  \dx\dt,\quad \text{ as }\eps\tto0.
	\end{align}
	Now we need to rewrite the limit in order to make the varifols $\omega$ appear, since we want to have a similar structure as the curvature term. In particular, we exploit the compatibility \eqref{eq:compatibility2}, which gives us, after applying the \RevX{Gau\ss\ Theorem},
	\begin{align}
	&\nonumber	-\int_0^\infty\int_\Omega \lambda\psi\Div\zeta  \dx\dt= \int_0^\infty\int_\Omega \lambda\zeta \cdot  \d\nabla \chi\dt\\&
	=\int_0^\infty\int_{\Omega\times \mathbb S^{d-1}} \lambda\zeta \cdot  p\d\mu\dt=\int_0^\infty\int_\Omega \lambda\langle\pi_{x,t}\rangle\cdot \zeta\d\omega,
		\end{align}
		where we recall that we set $\langle\pi_{x,t}\rangle=\int_{\mathbb S^{d-1}}p\d\pi_{x,t}$ for almost any $(x,t)\in \Omega\times(0,\infty)$, with $\pi_{x,t}$ as in Lemma \ref{lemmaslicing}. 
		
		All in all we have shown that 
		\begin{align}
	&\nonumber\liminf_{\eps\tto0}	\frac12\int_{0}^\infty\int_\Omega \frac1\eps\normmm{\bH_\eps-\lambda_\eps\sqrt{2W(u_\eps(\cdot,t))}\frac{\nabla \psi_\eps}{\normmm{\nabla \psi_\eps}}}^2\dx\dt\\&\geq \int_{\Omega\times(0,\infty)\times \mathbb S^{d-1}}(I_d-p\otimes p): \nabla \zeta\d\mu+\int_ {\Omega\times (0,\infty) } \lambda\langle\pi_{x,t}\rangle\cdot \zeta\d\omega-\frac12 \int_{\Omega\times(0,\infty)}\normmm{\zeta}^2\d\omega,\label{ineq1}
		\end{align}
	which holds for any $\zeta\in C_c^1(\Omega\times(0,T))$ and any $T\in(0,\infty]$. Observe now that we can define $L\in (C_c^1(\Omega\times(0,\infty)))'$ as to be a bounded linear functional such that
	$$
	L(\zeta):=\int_{\Omega\times(0,\infty)\times \mathbb S^{d-1}}(I_d-p\otimes p): \nabla \zeta\d\mu,
	$$
	and this operator, for any $T>0$, is such that $L\in (C_c^1(\Omega\times(0,T)))'$.
	Since the left-hand side of \eqref{ineq1} is finite by the energy estimate \eqref{eq:AC2},  then for any $\zeta\in C_c^1(\Omega\times(0,T))$,  such that $\norm{\zeta}_{L^2((0,T);L^2(\Omega,\d\omega_t;\R^d))}\leq 1$, it holds by Cauchy-Schwarz inequality,
	\begin{align}
	L(\zeta)&\leq C_0-\int_ {\Omega\times (0,T) } \lambda\langle\pi_{x,t}\rangle\cdot \zeta\d\omega+\frac12 \int_{\Omega\times(0,T)}\normmm{\zeta}^2\d\omega\nonumber\\&\nonumber \leq C_0+\int_0^T\left(\normmm{\lambda}\int_ {\Omega} \normmm{\zeta}\d\omega_t\right)\dt+\frac12 \int_{\Omega\times(0,T)}\normmm{\zeta}^2\d\omega\nonumber
    \\&\leq C_1\left(1+\left(\int_0^{T} \normmm{\lambda}^2\dt\right)^\frac12\right)\leq C_2(T),
	\end{align}
	for some $C_0,C_1,C_2>0$, with $C_1$ depending on $T$. Here we have crucially used the fact that $\lambda\in L^2_{loc}(0,\infty)$ due to \eqref{compactness_lambda}. As a consequence, we can uniquely extend $L$ by density as to be a linear bounded functional on $L^2((0,T);L^2(\Omega,\d\omega_t;\R^d))$ and thus, by Riesz Representation Theorem for Hilbert spaces, there exists a unique $\omega\llcorner\Omega\times(0,T)$-measurable ${\bH}_T:\Omega\times(0,T)\to \R^d$, such that
	\begin{align}
		L(\zeta)=-\int_{\Omega\times(0,T)}{\mathbf H}_T\cdot \zeta\d\omega,\label{repr}
	\end{align}
	for any $\zeta\in L^2((0,T);L^2(\Omega,\d\omega_t;\R^d))$ and any $T>0$. Thanks to the uniqueness of the Riesz representative on any interval $(0,T)$, we can now exploit a continuation argument, recalling that $L$ is actually defined over functions on $\Omega\times(0,\infty)$, to construct a global vector field $\bH:\Omega\times(0,\infty)\to \R^d$ which is $\omega$-measurable and such that $\bH_{|(0,T)}=\bH_T$ $\omega$-almost everywhere on $\Omega\times(0,\infty)$, for any $T>0$. As a consequence, it holds from \eqref{repr}, valid for any $T>0$,
		\begin{align}
		L(\zeta)=-\int_{\Omega\times(0,\infty)}{\mathbf H}\cdot \zeta\d\omega,\label{repr2}
	\end{align}
	for any $\zeta\in C_c(\Omega\times(0,\infty))$, so that $\bH$ corresponds exactly to the generalized mean curvature vector associated to $\mu$ defined in \eqref{eq:weakCurvature}. Concerning its integrability, we come back to \eqref{ineq1} and assume $T<\infty$. Using the representation \eqref{repr2} we thus have
	\begin{align}
		&\nonumber\liminf_{\eps\tto0}	\frac12\int_{0}^\infty\int_\Omega \frac1\eps\normmm{\bH_\eps-\lambda_\eps\sqrt{2W(u_\eps(\cdot,t))}\frac{\nabla \psi_\eps}{\normmm{\nabla \psi_\eps}}}^2\dx\dt\\&\geq \int_{\Omega\times(0,\infty)}(-\mathbf H+\lambda\langle\pi_{x,t}\rangle)\cdot \zeta\d\omega-\frac12 \int_{\Omega\times(0,\infty)}\normmm{\zeta}^2\d\omega,\label{ineq2}
	\end{align}
and we conclude the proof of \eqref{lowerbound} by taking the supremum over $\zeta\in C_c(\Omega\times(0,\infty))$. To obtain the local integrability in time of $\mathbf H$ we can combine the estimate \eqref{lowerbound} together with $\lambda\in L^2_{loc}(0,\infty)$, \eqref{energybound}-\eqref{eq:AC2}, after and application of Cauchy-Schwarz inequality:
\begin{align*}
	&\int_0^T\int_\Omega \normmm{\mathbf H}^2\d\omega\leq 2\int_{\Omega\times(0,\infty)}\normmm{\mathbf H-\lambda\langle\pi_{x,t}\rangle}^2\d\omega+2\int_0^T\int_\Omega \lambda^2\d\omega\\&
	\leq 
	2\liminf_{\eps\tto0}	\frac12\int_{0}^T\int_\Omega \frac1\eps\normmm{\bH_\eps-\lambda_\eps\sqrt{2W(u_\eps(\cdot,t))}\frac{\nabla \psi_\eps}{\normmm{\nabla \psi_\eps}}}^2\dx\dt+2\int_0^T\liminf_{\eps\downarrow0}E(u_\eps(\cdot,t))\lambda^2\dt\\&
	\leq (1+C(T))E_0,
\end{align*}
for any $T\in(0,\infty)$. This ends the proof of the proposition.
\end{proof}
We can finally put together all of these results above to conclude the proof of Theorem \ref{thm:existence}. 
\begin{proof}[Proof of Theorem \ref{thm:existence}]
The convergence results come directly from the previous propositions. In particular, convergence \eqref{energ} has been obtained in the proof of Proposition \ref{lemmaslicing} (see \eqref{conven1}). To check that the limit is a De Giorgi-type varifold solution, we check the items in Definition \ref{def:twoPhaseDeGiorgiVarifoldSolutions}. First, the velocity $V$ satisfies \eqref{eq:evolPhase} by Proposition \ref{propvel}, as well as property \eqref{eq:weakCurvature} is verified by Proposition \ref{propcurvature}. Then \eqref{squareint} and \eqref{conservationvol} are given in Proposition \ref{prop:compactness}. Also, the compatibility condition \eqref{eq:compatibility} is shown to hold in Lemma \ref{lemmaslicing}. Concerning the De Giorgi-type inequality \eqref{eq:DeGiorgiInequality}, we fix $T\in(0,\infty)$ and we aim at passing to the limit in 
  \begin{align}
	&E_\eps(u_\eps(\,\cdot\,,T))\nonumber\\& +\frac12 \int_0^T\int_{\Omega}  \eps(\partial_tu_\eps(\,\cdot\,,t))^2 \dx\dt+ \frac12\int_0^T\int_{\Omega}  \frac1\eps\normmm{\bH_\eps-\lambda_\eps\sqrt{2W(u_\eps(\cdot,t))}\frac{\nabla \psi_\eps}{\normmm{\nabla \psi_\eps}}}^2 \dx\dt\nonumber\\&\leq  E_\eps(u_\eps,0),
	\label{eq:ACfinal}
\end{align}
but this is now immediate, thanks to the bounds \eqref{eq:comp}, \eqref{Degiorgi}, \eqref{lowerbound} for the left-hand side, as well as \eqref{eq:Ilmanen initial 3} for the right-hand side.

The last property, i.e., the fact that the unoriented varifold  $\widehat\mu_t$, corresponding to $\mu_t$  (identifying antipodal points $\pm p$ on $\S^{d-1}$), is integer $(d-1)$-rectifiable for almost any $t\in(0,\infty)$, is already proven in \cite[Theorem 3.12]{Takasao} (see also \cite[Theorem 2.2]{Tonegawa} for the standard Allen-Cahn case).

In conclusion, by using this time \eqref{eq:AC4} we can pass to the limit in a similar way by means of  \eqref{eq:comp}, \eqref{Degiorgi}, \eqref{lowerbound} for the left-hand side (clearly the same estimates hold for the time interval $(s,t)$, $0\leq s<t$, in place of $(0,\infty)$) and using this time \eqref{sharp} for the right-hand side. After passing to the limit as $\eps\tto0$, we thus obtain the validity of \eqref{eq:DeGiorgiInequality} with $(0,T)$ replaced by $(s,t)$ for all $t>0$ and almost every $s<t$ \RevY{with} $s=0$ included. As a consequence, we immediately see that the total mass $t\mapsto \omega_t(\Omega)$ is nonincreasing. This concludes the proof of Theorem \ref{thm:existence}.
\end{proof}
		\section{Proof of Lemma \ref{lemma:further}}
           \label{sec:B}
In this section we show the consistency of our new notion of De Giorgi-type varifold solution.
\begin{itemize}
	\item [(i)] (Classical solutions are weak solutions). Most of the properties of Definition \ref{def:twoPhaseDeGiorgiVarifoldSolutions} can be shown to hold in a straightforward way from the regularity of the classical solution. The only property which might be less trivial is the De Giorgi-type sharp energy inequality \eqref{eq:DeGiorgiInequality}. Now, first notice that, since the evolution of $\mA$ is smooth (recall $\cI=\partial\mA$), it holds 
	\begin{align}
		\ddt \cH^{d-1}(\cI(t))=-\int_{\cI(t)} V_\cI\bn_\cI\cdot \bH_\cI,  \d\H^{d-1}
		\label{smoothevol}
	\end{align}  
	and, due to volume conservation, it holds 
	\begin{align}
		\int_{\cI(t)}V_\cI d\H^{d-1}=0,
		\label{conservvol}
	\end{align}
	so that, after an integration in time over $(0,T)$, we infer 
		\begin{align}
	 \cH^{d-1}(\cI(T))= \cH^{d-1}(\cI(0))-\int_0^T\int_{\cI(t)} V_\cI\bn_\cI\cdot (\bH_\cI-\lambda_\cI\bn_\cI)  \d\H^{d-1}\dt,
		\label{smoothevol1}
	\end{align}
	and thus, noticing that the velocity $V_\cI$ satisfies $V_\cI\bn_\cI= \bH_\cI-\lambda_\cI\bn_\cI$, one immediately gets inequality \eqref{eq:DeGiorgiInequality}.
	\item[(ii)](Smooth weak solutions are classical solutions). Let us consider $(\mu,\xhi)$ as a De~Giorgi-type varifold solution 
	such that there exists $T\in (0,\infty]$ so that $\chi(x,t)=\chi_{\mA(t)}(x)$
	and $\mu_t = |\nabla \chi(\cdot,t)| \otimes 
	\delta_{\frac{\nabla \chi(\cdot,t)}{|\nabla \chi(\cdot,t)|}}$ 
	for all $t\in [0,T]$ and some smoothly evolving family \RevX{$(\mA(t))_{t\in [0,T]}$}. Then most of the properties of a smooth mean curvature flow according to Definition \ref{strongsolmcf} are satisfied, and we only need to show \eqref{Vi}. \RevX{To this aim}, first recall that, thanks to the assumed characterization $\mu_t = |\nabla \chi(\cdot,t)| \otimes 
	\delta_{\frac{\nabla \chi(\cdot,t)}{|\nabla \chi(\cdot,t)|}}$, it holds $\bH=(\bH\cdot\bn)\bn$, where $\bn=\frac{\nabla \chi}{\normmm{\nabla \xhi}}$.
	Then, since the evolution of $\mA(t)$ is smooth, \eqref{smoothevol} holds, and, by the conservation of volume, also \eqref{conservvol} is satisfied, leading to 
			\begin{align}
		\cH^{d-1}(\partial\mA(T))= \cH^{d-1}(\partial\mA(0))-\int_0^T\int_{\partial\mA(t)} V\bn\cdot (\bH-\lambda\bn)  \d\H^{d-1}\dt.
		\label{smoothevol1bis}
	\end{align}
 Now, as a consequence of the assumptions on the solution, the De Giorgi inequality \eqref{eq:DeGiorgiInequality} reads
		\begin{align}
		\label{eq:DeGiorgiInequality2}
&	\cH^{d-1}(\partial\mA(T))+ \frac{1}{2} \int_{0}^{T} \int_{\partial\mA(t)} |V|^2 \d\H^{d-1} \dt
		+ \frac{1}{2} \int_{0}^{T} \int_{\partial\mA(t)} |\mathbf{H}-\lambda \bn|^2 \d\H^{d-1} \dt
		\nonumber\\&\leq  \cH^{d-1}(\partial\mA(0)),
	\end{align}
	and thus, completing the square, recalling $\normmm{\bn}=1$,
	\begin{align}
		\label{eq:DeGiorgiInequality2b}
		&	\nonumber\cH^{d-1}(\partial\mA(T))+ \frac{1}{2} \int_{0}^{T} \int_{\partial\mA(t)} |V\bn-\mathbf{H}+\lambda \bn|^2 \d\H^{d-1} \dt
		\\&+\int_0^T\int_{\partial\mA(t)} V\bn\cdot (\bH-\lambda\bn)  \d\H^{d-1}\dt
		\leq  \cH^{d-1}(\partial\mA(0)),
	\end{align}
	Subtracting \eqref{smoothevol1bis} to \eqref{eq:DeGiorgiInequality2b} we then get 
	\begin{align}
		\label{eq:DeGiorgiInequality3}
		&\int_{\partial\mA(t)} |V\bn-\mathbf{H}+\lambda \bn|^2 \d\H^{d-1} \dt\leq 0,
	\end{align}
	entailing the desired $V\bn=\mathbf{H}-\lambda \bn$ over $\partial\mA(t)$ for any $t\in(0,T)$. Clearly, recalling that \eqref{conservationvol} holds true, multiplying identity \eqref{eq:DeGiorgiInequality3} by $\bn$ and integrating over $\partial\mA(t)$, we also retrieve the explicit value of $\lambda$, which is, as expected,
	\begin{align*}
	\lambda(t)=\frac{\int_{\partial\mA(t)}\bH\cdot\bn\d\H^{d-1}}{\H^{d-1}(\partial\mA(t))},\quad\forall t\in[0,T),
	\end{align*} 
	concluding the proof of $(\mu,\xhi)$ as to be a classical smooth solution satisfying Definition \ref{strongsolmcf}. 
    \item[(iii)]  Most of the compactness proof follows identical as in \cite[Proof of Lemma 1, item (iv)]{LH}. Therefore, we only comment on the most relevant difference, which is the one concerning the mean curvature term. First, following \cite{LH}, we have that there exists $\omega:=\mathcal L^1\lcorner(0,\infty)\otimes (\omega_t)_{t\in(0,\infty)}$, with $(\omega_t)_{t\in(0,\infty)}$ as a weakly* $\mathcal L^1$-measurable family of Radon measures, such that, as $n\to\infty$, 
    \begin{align*}
        \omega_n\rightharpoonup \omega,\quad \text{ weakly* as Radon measures }\mathcal M(\Omega\times(0,\infty)),
    \end{align*}
    and, for almost any $t\geq0$,
    \begin{align*}
        \lim_{n\to\infty} (\omega_n)_t(\Omega)=\omega_t(\Omega).
    \end{align*}
    Moreover, again as in \cite{LH} (using the same arguments as in Proposition \ref{propvel}) we can prove the existence of $V\in L^2((0,\infty);L^2(\Omega,\d\omega_t))$ which is the Radon-Nykod\'{y}m derivative of $\partial_t\xhi$ with respect to $\omega_t$, so that, for almost any $T>0$,
    \begin{align*}
        \liminf_{n\to\infty}\int_{\Omega\times(0,T)} \normmm{V_n}^2\d\omega_n\geq \int_{\Omega\times(0,T)} \normmm{V}^2\d\omega,
    \end{align*}
    for almost any $T>0$. Furthermore, one can also show that, up to a subsequence, $\mu_n$ weakly* converges to some $\mu$ as Radon measures $\mathcal M(\Omega\times (0,\infty)\times\mathbb S^{d-1})$, with the corresponding mass measure as $\omega$. Note that to obtain this slicing result it is essential that the energy $t\mapsto(\omega_n)_t$ is nonincreasing in time for any $t$, as we need to apply \RevY{a similar poof as the one of} Lemma \ref{lemmaslicing}.
    
    Now, since $\int_0^T\lambda_n^2\dt\leq C(T)$ for almost any $T>0$, uniformly in $n$ by assumption, there exists $\lambda\in L^2_{loc}(0,\infty)$ such that
    \begin{align}
    \lambda_n\rightharpoonup \lambda,\quad \text{ weakly in }L^2(0,T),
        \label{convergenceweak}
    \end{align}
    as $n\to \infty$, for almost any $T>0$. Also, by a standard compactness argument (recall that $\normmm{\nabla \xhi_n(t)}\leq (\omega_n)_t$ in the sense of measures) there exists $\xhi\in L^\infty_{w*}((0,\infty); BV(\Omega))$ such that 
    \begin{align}
        \xhi_n(t)\to \xhi(t),\quad\text{strongly in }L^2(\Omega\times(0,T)),\label{chik}
    \end{align}
    for almost any $T>0$. Also, the compatibility condition \eqref{eq:compatibility} clearly holds in the limit as $n\to \infty$.
 We now show that 
    \begin{align}
    \int_0^T\int_\Omega \lambda_n\mediak\cdot \xi\d\omega_n\to \int_0^T\int_\Omega \lambda\media\cdot \xi\d\omega
        \label{intlambda}
    \end{align}
as $n\to\infty$ for any $\xi\in C_c^\infty(\Omega\times(0,T))$. Indeed, by the compatibility condition \eqref{eq:compatibility} and Gau\ss\ Theorem we have
\begin{align*}
   &\int_0^T\int_\Omega \lambda_n\mediak\cdot \xi\d\omega_n=\int_0^T\int_\Omega \lambda_n\bn_n \cdot \xi\d\normmm{\nabla \xhi_n(t)}\dt\\&=-\int_0^T\int_\Omega \lambda_n\chi_n \div\xi\dx\dt,
\end{align*}
and, analogously,
\begin{align*}
   &\int_0^T\int_\Omega \lambda\media\cdot \xi\d\omega=-\int_0^T\int_\Omega \lambda\chi \div\xi\dx\dt,
\end{align*}
and thus
\begin{align*}
    &\normmm{\int_0^T\int_\Omega \lambda_n\mediak\cdot \xi\d\omega-\int_0^T\int_\Omega \lambda\media\cdot \xi\d\omega}\\&
    \leq \normmm{\int_0^T\int_\Omega (\lambda_n-\lambda)\xhi\div \xi\dx\dt}+ \normmm{\int_0^T\int_\Omega \lambda_n(\xhi-\chi_n)\div\xi\dx\dt} \to0
\end{align*}
    as $n\to\infty$, since the first term in the right-hand side vanishes due to \eqref{convergenceweak}, whereas for the second one we have, by Cauchy-Schwarz inequality, 
\begin{align*}
 \normmm{\int_0^T\int_\Omega \lambda_n(\xhi-\chi_n)\div \xi\dx\dt} &\leq   \left({\int_0^T\lambda_n^2\dt}\right)^\frac12\left({\int_0^T\left(\int_\Omega\normmm{\xhi-\xhi_n}\dx\right)^2}\right)^\frac12\dt \\&  
\leq C(T)\left({\int_0^T\left(\int_\Omega\normmm{\xhi-\xhi_n}\dx\right)}\right)^\frac12\dt\to 0,
\end{align*}
as $n\to\infty$, thanks to \eqref{chik}. Thus convergence \eqref{intlambda} is proven. Finally we aim at proving 
    \begin{align}
\liminf_{n\to\infty}\frac12\int_{\Omega\times (0,T)}\normmm{\bH_n-\lambda_n\mediak}^2\d\omega\geq \frac12\int_{\Omega\times (0,T)}\normmm{\bH-\lambda\media}^2\d\omega.\label{ineqA}
    \end{align}
First, using \eqref{eq:weakCurvature}, we have
\begin{align}
& \nonumber\frac12\int_{\Omega\times (0,T)}\normmm{\bH_n-\lambda_n\mediak}^2\d\omega\\&\nonumber\geq  \int_{0}^{T} \int_{\Omega} (\mathbf{H}_n-\lambda_n\mediak) \cdot B \d(\omega_n)_t \dt -\frac12\int_{\Omega\times (0,T)}\normmm{B}^2\d(\omega_n)_t\dt\\&\nonumber
=- \int_{0}^{T} \int_{\Omega {\times} \S^{d-1}}
  				(I_d {-} p \otimes p) : \nabla B \d(\mu_n)_t dt-\int_0^T\int_\Omega \lambda_n\mediak\cdot B\d(\omega_n)_t\dt\\&\quad -\frac12\int_{\Omega\times (0,T)}\normmm{B}^2\d(\omega_n)_t\dt,\label{curvat}
\end{align}
for any $B \in C^\infty_{c}(\Omega {\times} (0,\infty);\R^d)$.
Then, thanks to \eqref{eq:DeGiorgiInequality} and the uniform bound on $\lambda_n$, we have the uniform control $\int_{0}^{T} \int_{\Omega} \normmm{\mathbf{H}_n-\lambda_n\mediak}^2 \d\omega_t\dt\leq C(T)$ for any $n\in \N$ and almost any $T>0$. Therefore, recalling \eqref{intlambda} and the weak* convergence of $\omega_n$ and $\mu_n$, we can pass to the limit as $n\to\infty$ and infer 
\begin{align*}
 C(T)&\geq  - \int_{0}^{T} \int_{\Omega {\times} \S^{d-1}}
  				(I_d {-} p \otimes p) : \nabla B \d\mu_t dt-\int_0^T\int_\Omega \lambda\media\cdot B\d\omega_t\dt\\&\quad -\frac12\int_{\Omega\times (0,T)}\normmm{B}^2\d\omega_t\dt,
\end{align*}
for some $C(T)\in(0,\infty)$ and for any $B \in C^\infty_{c}(\Omega {\times} (0,\infty);\R^d)$, for almost any $T>0$. As a consequence, by arguing exactly as in the proof of Proposition \ref{propcurvature} we can show that there exists $\bH\in L^2_{loc}((0,\infty);L^2(\Omega;\d\omega_t))$ such that 
            \begin{align}
  	\label{eq:weakCurvature2}
  				\int_{0}^{\infty} \int_{\Omega} \mathbf{H} \cdot B \d\omega_t \dt
  				= - \int_{0}^{\infty} \int_{\Omega {\times} \S^{d-1}}
  				(I_d {-} p \otimes p) : \nabla B \d\mu_t dt
  			\end{align}
  			for any $B \in C^\infty_{c}(\Omega {\times} (0,\infty);\R^d)$.
                In conclusion, taking first the $\liminf_{n\to\infty}$ in \eqref{curvat}, using \eqref{eq:weakCurvature2}, and then optimizing with respect to $B$ give \eqref{ineqA}.

                The proof is ended, since the above considerations show the validity of \eqref{eq:DeGiorgiInequality} in terms of $(\mu,\chi)$ as desired, for almost any $T>0$. All the other properties in Definition \ref{def:twoPhaseDeGiorgiVarifoldSolutions} are easily verified.
    \end{itemize}
\section{Proof of Theorem \ref{thmcalibration}. Existence of a gradient-flow calibration for volume preserving mean curvature flow.}
\label{sec:calibration}

Let $\cU_r(t) := \{x : \dist(x,\Sigmas(t)) < r \}$ denote the tubular neighborhood of $\Sigmas(t)=\partial\mA(t)$ with radius $r$ and let $P_{\cI(t)}(x)$ denote the nearest point projection onto $\Sigmas(t)$, whereas $s_{\cI}$ is the signed distance function defined as 
\begin{align}
s_{\cI(t)}(x,t)=\begin{cases}
	-\dist({x,\partial\Sigmas(t)}),\quad\text{ if }x\in \mA(t),\\
	\dist({x,\partial\Sigmas(t)}),\quad\text{ if }x\in \Omega\setminus \mA(t),
\end{cases}
\end{align}
By the assumed regularity on $\{\Sigmas(t)\}_{t\in[0,\Tstrong]}$, there exists  $\delta=\delta(\cI) > 0$ sufficiently small such that $P_{\cI(t)}$ is well defined and of class $C^1$ on $\cU_\delta(t)$, and injective for all $t \in [0, \Tstrong ]$ (see, e.g., \cite{AmbrosioDancer}). Notice that $\nabla \si(x,t)=-\bn_\cI(P_{\cI(t)}(x))$ for any $x\in \cU_\delta(t)$ and any $t\in[0,\Tstrong]$. Let us define also  a smooth cutoff function $ \zeta: \R \to [0,\infty)$, such that $\zeta (r) = 1 - (\frac{r}{\delta})^2$ for $r < \frac\delta2$, $\zeta = 0$ for $r > \delta$, and let us introduce $\xi(\cdot, t) := -\zeta (\si(\cdot, t))\nabla \si(\cdot, t)$, which then satisfies \eqref{xiprop}. Note that, as a consequence of the assumed regularity on the flow, it holds $ \xi\in C^{0,1}_c(\Omega\times[0,\Tstrong];\R^d)$, \extra{$\supp \xi(t)=\{x\in \Omega: \dist(x,\cI(t))\leq c\}$, with $c=\delta$}.

 Next, let $\vartheta$ be a smooth truncation of the identity function, say, an odd function (i.e., $\vartheta (r) = -\vartheta (-r)$), such that $\vartheta (r) = r$ for
$\abs{r} < \frac\delta2$ and $\vartheta (r) = \delta$ for $r\geq \delta$. We then set $\vartheta(x, t) := \vartheta (s_{\cI}(x, t))$, which corresponds to a smooth truncated version of the signed distance function. Additionally, we set $\lambdas(t):=\frac{\int_{\cI(t)}\div\xi\d\H^{d-1}}{\mathcal H^{d-1}(\cI(t))}$. 

We are left with the discussion of the approximated velocity vector field, which is the main novelty of the present paper.
 We can proceed in the construction of a vector $\bv=\bv(t)$ on $\mA(t)$, for any $t\in[0,\Tstrong]$ which will be then extended to construct the final vector $B$. Differently from \cite{Timvol}, here we do not consider the vector as to be the gradient of the solution of a Neumann-Laplace problem, but rather as the solution to an incompressible Stokes problem. This allows to impose a full Dirichlet boundary condition on $\bv(t)$, so that the final vector will point exactly in the normal direction on the interface, as desired. Let us then look for $\bv$ satisfying, for any $t\in[0,\Tstrong]$,
 \begin{align}
 	\begin{cases}
 		-\Delta\bv+\nabla \pi=0,&\quad\text{ in }\mA(t),\\
 		\Div \bv=0,&\quad\text{ in }\mA(t),\\
 		\bv=V_{\cI}\bn_{\cI},&\quad \text{ on } \cI(t).
 	\end{cases}\label{Stokes}
 \end{align}
 Since \RevY{$\mA(t)$ is a bounded domain (thus open and connected) for any $t\in[0,\Tstrong]$ by assumption, and since} $\cI(t)$ is of class $C^{2,\alpha}$, for $\alpha \in (0,1]$, there exists a unique solution $(\bv,\pi)$, with $\int_{\mA(t)}\pi\dx=0$, to the Stokes problem. Notice that this is possible since the compatibility condition 
 $$
 \int_{\cI(t)}(V_{\cI(t)}\bn_{\cI(t)})\cdot\bn_{\cI(t)}\d\H^{d-1}=0
 $$
 is satisfied due to the fact that $\cI(t)$ evolves under a volume-preserving mean curvature flow. 
 Then, the following Schauder-type estimates hold (see, for instance \cite{Agmon} or \cite[Lemma 5.3]{Dutto})
 \begin{align}
 	\nonumber&\norm{\bv(t)}_{C^{2,\alpha}(\overline{\mA(t)})}+	\norm{\pi}_{C^{1,\alpha}(\overline{\mA(t)})}\\&\leq C(\mA(t))\norm{V_{\cI(t)}(\cdot,t)\bn_{\cI(t)}(\cdot,t)}_{C^{2,\alpha}(\cI(t))}\leq C(\mA,\Tstrong),
 	\label{Schauder}
 \end{align}
 thanks to the regularity of the flow, which allows to control the constant uniformly in time. Now, since we have assumed that, for some $R_*>0$, $\mA(t)\subset B_{R_*}(0)$ for any $t\in[0,\Tstrong]$, a standard extension theorem (see \cite[Lemma 6.37]{Trudinger}), guarantees that there exists an extension $\overline\bv\in C^{2,\alpha}_c(B_{R_*}(0))$ such that 
   \begin{align}
 	&\norm{\overline\bv(t)}_{C^{2,\alpha}(\overline{B_{R_*}(0)})}\leq C(\mA,\Tstrong)\norm{\bv(t)}_{C^{2,\alpha}(\overline{\mA(t)})}\leq C(\mA,\Tstrong),
 	\label{Schauder3}
 \end{align}
 and $\overline\bv=\bv$ on $\mA(t)$, where again the constant can be bounded uniformly in time due to the regularity properties of $\{\mA(t)\}_{t\in[0,\Tstrong]}$. To define the final vector $B$, we consider a smooth bump function $h:\R\to\R$ such that $h\equiv1$ on $[-\delta/2,\delta/2]$, and $h\equiv 0$ on $(-\infty,-\delta)\cup (\delta,+\infty)$, and \RevY{we study the composition $h(\si)$}. Since \RevX{$\si(\cdot,t)\in C^{2,\alpha}(\cU_\delta(t))$}, it holds that $h(\si(\cdot,t))\in C^{2,\alpha}_c(\cU_\delta(t))$, for any $t\in[0,\Tstrong]$. Namely, since the evolution is smooth, we also have the uniform bound $\norm{h(\si(\cdot,t))}_{C^{2,\alpha}(\overline{\cU_\delta(t)})}\leq C(\mA,\Tstrong)$. Note that $ h(\si(\cdot,t))\equiv 1$ on $\cU_{\delta/2}(t)$ for all $t\in[0,\Tstrong]$. Then our vector field $B$ is set to be
\begin{align}
	B(t):=h(\si(\cdot,t))\overline\bv(\cdot,t)\in C^{1,1}_c(\cU_\delta(t)),
	\label{vectorB}
\end{align}
 satisfying, thanks to the fact that $\bv$ solves the Stokes problem \eqref{Stokes} and thanks to the construction of $h(\si)$,
 \begin{align}
 	(\Div B)(\cdot,t)=(\Div \overline\bv)(\cdot,t)=0,\quad\text{ in }{\cU_{\delta/2}(t)\cap \mA(t)},
 	\label{Divcondition}
 \end{align}
 clearly entailing also $\Div B=0$ on $\cI(t)$, by (Lipschitz) continuity. Furthermore, we have 
 \begin{align}
 	B(t)=V_{\cI(t)}\bn_{\cI(t)},\quad \text{ on }\cI(t),\label{normvelox}
 \end{align} 
 as well as the uniform in time control
 \begin{align}
 	&\norm{B(t)}_{C^{1,1}(\overline{\cU_{\delta}})}\leq C(\mA,\Tstrong),
 	\label{unifcontrol}
 \end{align}
 for any $t\in[0,\Tstrong]$, which comes from \eqref{Schauder3} together with the regularity of the flow and the fact that $h(\si(\cdot,t))\in C^{1,1}_c(\cU_\delta(t))$ uniformly in time. Note that, additionally, $\supp B(t)\subset \supp\xi(t)$ for any $t\in[0,\Tstrong]$, by construction.
 
 It is the right time to verify all the properties of the gradient-flow calibration. Properties \eqref{controltheta}, \eqref{identityn}, \eqref{short}, \eqref{lambdas}, and \eqref{varthetacontrol} are a straightforward consequence of the choice of $\xi,\vtheta,\lambdas$ above. Analogously, \eqref{regB}, \eqref{supports}, and \eqref{unifcontrol1} come from
 \eqref{vectorB} and \eqref{unifcontrol}. Additionally, property \eqref{div} can be inferred from \eqref{Divcondition}, entailing $\Div B=0$ on $\cI(t)$ and thus the result by Lipschitz continuity. Analogously, we can derive \eqref{normal}. Indeed, let us define the trivial extension $\overline B(x,t):=B(P_{\cI(t)}(x),t))$ over $\cU_\delta(t)$. Then, due to \eqref{normvelox}, $(Id-\nabla s_{\cI}\otimes \nabla s_{\cI})\overline B=(Id-\bn_\cI\circ P_\cI\otimes \bn_\cI\circ P_\cI)B\circ P_{\cI}=0$ for any $x\in \cU_\delta(t)$ and any $t\in[0,\Tstrong)$. This entails, by the Lipschitz continuity of $B$, using the uniform control \eqref{unifcontrol} and recalling $\abs{x-P_{\cI(t)}(x)}=\dist(x,\cI(t))$,  
 \begin{align*}
 \normmm{(Id-\nabla s_{\cI}\otimes \nabla s_{\cI})B}= \normmm{(Id-\nabla s_{\cI}\otimes \nabla s_{\cI})(B-\overline B)}\leq C(\Tstrong, \mA)\dist(x,\cI(t)),
 \end{align*} 
 for any $x\in \cU_\delta(t)$, giving \eqref{normal}. Property \eqref{thetabasic} can be obtained as in \cite{Timvol}. Indeed, thanks to \eqref{normvelox}, it is immediate to verify that 
 $$
 \partial_t\vtheta+(\overline B\cdot\nabla)\vtheta=0,\quad \forall x\in \cU_\delta(t),
 $$  
 so that 
  $$
 \partial_t\vtheta+(B\cdot\nabla)\vtheta=((B-\overline B)\cdot\nabla)\vtheta,\quad \forall x\in \cU_\delta(t).
 $$ 
 Therefore \eqref{thetabasic} directly comes from the (uniform in time) Lipschitz properties of $B$ and $\vtheta$, \RevX{ upon noticing that $\nabla \vartheta=0$ on $\Omega\setminus \cU_\delta$}. Similarly, we have 
 \begin{align*}
 	\pt\normmm{\xi}^2+(\overline B\cdot\nabla)\abs{\xi}^2=0
 \end{align*} 
 in $\cU_\delta(t)$, and thus again, arguing exactly as in \cite[Proof of Theorem 2.2]{Timvol},
  \begin{align*}
 \normmm{	\pt\normmm{\xi}^2+(B\cdot\nabla)\abs{\xi}^2}=\normmm{((B-\overline B)\cdot\nabla)\normmm{\xi}^2}\leq C\abs{s_\cI}^2,
 \end{align*} 
 giving \eqref{transport}. Furthermore, to obtain \eqref{transport}, we have 
 \begin{align*}
 	\RevX{\pt\xi+(B\cdot\nabla)\xi+(\nabla B)^T\xi=\pt\xi+(\overline B\cdot\nabla)\xi+(\nabla \overline B)^T\xi+((B-\overline B)\cdot\nabla)\xi+(\nabla (B-\overline B))^T\xi},
 \end{align*}
 and the first term in the right-hand side vanishes in $\cU_\delta(t)$, whereas the second one can be estimated as in \cite[Proof of Theorem 2.2]{Timvol}, giving in summary 
  \begin{align*}
 {\pt\xi+(B\cdot\nabla)\xi+(\nabla B)^T\xi}= f\xi+O({s_{\cI}}),
 \end{align*}
  where $f:=\xhi_{\cU_\delta(t)}\nabla s_\cI\cdot (\nabla \overline B)\nabla s_\cI\in L^\infty(\Omega\times(0,\Tstrong))$ thanks to the regularity of $s_\cI$ and \eqref{unifcontrol}. This entails \eqref{transport}. In conclusion, to obtain identity \eqref{geometric} note that, by the construction above, namely by properties \eqref{normvelox} and \eqref{identityn}, on $\cI(t)$ this identity becomes exactly \eqref{Vi}, which holds true as a direct consequence of the fact that $\mA(t)$ evolves according to a volume-preserving mean curvature flow. Therefore, \eqref{geometric} is obtained also on $\cU_\delta(t)$ by the (uniform in time) Lipschitz continuity of all functions involved. This concludes the proof of the existence of a gradient-flow calibration in the sense of Definition \ref{calibration}.
		\section{Weak-strong uniqueness: Proof of Theorem \ref{theo:weakStrongUniqueness}}
		
	\label{sec:C}
	In this section, we aim at proving the weak-strong uniqueness result of Theorem \ref{theo:weakStrongUniqueness}. To this aim, we will resort to the new gradient-flow calibration (see Theorem \ref{thmcalibration}) to approximate the smooth volume-preserving mean curvature flow. 
	 Let then $\Tstrong\in(0,\infty)$ be a finite time horizon, and let $(\mA(t))_{t\in [0,T_{strong}]}$ be a smoothly evolving family of bounded domains in $\Omega$ whose associated family of interfaces $(\mathcal I(t))_{t\in[0,T_{strong}]}$, $\mathcal I(t)=\partial \mA(t)$ is of class $C^{2,\alpha}$, for some $\alpha\in(0,1)$, for all $t\in[0,\Tstrong]$, assumed to evolve by a volume-preserving mean curvature flow as in Definition \ref{strongsolmcf}. As already discussed in the introduction and in Section \ref{sec:main}, the calibration proposed in \cite{Timvol} as well as the one introduced in \cite{Timkroemer} are not completely satisfying, in the sense that they construct an approximated velocity vector field $B$ which is not pointing in the normal direction on the interface $\cI(t)$. This is intrinsically assumed in \cite{Timkroemer}, since $B$ admits by construction a tangential vector field on the interface, whereas in \cite{Timvol}, where the main property of $B$ is to be divergence free in the domain $\mA(t)$, it is just a consequence of the choice of the construction of $B$ in $\mA(t)$. We have shown in Theorem \ref{thmcalibration} that changing this construction allows to retrieve the desired property while preserving the divergence-free condition of $B$ inside $\mA(t)$.
	With this new gradient-flow calibration for the strong solution at hand, we can build up some measures to quantify the distance between a De Giorgi-type varifold solution $(\xhi_{A(t)},\mu)$ according to Definition \ref{def:twoPhaseDeGiorgiVarifoldSolutions} (where as usual we set $\bn:=\frac{\nabla \xhi_{A(t)}}{\normmm{\nabla \chi_{A(t)}}}$ as to be the inward measure-theoretic normal to $\p^*A(t)$ for almost any $t>0$), and a strong solution for the volume-preserving mean curvature flow. Namely, as introduced in \eqref{eq:relativeEntropyIntro} and \eqref{eq:bulkErrorIntro}, we set the relative entropy functional 
	\begin{align}
	\label{eq:relativeEntropy1}
	\Er := \int_{\Omega {\times} \S^{d-1}} 
	1 - p \cdot \xi(x,t) \d\mu_t(x,p) \geq 0, 
	\quad t \in [0,\Tstrong],
\end{align}
as well as the bulk error
\begin{align}
	\label{eq:bulkError1}
	\Eb := \int_{\Omega} 
	\abs{\chi_{A(t)}(x) {-} \chi_{\mA(t)}(x)} \abs{\vartheta(x,t)}\dx \geq 0, 
	\quad t \in [0,\Tstrong].
\end{align}	
Note that it can be easily seen from the sign properties of $\vartheta$ that the bulk error can be rewritten also withour the moduli, namely 
\begin{align}
	\label{eq:bulkError2}
	\Eb = \int_{\Omega} 
	({\chi_{A(t)}(x) {-} \chi_{\mA(t)}(x)}) {\vartheta(x,t)}\dx \geq 0, 
	\quad t \in [0,\Tstrong].
\end{align}	

Exploiting the properties of $\xi$ and $\vartheta$ in Definition \ref{calibration}, we observe that the following fundamental coercivity properties hold, which are the same as the ones for the standard mean curvature flow (see, e.g., \cite[Section 4]{LH})
\begin{align}
	\label{coerciv1}
	\int_{\Omega {\times} \S^{d-1}} \abs{p {-} \xi(\cdot,t)}^2 \d\mu_t 
	&\leq 2\Er,
	\\ \label{coerciv2}
	\int_{\Omega} \min\{1,\RevX{\frac1{c^2}}\ \dist(\cdot,\mathcal{I}(t))^2\} \d\omega_t
	&\leq 2\Er,
	\\ \label{coerciv4}
	\int_{\Omega} \abs{\chi_{A(t)} {-} \chi_{\mA(t)}}
	\min\{1,\RevX{\frac1c}\ \dist(\cdot,\mathcal{I}(t))\} \dx
	&\leq \Eb.
\end{align}
Note that the first two properties are a consequence of the fact that $1-p\cdot \xi\geq \frac12\normmm{\xi-p}^2+\frac12(1-\normmm{\xi})$ and property \eqref{short}, whereas the last two are a simple consequence of the sign condition on the weight $\vartheta$ given in property (vi) of Definition \ref{calibration} and of property \eqref{varthetacontrol}, respectively.

We also point out that the relative entropy is a very flexible functional, as it admits equivalent expressions. To compute its time evolution, we rewrite it exploiting the compatibility condition~\eqref{eq:compatibility} for the De Giorgi-type varifold solution, together with an integration by parts, to obtain
\begin{align}
	\Er
	\label{eq:relEntropy2}
	= \int_{\Omega} 1 \d\omega_t + \int_{\Omega} \chi(\cdot,t)
	(\Div \xi)(\cdot,t) \dx.
\end{align}
In conclusion, to collect other coercivity properties, we add and subtract
the quantity $\rho=\frac{\d\normmm{\nabla \xhi_A(t)}}{\d\omega_t}$ defined in Remark \ref{compatib}, and exploit the compatibility \eqref{eq:compatibility} once more, to infer
\begin{align}
	\label{relEntropy3}
	\Er
	&= \int_{\Omega} 1 - \rho(\cdot,t) \d\omega_t 
	+ \int_{\p^*A(t)} 1 - \bn(\cdot,t)\cdot\xi(\cdot,t) \d\H^{d-1}.
\end{align}
As a consequence, we get the following basic further coercivity results:
\begin{align}
	\label{coerciv5}
	\int_{\Omega} 1 - \rho(\cdot,t) \d\omega_t
	&\leq \Er,\\
		\label{coerciv5b}
	\int_{\pA(t)} 1 - \bn(\cdot,t)\cdot \xi(\cdot,t)  \d\H^{d-1} 
	&\leq \Er,
	\\ \label{coerciv6}
	\int_{\p^*A(t)} \abs{\bn(\cdot,t){-}\xi(\cdot,t)}^2 \d\H^{d-1}
	&\leq 2\Er,
	\\ \label{coerciv7}
	\int_{\p^*A(t)} \min\{1,\RevX{\frac1{c^2}}\ \dist(\cdot,\mathcal{I}(t))^2\}\d\H^{d-1}\leq \int_{\RevX{{\p^*A(t)}}} (1-\normmm{\xi})\d\H^{d-1} 
	&\leq 2\Er.
\end{align}
Note that \eqref{coerciv5} comes from \eqref{relEntropy3} the fact that $\rho\in[0,1]$, and $\normmm{\bn\cdot \xi}\leq 1$ over $\partial^*A(t)$. Then, \eqref{coerciv6} is a consequence of the simple property that $\normmm{\xi}\leq 1$ and thus $ \abs{\bn(\cdot,t){-}\xi(\cdot,t)}^2\leq 2(1-\bn(\cdot,t)\cdot\xi(\cdot,t))$ on $\partial^*A(t)$. In conclusion, \eqref{coerciv7} is entailed by $1-\xi\cdot 	\bn\geq \frac12\normmm{\bn-\xi}^2+\frac12(1-\normmm{\xi})$ and \eqref{short}, similarly to \eqref{coerciv2}. \extra{In conclusion, we observe that, by the definition of $\rho\in[0,1]$ and the fact that $\xi\cdot\bn\leq 1$ over $\partial^*A(t)$, we have
\begin{align}
\nonumber&\int_{\Omega\setminus\supp\xi(t)} 1\d\omega_t\\&\leq \int_{\Omega\setminus\supp\xi(t)} 1\d\omega_t+\int_{\supp\xi(t)} (1-\rho(\cdot,t))\d\omega_t+\int_{\pA(t)\cap \supp\xi(t)}(1-\bn(\cdot,t)\cdot\xi(\cdot,t))\d\cH^{d-1}\nonumber\\&
= \Er.
\label{outsupp}
\end{align}
}
We can finally proceed with the proof of Theorem \ref{theo:weakStrongUniqueness}, which takes inspiration from \cite{LH}, but incorporates the difficulties related to the presence of the Lagrange multiplier. In this sense, we are peforming the lift of the weak-strong uniqueness result in \cite{Timvol} to the varifold setting.
In the following we set $T\in(0,\Tstrong]$.
First, from \eqref{eq:evolPhase}, which is satisfied by the De Giorgi-type varifold solution, by setting $\zeta\equiv 1$, we infer, recalling mass conservation \eqref{conservation}, \RevX{for almost any $T>0$},
\begin{align*}
	0=\int_\Omega\xhi(\cdot, T)\dx- \int_\Omega\xhi_0\dx=-\int_0^T\int_\Omega V\d\omega_t\dt.
\end{align*}
\RevX{This actually entails that 
\begin{align*}
	\int_0^T\int_\Omega V\d\omega_t\dt=0,\quad \forall T>0,
\end{align*}
since $V\in L^2((0,\infty);L^2(\Omega,\d\omega_t))$ and thus the function $t\mapsto \int_0^t\int_\Omega V(x,s)\d\omega_s\d s$ is absolutely continuous on $[0,\infty)$.} \RevX{Therefore, we can localize in time} and deduce
\begin{align}
	\int_\Omega V\d\omega_t=0,\quad \text{ for almost any }t>0,
\label{localized}
\end{align}
so that we also have, for $T\in(0,\Tstrong]$,
\begin{align}
	\int_0^T\int_\Omega \lambdas V\d\omega_t\dt=0.
	\label{localizedlambda}
\end{align}
Now, recalling the alternative definition of the relative entropy \eqref{relEntropy3}, we can exploit the De Giorgi-type inequality \eqref{eq:DeGiorgiInequality} and again \eqref{eq:evolPhase}, with $\zeta=\Div\xi$, to obtain, after an integration by parts, and also adding \eqref{localizedlambda}, 
\begin{align}
\nonumber\Erel{T}&\leq E[\mu_0,\xhi_0|\mA(0)]- \frac{1}{2} \int_{0}^{T} \int_{\Omega} |V|^2 \d\omega_t \dt
- \frac{1}{2} \int_{0}^{T} \int_{\Omega} \abs{\bH-\lambda\media}^2 \d\omega_t \dt
\\&
- \int_{0}^{T} \int_{\Omega} V(\Div\xi -\lambdas)\d\omega_t \dt
- \int_{0}^{T} \int_{\pA(t)} \bn\cdot\partial_t\xi \d\H^{d-1} \dt.	
	\label{estA}
\end{align}
We now need to estimate all these terms separately. First, we aim at making the approximate transport equations \eqref{transport}-\eqref{almosttranported} appear, and thus we obtain
\begin{align}
&\nonumber	- \int_{0}^{T} \int_{\pA(t)} \bn\cdot\partial_t\xi \d\H^{d-1} \dt\\&\nonumber=- \int_{0}^{T} \int_{\pA(t)} \bn\cdot(\partial_t\xi+(B\cdot\nabla )\xi+(\nabla B)^T\xi) \d\H^{d-1} \dt\\&\nonumber\quad+\int_{0}^{T} \int_{\pA(t)} \bn\cdot((B\cdot\nabla )\xi+(\nabla B)^T\xi) \d\H^{d-1} \dt\\&\nonumber
=-\int_{0}^{T} \int_{\pA(t)} (\bn-\xi)\cdot(\partial_t\xi+(B\cdot\nabla )\xi+(\nabla B)^T\xi) \d\H^{d-1} \dt\\\nonumber&\quad-\int_{0}^{T} \int_{\pA(t)} \xi\cdot(\partial_t\xi+(B\cdot\nabla )\xi) \d\H^{d-1} \dt\\&\quad\nonumber+\int_{0}^{T} \int_{\pA(t)} \xi\otimes(\bn-\xi):\nabla B \d\H^{d-1} \dt\\&\quad\nonumber+\int_{0}^{T} \int_{\pA(t)} \bn\cdot(B\cdot\nabla )\xi \d\H^{d-1} \dt\\&\nonumber
\leq \int_{0}^{T} \int_{\pA(t)\RevX{\cap \supp\xi(t)}} \left(({\bn-\xi})\cdot f{\xi}+\normmm{\bn-\xi}O(\dist(\cdot,\mathcal I(t))) \right)\d\H^{d-1} \dt\\&\quad\nonumber+
\int_{0}^{T} \int_{\pA(t)\RevX{\cap \supp\xi(t)}} O(\dist(\cdot,\mathcal I(t))^2)\d\H^{d-1} \dt\\&\quad\nonumber+\int_{0}^{T} \int_{\pA(t)} \xi\otimes(\bn-\xi):\nabla B \d\H^{d-1} \dt+\int_{0}^{T} \int_{\pA(t)} \bn\cdot(B\cdot\nabla )\xi \d\H^{d-1} \dt,\\&\nonumber
{\leq}  \int_{0}^{T} \int_{\pA(t)} \norm{f}_{L^\infty(\Omega\times(0,\Tstrong
	))}((1-\normmm{\xi}^2)+(1-\bn\cdot \xi))\d\H^{d-1} \dt\\&\quad\nonumber+C\int_{0}^{T} \int_{\pA(t)\RevX{\cap \supp\xi(t)}}O(\dist(\cdot,I(t))^2) \d\H^{d-1} \dt\\&\quad\nonumber
	+\int_{0}^{T} \int_{\pA(t)} \xi\otimes(\bn-\xi):\nabla B \d\H^{d-1} \dt+\int_{0}^{T} \int_{\pA(t)} \bn\cdot(B\cdot\nabla )\xi \d\H^{d-1} \dt\\&
	\leq \int_{0}^{T} \int_{\pA(t)} \xi\otimes(\bn-\xi):\nabla B \d\H^{d-1} \dt\nonumber\\&\quad +\int_{0}^{T} \int_{\pA(t)} \bn\cdot(B\cdot\nabla )\xi \d\H^{d-1} \dt+C\int_0^T\Erel{t}\dt,\label{last}
\end{align}
where we used \eqref{transport}-\eqref{almosttranported}, and, to deal with the extra term $f$ apprearing in \eqref{almosttranported}, we made use of the inequality $(\bn-\xi)\cdot \xi= (1-\normmm{\xi}^2)+(\bn\cdot \xi-1)\leq (1-\normmm{\xi}^2)+(1-\bn\cdot \xi)$, as well as the inequality 
\begin{align}
1-\normmm{\xi}^2\leq 2(1-\normmm{\xi})\leq 2(1-\bn\cdot\xi)\label{ctrlxi},
\end{align}
 which comes from $\normmm{\xi}\leq 1$. The coercivity estimates \eqref{coerciv5b},\eqref{coerciv7}, as well as the length control on $\xi$ given in \eqref{short} then allowed to obtain the final inequality in \eqref{last}. \RevX{Note that, since by construction it holds $\supp \xi(t)=\{x\in \Omega:\ \dist(x,\cI(t))\leq c\}$, we made use of the fact that 
 \begin{align*}
     \int_{0}^{T} \int_{\pA(t)\RevX{\cap \supp\xi(t)}}O(\dist(\cdot,I(t))^2) \d\H^{d-1} \dt\leq c^2C\int_{0}^{T} \int_{\pA(t)}\min\{1,\frac1{c^2}\dist(\cdot,I(t))^2\} \d\H^{d-1} \dt. 
 \end{align*}
 }
Then, concerning the first term in the right-hand side of \eqref{last}, by means of the compatibility condition \eqref{eq:compatibility}, and then exploiting the identity \eqref{eq:RadonNikodymProperty}, 
\begin{align}
	&\nonumber \int_{0}^{T} \int_{\pA(t)} \xi\otimes(\bn-\xi):\nabla B \d\H^{d-1} \dt\\&
	=\int_{0}^{T} \int_{\Omega\times \mathbb S^{d-1}} \xi\otimes(p-\xi):\nabla B \d\mu_t\dt +\int_{0}^{T} \int_{\Omega} (1-\rho)\xi\otimes\xi:\nabla B \d\omega_t \dt\nonumber\\&
	\stackrel{\eqref{coerciv5}} {\leq}\int_{0}^{T} \int_{\pA(t)} \xi\otimes(p-\xi):\nabla B \d\mu_t\dt +C\int_0^T\Erel{t}\dt.
	\label{ctrl}
\end{align}
About the second term in the right-hand side of \eqref{last}, we can write, recalling again  \eqref{eq:compatibility} and the identity \eqref{eq:RadonNikodymProperty}, 
\begin{align}
&\nonumber\int_{0}^{T} \int_{\pA(t)} \bn\cdot(B\cdot\nabla )\xi \d\H^{d-1} \dt\\&\nonumber=\int_{0}^{T} \int_{\pA(t)} \bn\cdot\Div(\xi\otimes B ) \d\H^{d-1} \dt-\int_{0}^{T} \int_{\pA(t)} \bn\cdot\xi\Div B \d\H^{d-1} \dt\\&\nonumber
=\int_{0}^{T} \int_{\pA(t)} \bn\cdot\Div(\xi\otimes B ) \d\H^{d-1} \dt-\int_{0}^{T} \int_{\pA(t)} (\bn\cdot\xi-1)\Div B \d\H^{d-1} \dt\\&
\quad\nonumber-\int_{0}^{T} \int_{\Omega} (\rho-1)\Div B \d\omega_t \dt\\&
\quad-\nonumber \int_{0}^{T} \int_{\Omega\times \mathbb S^{d-1}}(Id-p\otimes p):\nabla B \d\mu_t \dt\\&
\quad  - \int_{0}^{T} \int_{\Omega\times \mathbb S^{d-1}}(p\otimes p):\nabla B \d\mu_t \dt\nonumber
\\&\stackrel{\text{\eqref{coerciv5}-\eqref{coerciv5b}}}{\leq} \int_{0}^{T} \int_{\pA(t)} \bn\cdot\Div(\xi\otimes B ) \d\H^{d-1} \dt-\nonumber \int_{0}^{T} \int_{\Omega\times \mathbb S^{d-1}}(Id-p\otimes p):\nabla B \d\mu_t \dt\\&
\quad  - \int_{0}^{T} \int_{\Omega\times \mathbb S^{d-1}}(p\otimes p):\nabla B \d\mu_t \dt+C\int_0^T\Erel{t}\dt.\label{ptensorp}
\end{align}
\RevX{Here we have used the fact that, thanks to the disintegration $\mu_t=\omega_t\otimes (\pi_{x,t})_{x\in \Omega}$, and since $B$ does not depend on $p\in \mathbb S^{d-1}$, we have
\begin{align*}
    &\int_{\Omega\times \mathbb S^{d-1}}Id:\nabla B(x,t)\d\mu_t(x,p)\\&=\int_\Omega Id:\nabla B(x,t)\left(\int_{\mathbb S^{d-1}}1\d\pi_{x,t}(p)\right)\d\omega_t(x)=\int_\Omega Id:\nabla B(x,t)\d\omega_t(x),
\end{align*}
where the last identity is valid since $\pi_{x,t}$ is a probability measure on $\mathbb S^{d-1}$.
}
Now we can use the symmetric relation $\Div(\Div(\xi\otimes B))=\Div(\Div (B\otimes \xi))$ to obtain the identity, after integration by parts, $\bn\cdot  \Div (B\otimes \xi)=\bn\cdot  \Div (\xi\otimes B)=(\bn \cdot B)\Div \xi+\bn\otimes \xi:\nabla B$, so that, using again the compatibility \eqref{eq:compatibility}, we infer 
\begin{align*}
	&\int_{0}^{T} \int_{\pA(t)} \bn\cdot\Div(\xi\otimes B ) \d\H^{d-1} \dt\\&=\int_{0}^{T} \int_{\pA(t)} (\bn \cdot B)\Div \xi \d\H^{d-1} \dt+\int_{0}^{T} \int_{\pA(t)} \bn\otimes \xi:\nabla B\d\H^{d-1} \dt\\&
=	\int_{0}^{T} \int_{\pA(t)} (\bn \cdot B)\Div \xi \d\H^{d-1} \dt+\int_{0}^{T} \int_{\Omega\times \bS} p\otimes \xi:\nabla B\d\mu_t\dt,
\end{align*}
so that, putting together this last identity with estimates \eqref{last}, \eqref{ctrl}, and \eqref{ptensorp}, we end up with 
\begin{align}
		&- \int_{0}^{T} \int_{\pA(t)} \bn\cdot\partial_t\xi \d\H^{d-1} \dt\nonumber\\&
		\leq\nonumber \int_{0}^{T} \int_{\pA(t)} (\bn \cdot B)\Div \xi \d\H^{d-1} \dt-\int_{0}^{T} \int_{\Omega\times \mathbb S^{d-1}}(Id-p\otimes p):\nabla B \d\mu_t \dt\\&\quad -\int_{0}^{T} \int_{\Omega\times \bS} (p-\xi)\otimes (p-\xi):\nabla B\d\mu_t\dt+C\int_0^T\Erel{t}\dt.
	\label{final1}
\end{align}
Exploiting then \eqref{coerciv1}, we can deduce \RevX{that \eqref{estA}} can be upgraded to 
\begin{align}
	\nonumber\Erel{T}&\leq E[\mu_0,\xhi_0|\mA(0)]- \frac{1}{2} \int_{0}^{T} \int_{\Omega} |V|^2 \d\omega_t \dt
	- \frac{1}{2} \int_{0}^{T} \int_{\Omega} \abs{\bH-\lambda\media}^2 \d\omega_t \dt
	\\&\nonumber
	\quad+\int_{0}^{T} \int_{\pA(t)} (\bn \cdot B)\Div \xi \d\H^{d-1} \dt-\int_{0}^{T} \int_{\Omega\times \mathbb S^{d-1}}(Id-p\otimes p):\nabla B \d\mu_t \dt\\&
	\quad - \int_{0}^{T} \int_{\Omega} V(\Div\xi -\lambdas)\d\omega_t \dt+C\int_0^T\Erel{t}\dt.
	\label{estB}
\end{align}
We can now pass to discuss the controls given by the dissipative terms in the De Giorgi inequality \eqref{eq:DeGiorgiInequality}. We start from the velocity one, where the presence of the Lagrange multiplier $\lambdas$ has to be carefully taken into account. Namely, from \eqref{estB}, we have, recalling the compatibility condition \eqref{eq:compatibility}, the identity \eqref{eq:RadonNikodymProperty}, together with \eqref{geometric},

\begin{align}
&\nonumber- \frac{1}{2} \int_{0}^{T} \int_{\Omega} |V|^2 \d\omega_t \dt
	- \int_{0}^{T} \int_{\Omega} V(\Div\xi -\lambdas)\d\omega_t \dt+\int_{0}^{T} \int_{\pA(t)} (\bn \cdot B)\Div \xi \d\H^{d-1} \dt\nonumber\\&\nonumber
	=- \frac{1}{2} \int_{0}^{T} \int_{\Omega} |V+\Div \xi-\lambdas|^2 \d\omega_t \dt+\frac12 \int_{0}^{T} \int_{\RevX{\supp\xi(t)}} (\Div\xi-\lambdas)^2 \d\omega_t \dt\\&\quad +\frac12 \int_{0}^{T} \int_{\RevX{\Omega\setminus\supp\xi(t)}}\normmm{\lambdas}^2 \d\omega_t \dt\nonumber+\int_{0}^{T} \int_{\pA(t)} (\bn \cdot B)\Div \xi \d\H^{d-1} \dt\\&\nonumber=
	- \frac{1}{2} \int_{0}^{T} \int_{\Omega} |V+\Div \xi-\lambdas|^2 \d\omega_t \dt+\frac12 \int_{0}^{T} \int_{\RevX{\supp \xi(t)}} (\Div\xi-\lambdas+B\cdot\xi)^2 \d\omega_t \dt\\&\nonumber\quad+\frac12 \int_{0}^{T}\normmm{\lambdas}^2 \int_{\RevX{\Omega\setminus\supp\xi(t)}}1 \d\omega_t \dt- \int_{0}^{T} \int_{\Omega} (\Div\xi-\lambdas)B\cdot \xi \d\omega_t \dt\\&\nonumber\quad +\int_{0}^{T} \int_{\pA(t)} (\bn \cdot B)\Div \xi \d\H^{d-1} \dt-\frac12 \int_{0}^{T} \int_{\Omega} (B\cdot \xi)^2 \d\omega_t \dt\\&\nonumber
	=	- \frac{1}{2} \int_{0}^{T} \int_{\Omega} |V+\Div \xi-\lambdas|^2 \d\omega_t \dt+\frac12 \int_{0}^{T} \int_{\RevX{\supp \xi(t)}} (\Div\xi-\lambdas+B\cdot\xi)^2 \d\omega_t \dt\\&\nonumber\quad +\frac12 \int_{0}^{T}\normmm{\lambdas}^2 \int_{\RevX{\Omega\setminus\supp\xi(t)}}1 \d\omega_t \dt+\int_{0}^{T} \int_{\Omega} (\rho-1)(\Div\xi) B\cdot \xi \d\omega_t \dt\\&\quad +\int_0^T\int_{\pA(t)}(\Div\xi) B\cdot (\bn-\xi) \d\H^{d-1} \nonumber\dt\nonumber+\int_{0}^{T} \int_{\Omega} \lambdas B\cdot \xi \d\omega_t \dt-\frac12 \int_{0}^{T} \int_{\Omega} (B\cdot \xi)^2 \d\omega_t \dt\\&
		\stackrel{\text{\eqref{coerciv2}-\eqref{coerciv5}, \eqref{outsupp}}}{\leq} 	- \frac{1}{2} \int_{0}^{T} \int_{\Omega} |V+\Div \xi-\lambdas|^2 \d\omega_t \dt +\int_0^T\int_{\pA(t)}(\Div\xi) B\cdot (\bn-\xi) \d\H^{d-1} \nonumber\dt\\&\quad+\int_{0}^{T} \int_{\Omega} \lambdas B\cdot \xi \d\omega_t \dt-\frac12 \int_{0}^{T} \int_{\Omega} (B\cdot \xi)^2 \d\omega_t \dt+C\int_0^T(1+\RevX{\normmm{\lambda^*}^2})\Erel{t}\dt.
			\label{fundamental1}
\end{align}
 It is the appropriate time to concentrate on the mean curvature term in the inequality \eqref{estB}. Note that it is the only occasion in all the proof where we need to carefully exploit the novel gradient-flow calibration proposed in Theorem \ref{calibration}, which guarantees the fundamental new control \eqref{normal}, corresponding to the fact that $B$ points in the normal direction on $\cI(t)$. We then have, by \eqref{supports} and \eqref{xiprop},
 \begin{align}
   & \nonumber - \frac{1}{2} \int_{0}^{T} \int_{\Omega} \abs{\bH-\lambda\media}^2 \d\omega_t \dt-\frac12 \int_{0}^{T} \int_{\Omega} (B\cdot \xi)^2 \d\omega_t \dt\\
   &\nonumber\leq - \frac{1}{2} \int_{0}^{T} \int_{\supp B(t)} \abs{\bH-\lambda\media}^2 \d\omega_t \dt-\frac12 \int_{0}^{T} \int_{\Omega} (B\cdot \xi)^2 \d\omega_t \dt\\&
   =- \frac{1}{2} \int_{0}^{T} \nonumber\int_{\supp B(t)} \abs{(Id-\nabla s_\cI\otimes \nabla s_\cI)(\bH-\lambda\media)}^2 \d\omega_t \dt\\&\nonumber\quad-\frac12\int_0^T\int_{\supp B(t)}\normmm{\nabla s_\cI\cdot (\bH-\lambda\media)}^2(1-\normmm{\xi}^2)\d\omega_t\dt\\&\quad-\frac12\int_0^T\int_{\supp B(t)} \normmm{(\bH-\lambda\media)\cdot\xi-B\cdot \xi)}^2\d\omega_t\dt-\int_0^T\int_\Omega(B\cdot \xi)(\bH-\lambda\media)\cdot\xi\d\omega_t\dt.\label{ctrl1}
  \end{align}
Then, recalling \eqref{supports}, we obtain
\begin{align}
	&\nonumber -\int_{0}^{T} \int_{\Omega\times \mathbb S^{d-1}}(Id-p\otimes p):\nabla B \d\mu_t \dt\stackrel{\eqref{eq:weakCurvature}}{=}\int_0^T\int_\Omega \bH\cdot B\d\omega_t\dt \\&\nonumber=\int_0^T\int_{\supp B(t)} (\bH-\lambda\media)\cdot B\d\omega_t\dt+\int_0^T\int_\Omega \lambda\media\cdot B\d\omega_t\dt\\&\nonumber
    \stackrel{\eqref{xiprop}}{=}\int_0^T\int_{\supp B(t)}(Id-\nabla s_\cI\otimes \nabla s_\cI)(\bH-\lambda\media)):(Id-\nabla s_\cI\otimes \nabla s_\cI)B\d\omega_t\dt\\&\nonumber\quad+\int_0^T\int_\Omega(B\cdot \xi)(\bH-\lambda\media)\cdot\xi\d\omega_t\dt \\&\quad 
    +\int_0^T\int_{\supp B(t)}(1-\normmm{\xi}^2)(B\cdot \nabla s_\cI)(\bH-\lambda\media)\cdot\nabla s_\cI\d\omega_t\dt+\int_0^T\int_\Omega \lambda\media\cdot B\d\omega_t\dt.\label{ctrl2}
    \end{align}
Summing up \eqref{ctrl1} together with \eqref{ctrl2}, using Cauchy-Schwarz and Young's inequalities, we obtain
\begin{align}
    \nonumber &- \frac{1}{2} \int_{0}^{T} \int_{\Omega} \abs{\bH-\lambda\media}^2 \d\omega_t \dt-\frac12 \int_{0}^{T} \int_{\Omega} (B\cdot \xi)^2 \d\omega_t \dt-\int_{0}^{T} \int_{\Omega\times \mathbb S^{d-1}}(Id-p\otimes p):\nabla B \d\mu_t \dt\\&
    \leq\nonumber - \frac{1}{2} \int_{0}^{T} \nonumber\int_{\supp B(t)} \abs{(Id-\nabla s_\cI\otimes \nabla s_\cI)(\bH-\lambda\media)}^2 \d\omega_t \dt\\&\quad \nonumber-\frac12\int_0^T\int_{\supp B(t)}\normmm{\nabla s_\cI\cdot (\bH-\lambda\media)}^2(1-\normmm{\xi}^2)\d\omega_t\dt\\&\nonumber \quad-\frac12\int_0^T\int_{\supp B(t)} \normmm{(\bH-\lambda\media)\cdot\xi-B\cdot \xi)}^2\d\omega_t\dt\\&\nonumber
    \quad +\int_0^T\int_{\supp B(t)}(Id-\nabla s_\cI\otimes \nabla s_\cI)(\bH-\lambda\media)):(Id-\nabla s_\cI\otimes \nabla s_\cI)B\d\omega_t\dt \nonumber\\&\nonumber\quad 
    +\int_0^T\int_{\supp B(t)}(1-\normmm{\xi}^2)(B\cdot \nabla s_\cI)(\bH-\lambda\media)\cdot\nabla s_\cI\d\omega_t\dt+\int_0^T\int_\Omega \lambda\media\cdot B\d\omega_t\dt\nonumber
    \\&
    \leq - \frac{1}{4} \int_{0}^{T} \nonumber\int_{\supp B(t)} \abs{(Id-\nabla s_\cI\otimes \nabla s_\cI)(\bH-\lambda\media)}^2 \d\omega_t \dt\\&\quad \nonumber-\frac14\int_0^T\int_{\supp B(t)}\normmm{\nabla s_\cI\cdot (\bH-\lambda\media)}^2(1-\normmm{\xi}^2)\d\omega_t\dt\\&\nonumber \quad-\frac12\int_0^T\int_{\supp B(t)} \normmm{(\bH-\lambda\media)\cdot\xi-B\cdot \xi)}^2\d\omega_t\dt +C\int_0^T\int_{\Omega}\normmm{(Id-\nabla s_\cI\otimes \nabla s_\cI)B}^2\d\omega_t\dt \\&\quad 
    +\int_0^T\int_{\Omega}(1-\normmm{\xi}^2)(B\cdot \nabla s_\cI)^2\d\omega_t\dt+\int_0^T\int_\Omega \lambda\media\cdot B\d\omega_t\dt\nonumber\\&
\stackrel{\eqref{coerciv5b},\eqref{ctrlxi}}{\leq} \nonumber C\int_0^T\int_{\Omega}\normmm{(Id-\nabla s_\cI\otimes \nabla s_\cI)B}^2\d\omega_t\dt+\int_0^T\int_\Omega \lambda\media\cdot B\d\omega_t\dt+C\int_0^T\Erel{t}\dt\nonumber\\&
    \stackrel{\eqref{normal}}{\leq} \int_0^T\int_\Omega \lambda\media\cdot B\d\omega_t\dt+C\int_0^T\Erel{t}\dt.\label{f1}
\end{align}
Note that, additionally to \eqref{supports}, in the last step we crucially used property \eqref{normal}, which is the main novelty of the calibration given in Theorem \ref{calibration}, as well as the coercivity property \eqref{coerciv2}.
Then, coming back to \eqref{fundamental1}, we recall, by \eqref{ctrlxi}, that 
\begin{align}
\label{fondamentalexi}
\normmm{\xi\cdot(\bn-\xi)}\leq (1-\normmm{\xi}^2)+(1-\xi\cdot \bn)\leq 3(1-\bn\cdot \xi).
\end{align}
As a consequence, we can write
\begin{align}
&\nonumber\int_0^T\int_{\pA(t)}(\Div\xi) B\cdot (\bn-\xi) \d\H^{d-1} \nonumber\dt\\&\nonumber=\int_0^T\int_{\pA(t)}(\Div\xi) (\bn-\xi)(Id-\xi\otimes \xi)B \d\H^{d-1} \nonumber\dt\\&\quad +\int_0^T\int_{\pA(t)}(\Div\xi)(B\cdot \xi) (\bn-\xi)\cdot \xi \d\H^{d-1} \dt\nonumber\\&
\stackrel{\eqref{fondamentalexi}}{\leq} \int_0^T\int_{\pA(t)}(\Div\xi) (\bn-\xi)(Id-\xi\otimes \xi)B \d\H^{d-1} \dt+ C\int_0^T\Erel{t}\dt.
    \label{diff1}
\end{align}
We also observe that, by \eqref{xiprop} and \eqref{ctrlxi}, it holds 
\begin{align}
\normmm{\nabla s_\cI\otimes \nabla s_\cI-\xi\otimes \xi}\leq \sqrt{ 1+\normmm{\xi}^4-2\normmm{\xi}^2}=1-\normmm{\xi}^2\leq 2(1-\bn\cdot \xi),
    \label{ctrlxi1}
\end{align}
on $ \supp \xi(t)$, so that, by Cauchy-Schwarz and Young's inequalities, and again the crucial novel property \eqref{normal},
\begin{align}
    &\nonumber\int_0^T\int_{\pA(t)}(\Div\xi) (\bn-\xi)(Id-\xi\otimes \xi)B \d\H^{d-1} \dt\\&
    \leq \int_0^T\int_{\pA(t)}(\Div\xi) (\bn-\xi)(Id-\nabla s_\cI\otimes \nabla s_\cI)B \d\H^{d-1} \dt\nonumber\\&\nonumber
    \quad +\int_0^T\int_{\pA(t)}(\Div\xi) (\bn-\xi)(\nabla s_\cI\otimes \nabla s_\cI-\xi\otimes \xi)B \d\H^{d-1} \dt\\&\nonumber
   {\leq} \int_0^T\int_{\pA(t)}\normmm{\Div\xi} \normmm{\bn-\xi}^2 \d\H^{d-1} \dt\nonumber\\&\nonumber\quad +C\int_0^T\int_{\pA(t)}\normmm{\Div\xi} \normmm{(Id-\nabla s_\cI\otimes \nabla s_\cI)B }^2\d\H^{d-1} \dt\\&
    \quad +C\int_0^T\int_{\pA(t)}\normmm{\Div\xi} (1-\xi\cdot\bn) \d\H^{d-1} \dt\nonumber\\&
     \stackrel{\text{\eqref{coerciv5b}-\eqref{coerciv7},\eqref{normal}}}{\leq} C\int_0^T\Erel{t}\dt.\label{copia}
\end{align}
In conclusion, we deal with the Lagrange multiplier terms. First, from \eqref{fundamental1}, by means of \eqref{eq:compatibility},
\begin{align}
    &\nonumber\int_{0}^{T} \int_{\Omega} \lambdas B\cdot \xi \d\omega_t \dt\\&\nonumber
    = \int_{0}^{T} \int_{\Omega} \lambdas (1-\rho)B\cdot \xi \d\omega_t \dt+\int_{0}^{T} \int_{\pA(t)} \lambdas B\cdot (\xi-\bn) \d\H^{d-1} \dt+\int_{0}^{T} \int_{\pA(t)} \lambdas B\cdot \bn \d\H^{d-1} \dt\\&
    \nonumber=\int_{0}^{T} \int_{\Omega} \lambdas (1-\rho)B\cdot \xi \d\omega_t \dt+\int_{0}^{T} \int_{\pA(t)} \lambdas  (\xi-\bn)\cdot (Id-\nabla s_\cI\otimes \nabla s_\cI)B \d\H^{d-1} \dt\\&\nonumber\quad \RevX{+}\int_{0}^{T} \int_{\pA(t)} \lambdas  (\xi\cdot B)(\xi-\bn)\cdot\xi  \d\H^{d-1} \dt\\&\quad \RevX{-}\int_{0}^{T} \int_{\pA(t)} \lambdas  (\xi-\bn)\cdot (\xi\otimes \xi-\nabla s_\cI\otimes \nabla s_\cI)B \d\H^{d-1}
     +\int_{0}^{T} \int_{\pA(t)} \lambdas B\cdot \bn \d\H^{d-1} \dt\nonumber\\&
    \stackrel{\eqref{fondamentalexi},\eqref{ctrlxi1}}{\leq} C\int_0^T\normmm{\lambdas(t)}\Erel{t}\dt +\int_{0}^{T} \int_{\pA(t)} \lambdas B\cdot \bn \d\H^{d-1} \dt,\label{compat0}
\end{align}
 where in the last step we used the coercivity properties \eqref{coerciv5}, \eqref{coerciv5b}, \eqref{coerciv6}, as well as \eqref{normal}, after an application of Cauchy-Schwarz and Young's inequality (see also the estimate for \eqref{copia} above for a similar argument). In conclusion, concerning $\lambda$, we have from \eqref{f1}, thanks to the compatibility condition \eqref{eq:compatibility} (see, in particular, \eqref{equivalence}),
\begin{align}
\int_0^T\int_\Omega \lambda \media \cdot B\d\omega_t\dt= \int_0^T\int_{\pA(t)}\lambda  B\cdot \bn\d\H^{d-1}\dt
    \label{compat1}
\end{align}
We can now plug estimates \eqref{fundamental1}, \eqref{f1}, \eqref{diff1}, \eqref{copia}, \eqref{compat0} and \eqref{compat1} into \eqref{estB}, to obtain  

\begin{align}
	&\nonumber\Erel{T}\\&\nonumber\leq E[\mu_0,\xhi_0|\mA(0)]	- \frac{1}{2} \int_{0}^{T} \int_{\Omega} |V+\Div \xi-\lambdas|^2 \d\omega_t \dt\\&
   \quad +  \int_0^T(\lambdas(t)+\lambda(t))\int_{\pA(t)} B\cdot \bn\d\H^{d-1}\dt+C\int_0^T(1+\RevX{\normmm{\lambdas(t)}{{^2}}})\Erel{t}\dt.
	\label{estC}
\end{align}
We are left with the last term, related again to the Lagrange multipliers, which is the nonlocal one, coupling the relative entropy inequality with the bulk error inequality. Namely, first observe that by the Gau\ss\  Theorem, recalling the properties \eqref{identityn}, \eqref{geometric}, and then using the representation of $\lambdas$ given in \eqref{lambdas}, for any $t\in[0,T]$,
\begin{align}
&\nonumber\int_\Omega \xhi_{\mA(t)}\Div B\dx =-\int_{\cI(t)}B\cdot \bn_{\cI}\d\H^{d-1}\\&=-\int_{\cI(t)}B\cdot \xi\d\H^{d-1}=\int_{\cI(t)}(-\Div\xi+\lambdas)\d\H^{d-1}=0.\label{Lagrange}
\end{align}
Therefore, again Gau\ss\ Theorem gives
\begin{align}
	\nonumber&\int_{0}^{T} (\lambdas(t) +\lambda(t))\int_{\pA(t)} B\cdot \bn \d\H^{d-1} \dt\\&
	\nonumber =-\int_{0}^{T} (\lambdas(t) +\lambda(t)) \int_{\Omega}\xhi_{A(t)}\Div B \dx \dt\nonumber
	=-\int_{0}^{T} \int_{\Omega} (\lambdas(t) +\lambda(t))(\xhi_{A(t)}-\xhi_{\mA(t)})\Div B \dx \dt\\&\nonumber
	\leq  \int_{0}^{T} \int_{\Omega} (\normmm{\lambdas(t)} +\normmm{\lambda(t)})\normmm{\xhi_{A(t)}-\xhi_{\mA(t)}}\normmm{\Div B} \dx \dt\\&
\leq  \int_{0}^{T} (\normmm{\lambdas(t)} +\normmm{\lambda(t)})\int_{\RevX{\supp B(t)}} \normmm{\xhi_{A(t)}-\xhi_{\mA(t)}}O(\dist(\cdot,\cI(t))) \dx \dt\nonumber\\&\leq  \int_{0}^{T} (\normmm{\lambdas(t)} +\normmm{\lambda(t)}) \Ebulk{t} \dt,\label{estlambda}
\end{align}
where in the last step we used for the first time in the proof the fundamental approximately-zero divergence condition of $B$ given in \eqref{div}, as well as the coercivity condition \eqref{coerciv4}, \RevX{together with the fact that $\supp B(t)\subset \supp \xi(t)$ for any $t\in[0,T]$}. We stress again that, as anticipated when stating the main results, here we did not make use of any estimate on the closeness of the Lagrange multipliers $\lambda$ and $\lambdas$, which are only required to belong to $L^2(0,T)$ to apply Gronwall's Lemma later on.

To close estimate \eqref{estC} it is then imperative to find a control over the bulk error. To this aim, we use the approximate transport equation \eqref{thetabasic} for the weight $\vartheta$. \RevX{First, exploiting the evolution equation \eqref{eq:evolPhase} with $\zeta=\vartheta$, we obtain 
\begin{align}
 \int_{\Omega} 
	({\chi_{A(T)}\vartheta(T) {-} \chi_{A(0)}} {\vartheta(0)})\dx=\int_0^T\int_\Omega \xhi_{A}\pt\vartheta\dx\dt-\int_0^T\int_\Omega V\vartheta \d\omega_t\dt.\label{I1}  
\end{align}
Now observe that $\partial_t\xhi_{\mA(t)}$ is supported on $\cI(t)$ and  $\vartheta(\cdot,t)=0$ on $\cI(t)$, so that, integrating by parts in time, we have 
\begin{align}
   \nonumber\int_0^T\int_\Omega \xhi_{\mA(t)}\pt\vartheta\dx\dt&\nonumber=-\int_0^T\int_\Omega \partial_t\xhi_{\mA(t)}\vartheta(t)\dx\dt+\int_\Omega \xhi_{\mA(T)}\vartheta(T)\dx-\int_\Omega \xhi_{\mA(0)}\vartheta(0)\dx\\&= \int_\Omega \xhi_{\mA(T)}\vartheta(T)\dx-\int_\Omega \xhi_{\mA(0)}\vartheta(0)\dx.\label{I2}
\end{align}
Therefore, from the two identities \eqref{I1}-\eqref{I2} we obtain
\begin{align*}
 &\int_{\Omega} 
	({\chi_{A(T)}-\xhi_{\mA(T)})\vartheta(T) }\dx\\&=\int_\Omega(\xhi_{A(0)}-\xhi_{\mA(0)})\vartheta(0)\dx +\int_0^T\int_\Omega (\xhi_{A(t)}-\xhi_{\mA(t)})\pt\vartheta\dx\dt-\int_0^T\int_\Omega V\vartheta \d\omega_t\dt.  
\end{align*} 
As a consequence, considering the alternative formulation of the bulk error \eqref{eq:bulkError2} evaluated at times $T$ and $0$, we obtain 
}
\begin{align}
	\Ebulk{T}= E_{\mathrm{bulk}}[\xhi_0|\mA(0)]-	\int_0^T\int_\Omega V\vartheta \d\omega_t\dt +\int_{0}^{T}\int_{\Omega}({\xhi_{A(t)}-\xhi_{\mA(t)}})\pt\vartheta \dx \dt.
	\label{estE}
\end{align}
\RevY{Then, again recalling that $\vartheta(\cdot,t)=0$ on $\cI(t)$, together with the alternative definition  \eqref{eq:bulkError2} and by means of the Gau\ss\ Theorem, we infer} 
\begin{align}
&\int_{0}^{T}\int_{\Omega}({\xhi_{A(t)}-\xhi_{\mA(t)}})\pt\vartheta \dx \dt\nonumber\\&\stackrel{\eqref{thetabasic}}{=}\int_{0}^{T}\int_{\Omega}(\xhi_{A(t)}-\xhi_{\mA(t)})\RevX{O\left(\min\{1,\frac1{c}\dist(\cdot,I(t))\}\right) }\dx \dt-\int_{0}^{T}\int_{\Omega}(\xhi_{A(t)}-\xhi_{\mA(t)})B\cdot\nabla\vartheta \dx \dt\nonumber\\&
\stackrel{\eqref{coerciv4}}{\leq} C\int_{0}^{T} \Ebulk{t} \dt-\int_{0}^{T}\int_{\Omega}(\xhi_{A(t)}-\xhi_{\mA(t)})\Div(\vartheta B) \dx \dt\nonumber\\&\quad+\int_{0}^{T}\int_{\RevX{\supp B(t)}}(\xhi_{A(t)}-\xhi_{\mA(t)})\vartheta\Div B \dx \dt\nonumber
\\&\leq \int_{0}^{T}\int_{\pA(t)}\vartheta B\cdot\bn \d\H^{d-1} \dt+C\int_{0}^{T} \Ebulk{t} \dt\label{diff}.
\end{align} 
In conclusion, we need to treat the first term appearing in the right-hand side of the previous inequality. To this aim, the idea is to reconstruct the term $V+\Div \xi-\lambdas$, which appears squared with the negative sign in the estimate \eqref{estC}, as well as use the geometric evolution \eqref{geometric} once more. We thus write, using again the compatibility \eqref{eq:RadonNikodymProperty} and adding zero,
\begin{align}
&\nonumber -\int_0^T \int_\Omega V\vartheta \d\omega_t\dt+\int_{0}^{T}\int_{\pA(t)}\vartheta B\cdot\bn \d\H^{d-1} \dt\\&
\nonumber = -\int_0^T \int_\Omega V\vartheta \d\omega_t\dt+\int_{0}^{T}\int_{\pA(t)}\vartheta B\cdot(\bn-\xi) \d\H^{d-1} \dt\\&\nonumber
\quad +\int_0^T\int_{\pA(t)\RevX{\cap \supp \xi(t)}}\vartheta (B\cdot\xi+\div\xi-\lambdas)\d\H^{d-1}\dt+\int_0^T\int_{\RevX{\supp \xi(t)}}\vartheta (1-\rho) (\div\xi-\lambdas)\d\omega_t\dt\\&\nonumber\quad -\int_0^T\int_{\RevX{\supp \xi(t)}}\vartheta (\div\xi-\lambdas)\d\omega_t\dt\\&
=\underbrace{\int_{0}^{T}\int_{\pA(t)}\vartheta B\cdot(\bn-\xi) \d\H^{d-1} \dt}_{J_1} +\underbrace{\int_0^T\int_{\pA(t)\RevX{\cap \supp \xi(t)}}\vartheta (B\cdot\xi+\div\xi-\lambdas)\d\H^{d-1}\dt}_{J_2}\nonumber\\&\quad\underbrace{+\int_0^T\int_{\RevX{\supp \xi(t)}}\vartheta (1-\rho) (\div\xi-\lambdas)\d\omega_t\dt}_{J_3} \underbrace{-\int_0^T\int_{\RevX{\supp \xi(t)}}\vartheta (V+\div\xi-\lambdas)\d\omega_t\dt}_{J_4}.\label{complete}
	\end{align}
We estimate all these terms one by one. First, by Cauchy-Schwarz and Young's inequalities, using the Lipschitz property \eqref{varthetacontrol} of $\vartheta$ and the coercivity estimates \eqref{coerciv6} and \eqref{coerciv7},
\begin{align*}
	\normmm{J_1}\leq \int_{0}^{T}\int_{\pA(t)}\normmm{B}^2\vartheta^2 \d\H^{d-1} \dt+C\int_{0}^{T}\int_{\pA(t)}\normmm{\bn-\xi}^2 \d\H^{d-1} \dt \leq C\int_0^T\Erel\dt.
\end{align*}
Then, again by \eqref{varthetacontrol}, the geometric evolution equation \eqref{geometric}, and the coercivity \eqref{coerciv7}, it holds
\begin{align*}
	\normmm{J_2}\leq \int_{0}^{T}\int_{\pA(t)\RevX{\cap \supp \xi(t)}}O(\dist(\cI(\cdot,t))^2) \d\H^{d-1} \dt\leq C\int_0^T\Erel{t}\dt.
\end{align*}
\RevX{About $J_3$, we exploit the coercivity estimate  \eqref{coerciv5}, to get, recalling the regularity of $\xi$, 
\begin{align*}
	\normmm{J_3}\leq C\int_0^T(1+\normmm{\lambdas(t)})\int_{\Omega} (1-\rho)\d\omega_t\dt\leq  C\int_0^T(1+\normmm{\lambdas(t)})\Erel{t}\dt.
\end{align*}}
In conclusion, we have, Cauchy-Schwarz and Young's inequalities, together with the Lipschitz property \eqref{varthetacontrol} and the coercivty property \eqref{coerciv2},
\begin{align*}
	\normmm{J_4}&\leq \frac14 \int_0^T\int_\Omega \normmm{V+\div\xi-\lambdas}^2\d\omega_t\dt+C\int_0^T\int_{\RevX{\supp \xi(t)}} \normmm{\vartheta}^2\d\omega_t\dt\\&\leq \frac14 \int_0^T\int_\Omega \normmm{V+\div\xi-\lambdas}^2\d\omega_t\dt+C\int_0^T\Erel{t}\dt.
\end{align*}
To sum up, coming back to \eqref{estE}, we can use \eqref{diff}, \eqref{complete}, and the estimates over $J_i$, $i=1,\ldots,4$, to infer 
\begin{align}
		\nonumber&\Ebulk{T}\\&\leq E_{\mathrm{bulk}}[\xhi_0|\mA(0)]+\frac14 \int_0^T\int_\Omega \normmm{V+\div\xi-\lambdas}^2\d\omega_t\dt\nonumber\\&\quad +C\int_0^T(1+\normmm{\lambdas(t)})\Erel{t}\dt+C\int_0^T\Ebulk{t}\dt.
	\label{estF}
\end{align}
Summing this inequality with \eqref{estC}, and using also \eqref{estlambda}, we finally end up with 
\begin{align*}
	&\Erel{T}+\Ebulk{T}+\frac{1}{4} \int_{0}^{T} \int_{\Omega} |V+\Div \xi-\lambdas|^2\\&\nonumber\leq E[\mu_0,\xhi_0|\mA(0)]+E_{\mathrm{bulk}}[\xhi_0|\mA(0)]	\nonumber\\&\quad +C\int_0^T(1+\RevX{\normmm{\lambdas(t)}^2})\Erel{t}\dt+C\int_0^T(1+\normmm{\lambda(t)}+\normmm{\lambdas(t)})\Ebulk{t}\dt,
\end{align*}
so that, since $\lambdas,\lambda\in L^2(0,T)$, we can apply Gronwall's Lemma and conclude the proof of \eqref{eq:stabilityglobal}. Weak-strong uniqueness is then an immediate consequence of this result.
\appendix
\section{An alternative proof for weak-strong uniqueness for integral varifolds}
\label{App:extra}
In this Appendix we aim at proving the weak-strong uniqueness result of Theorem \ref{theo:weakStrongUniqueness} without making use of property \eqref{normal} of the gradient-flow calibration for the strong solution to volume-preserving mean curvature flow, whose existence is guaranteed by Theorem \ref{theo:weakStrongUniqueness}.  To this aim we make the further assumption on the weak (varifold) solution that the unoriented varifold $\widehat\mu_t$, corresponding to $\mu_t$  (identifying antipodal points $\pm p$ on $\S^{d-1}$), has integer multiplicity for almost any $t\in(0,\infty)$. The importance of the integer multiplicity of the varifold lies in the possibility of applying Brakke's orthogonality theorem so that Remark \ref{brakke} is in force. We will then exploit identity \eqref{basciWSuniq} in the rest of the following weak-strong uniqueness proof. We point out that, to perform this argument under this extra assumption on the integrality of the varifold, one can also use the gradient-flow calibration proposed in \cite{Timvol,Timkroemer}, which does not need to have property \eqref{normal} for the extended velocity vector field $B$.

We only discuss the parts which are considerably different from the techniques adopted in Section \ref{sec:C}.
Following the same arguments, we can obtain \eqref{estA}, and then, using  \eqref{ctrl}, \eqref{ptensorp}, \eqref{final1}, and \eqref{f1} (see also \eqref{estB}) we upgrade it to
\begin{align}
	\nonumber\Erel{T}&\leq E[\mu_0,\xhi_0|\mA(0)]- \frac{1}{2} \int_{0}^{T} \int_{\Omega} |V|^2 \d\omega_t \dt
	- \frac{1}{2} \int_{0}^{T} \int_{\Omega} \abs{\bH-\lambda\media}^2 \d\omega_t \dt
	\\&\nonumber
	\quad+\int_{0}^{T} \int_{\pA(t)} (\bn \cdot B)\Div \xi \d\H^{d-1} \dt-\int_{0}^{T} \int_{\Omega\times \mathbb S^{d-1}}(Id-p\otimes p):\nabla B \d\mu_t \dt\\&
	\quad - \int_{0}^{T} \int_{\Omega} V(\Div\xi -\lambdas)\d\omega_t \dt+C\int_0^T\Erel{t}\dt.
	\label{estB1}
\end{align}
Concerning the dissipative terms in the De Giorgi's inequality \eqref{eq:DeGiorgiInequality}, from \eqref{estB} we have, recalling the compatibility condition \eqref{eq:compatibility}, and the identity \eqref{eq:RadonNikodymProperty}, 

\begin{align}
	&\nonumber- \frac{1}{2} \int_{0}^{T} \int_{\Omega} |V|^2 \d\omega_t \dt
	- \int_{0}^{T} \int_{\Omega} V(\Div\xi -\lambdas)\d\omega_t \dt+\int_{0}^{T} \int_{\pA(t)} (\bn \cdot B)\Div \xi \d\H^{d-1} \dt\nonumber\\&\nonumber
	=- \frac{1}{2} \int_{0}^{T} \int_{\Omega} |V+\Div \xi-\lambdas|^2 \d\omega_t \dt+\frac12 \int_{0}^{T} \int_{\RevX{\supp \xi(t)}} (\Div\xi-\lambdas)^2 \d\omega_t \dt\\&\nonumber\quad+\frac12 \int_{0}^{T}\normmm{\lambdas}^2 \int_{\RevX{\Omega\setminus\supp\xi(t)}}1 \d\omega_t \dt+\int_{0}^{T} \int_{\pA(t)} (\bn \cdot B)\Div \xi \d\H^{d-1} \dt\\&\nonumber=
	- \frac{1}{2} \int_{0}^{T} \int_{\Omega} |V+\Div \xi-\lambdas|^2 \d\omega_t \dt+\frac12 \int_{0}^{T} \int_{\RevX{\supp \xi(t)}} (\Div\xi-\lambdas+B\cdot\xi)^2 \d\omega_t \dt\\&\nonumber\quad+\frac12 \int_{0}^{T}\normmm{\lambdas}^2 \int_{\RevX{\Omega\setminus\supp\xi(t)}}1 \d\omega_t \dt- \int_{0}^{T} \int_{\Omega} (\Div\xi-\lambdas)B\cdot \xi \d\omega_t \dt\\&\quad\nonumber+\int_{0}^{T} \int_{\pA(t)} (\bn \cdot B)\Div \xi \d\H^{d-1} \dt-\frac12 \int_{0}^{T} \int_{\Omega} (B\cdot \xi)^2 \d\omega_t \dt\\&\nonumber
	=	- \frac{1}{2} \int_{0}^{T} \int_{\Omega} |V+\Div \xi-\lambdas|^2 \d\omega_t \dt+\frac12 \int_{0}^{T} \int_{\RevX{\supp \xi(t)}} (\Div\xi-\lambdas+B\cdot\xi)^2 \d\omega_t \dt\\&\nonumber\quad +\frac12 \int_{0}^{T}\normmm{\lambdas}^2 \int_{\RevX{\Omega\setminus\supp\xi(t)}}1 \d\omega_t \dt+\int_{0}^{T} \int_{\Omega} (\Div\xi) B\cdot (\rho\bn-\xi) \d\omega_t \dt\nonumber\\&\quad +\int_{0}^{T} \int_{\Omega} \lambdas B\cdot \xi \d\omega_t \dt-\frac12 \int_{0}^{T} \int_{\Omega} (B\cdot \xi)^2 \d\omega_t \dt\nonumber\\&\nonumber
	=	- \frac{1}{2} \int_{0}^{T} \int_{\Omega} |V+\Div \xi-\lambdas|^2 \d\omega_t \dt+\frac12 \int_{0}^{T} \int_{\RevX{\supp \xi(t)}} (\Div\xi-\lambdas+B\cdot\xi)^2 \d\omega_t \dt\\&\quad +\frac12 \int_{0}^{T}\normmm{\lambdas}^2 \int_{\RevX{\Omega\setminus\supp\xi(t)}}1 \d\omega_t \dt+\int_{0}^{T} \int_{\Omega\times \bS} (-\Div\xi+\lambdas) B\cdot (\xi-p) \d\mu_t \dt\nonumber\\&\nonumber\quad +\int_{0}^{T} \int_{\Omega\times \bS} \lambdas B\cdot p \d\mu_t \dt-\frac12 \int_{0}^{T} \int_{\Omega} (B\cdot \xi)^2 \d\omega_t \dt\nonumber\\&
	\stackrel{\eqref{geometric},\eqref{coerciv2},\eqref{outsupp}}{\leq} 	- \frac{1}{2} \int_{0}^{T} \int_{\Omega} |V+\Div \xi-\lambdas|^2 \d\omega_t \dt +\int_{0}^{T} \int_{\Omega\times \bS} (-\Div\xi+\lambdas) B\cdot (\xi-p) \d\mu_t \dt\nonumber\\&\quad+\int_{0}^{T} \int_{\Omega\times \bS} \lambdas B\cdot p \d\mu_t \dt-\frac12 \int_{0}^{T} \int_{\Omega} (B\cdot \xi)^2 \d\omega_t \dt+C\int_0^T\RevX{(1+\normmm{\lambdas}^2)}\Erel{t}\dt.
	\label{fundamental2}
\end{align}
The main difference in the argument with respect to Section \ref{sec:C} is the estimate for the mean curvature term in the inequality \eqref{estB}, which is the only part of the proof strongly relying on property \eqref{normal}. Here we thus need to significantly modify the approach, since we do not want to use property \eqref{normal} of the vector $B$, i.e., we do not assume that $B$ is orthogonal to the interface $\cI(t)$. Let us first observe that, thanks to \eqref{equiv2} as well as to \eqref{basciWSuniq}, which holds due to the integrality assumption of the varifold, we have $\rho(x,t)^2=\normmm{\media}^2$ and 

\begin{align}
	\nonumber &	- \frac{1}{2} \int_{0}^{T} \int_{\Omega} \abs{\bH-\lambda\media}^2 \d\omega_t \dt\\& \nonumber=	- \frac{1}{2} \int_{0}^{T} \int_{\Omega} \abs{\bH-\lambda\media}^2 (1-\rho^2)\d\omega_t \dt	- \frac{1}{2} \int_{0}^{T} \int_{\Omega} \abs{\bH-\lambda\media}^2 \normmm{\media}^2\d\omega_t \dt\\&
	= 	- \frac{1}{2} \int_{0}^{T} \int_{\Omega} \abs{\bH-\lambda\media}^2 (1-\rho^2)\d\omega_t \dt	- \frac{1}{2} \int_{0}^{T} \int_{\Omega} \abs{(\bH-\lambda\media)\cdot \media}^2 \d\omega_t \dt.
	\label{curvature}
\end{align}
Additionally, recalling property \eqref{eq:weakCurvature}, together with \eqref{basciWSuniq}, entailing $ (\bH-\lambda\media)\cdot (Id-\media\otimes \media) B=0$,
\begin{align}
	&\nonumber -\int_{0}^{T} \int_{\Omega\times \mathbb S^{d-1}}(Id-p\otimes p):\nabla B \d\mu_t \dt=\int_0^T\int_\Omega \bH\cdot B\d\omega_t\dt \\&\nonumber=\int_0^T\int_\Omega (\bH-\lambda\media)\cdot B\d\omega_t\dt+\int_0^T\int_\Omega \lambda\media\cdot B\d\omega_t\dt\\&\nonumber 
	=\int_0^T\int_\Omega (\bH-\lambda\media)\cdot (Id-\media\otimes \media) B\d\omega_t\dt\\&\nonumber\quad +\int_0^T\int_\Omega ((\bH-\lambda\media)\cdot \media)(\media\cdot B)\d\omega_t\dt+\int_0^T\int_\Omega \lambda\media\cdot B\d\omega_t\dt\\&
	=\int_0^T\int_\Omega ((\bH-\lambda\media)\cdot \media)(\media\cdot B)\d\omega_t\dt+\int_0^T\int_\Omega \lambda\media\cdot B\d\omega_t\dt,
\end{align}
so that, also using \eqref{curvature}, we obtain, also applying Cauchy-Schwarz and Young's inequalities, 
\begin{align*}
	\nonumber	&- \frac{1}{2} \int_{0}^{T} \int_{\Omega} \abs{\bH-\lambda\media}^2 \d\omega_t \dt-\int_{0}^{T} \int_{\Omega\times \mathbb S^{d-1}}(Id-p\otimes p):\nabla B \d\mu_t \dt\\& \nonumber=
	- \frac{1}{2} \int_{0}^{T} \int_{\Omega} \abs{\bH-\lambda\media}^2 (1-\rho^2)\d\omega_t \dt-\frac12	\int_0^T\int_\Omega \normmm{(\bH-\lambda\media)\cdot \media- B\cdot \xi}^2\d\omega_t\dt\nonumber\\&
	\nonumber\quad +\frac12	\int_0^T\int_\Omega ({B\cdot \xi})^2\d\omega_t\dt+\int_0^T\int_{\Omega\times\mathbb S^{d-1}} \left((\bH-\lambda\media)\cdot \media-B\cdot\xi\right)((p-\xi)\cdot B)\d\mu_t\dt\\&\quad +\int_0^T\int_\Omega \left(B\cdot\xi\right)((\media-\xi)\cdot B)\d\omega_t\dt+\int_0^T\int_\Omega \lambda\media\cdot B\d\omega_t\dt\nonumber\\&
	\leq - \frac{1}{2} \int_{0}^{T} \int_{\Omega} \abs{\bH-\lambda\media}^2 (1-\rho^2)\d\omega_t \dt-\frac12	\int_0^T\int_\Omega \normmm{(\bH-\lambda\media)\cdot \media- B\cdot \xi}^2\d\omega_t\dt\nonumber\\&
	\nonumber\quad +	\frac12\int_0^T\int_\Omega ({B\cdot \xi})^2\d\omega_t\dt+\frac14\int_0^T\int_\Omega \normmm{(\bH-\lambda\media)\cdot \media- B\cdot \xi}^2\d\omega_t\dt\nonumber\\&\quad +\int_0^T\int_{\Omega\times \bS} \normmm{p-\xi}^2\normmm{B}^2\d\mu_t\dt+\int_0^T\int_{\Omega\times \bS} \left(B\cdot\xi\right)((p-\xi)\cdot B)\d\mu_t\dt \\&\quad \nonumber
	+\int_0^T\int_{\Omega\times \bS} \lambda p\cdot B\d\mu_t\dt,
\end{align*}
where we used in the last step $\int_0^T\int_\Omega \lambda\media\cdot B\d\omega_t\dt=\int_0^T\int_{\Omega\times \bS} \lambda p\cdot B\d\mu_t\dt$.
Then, we infer 
\begin{align}
	\nonumber	&- \frac{1}{2} \int_{0}^{T} \int_{\Omega} \abs{\bH-\lambda\media}^2 \d\omega_t \dt-\int_{0}^{T} \int_{\Omega\times \mathbb S^{d-1}}(Id-p\otimes p):\nabla B \d\mu_t \dt\\& \nonumber
	\stackrel{\eqref{coerciv1}}{\leq} - \frac{1}{2} \int_{0}^{T} \int_{\Omega} \abs{\bH-\lambda\media}^2 (1-\rho^2)\d\omega_t \dt-\frac14	\int_0^T\int_\Omega \normmm{(\bH-\lambda\media)\cdot \media- B\cdot \xi}^2\d\omega_t\dt\nonumber\\&
	\nonumber\quad +	\frac12\int_0^T\int_\Omega ({B\cdot \xi})^2\d\omega_t\dt+\int_0^T\int_{\Omega\times \bS} \left(B\cdot\xi\right)((p-\xi)\cdot B)\d\mu_t\dt 
	\nonumber
	\\&\quad+\int_0^T\int_{\Omega\times \bS} \lambda p\cdot B\d\mu_t\dt+C\int_0^T\Erel{t}\dt.\label{mostdifficult}
\end{align}
We can now combine \eqref{estB1}, \eqref{fundamental2}, and \eqref{mostdifficult}, to obtain, after an application of Cauchy-Schwarz and Young's inequalities and exploiting the approximate geometric evolution equation, 

\begin{align}
	&\nonumber\Erel{T}\\&\nonumber\leq E[\mu_0,\xhi_0|\mA(0)]	- \frac{1}{2} \int_{0}^{T} \int_{\Omega} |V+\Div \xi-\lambdas|^2 \d\omega_t \dt\nonumber\\&\quad \nonumber- \frac{1}{2} \int_{0}^{T} \int_{\Omega} \abs{\bH-\lambda\media}^2 (1-\rho^2)\d\omega_t \dt-\frac14	\int_0^T\int_\Omega \normmm{(\bH-\lambda\media)\cdot \media- B\cdot \xi}^2\d\omega_t\dt\\&\nonumber \quad-\int_{0}^{T} \int_{\Omega\times \bS} (B\cdot \xi+\Div\xi-\lambdas) B\cdot (\xi-p) \d\mu_t \dt\\&\quad+\int_{0}^{T} \int_{\Omega\times \bS} (\lambdas +\lambda)B\cdot p \d\mu_t \dt+C\int_0^T\RevX{(1+\normmm{\lambdas}^2)}\Erel{t}\nonumber\dt\\&
	\stackrel{\eqref{geometric}}{\leq} E[\mu_0,\xhi_0|\mA(0)]	- \frac{1}{2} \int_{0}^{T} \int_{\Omega} |V+\Div \xi-\lambdas|^2 \d\omega_t \dt\nonumber\\&\quad \nonumber- \frac{1}{2} \int_{0}^{T} \int_{\Omega} \abs{\bH-\lambda\media}^2 (1-\rho^2)\d\omega_t \dt-\frac14	\int_0^T\int_\Omega \normmm{(\bH-\lambda\media)\cdot \media- B\cdot \xi}^2\d\omega_t\dt\\&\nonumber \quad+\int_{0}^{T} \int_{\RevX{\supp B(t)}\times \bS} O(\dist(\cdot,\cI(t))^2)\d\mu_t \dt+\int_{0}^{T} \int_{\Omega\times \bS} \normmm{B}^2\normmm{\xi-p}^2 \d\mu_t \dt\\&\quad+\int_{0}^{T} \int_{\Omega\times \bS} (\lambdas +\lambda)B\cdot p \d\mu_t \dt+C\int_0^T\RevX{(1+\normmm{\lambdas}^2)}\Erel{t}\nonumber\dt\\&
	\stackrel{\eqref{coerciv1}}{\leq} \nonumber E[\mu_0,\xhi_0|\mA(0)]	- \frac{1}{2} \int_{0}^{T} \int_{\Omega} |V+\Div \xi-\lambdas|^2 \d\omega_t \dt\nonumber\\&\quad \nonumber- \frac{1}{2} \int_{0}^{T} \int_{\Omega} \abs{\bH-\lambda\media}^2 (1-\rho^2)\d\omega_t \dt-\frac14	\int_0^T\int_\Omega \normmm{(\bH-\lambda\media)\cdot \media- B\cdot \xi}^2\d\omega_t\dt\\& \quad+\int_{0}^{T} \int_{\Omega\times \bS} (\lambdas +\lambda)B\cdot p \d\mu_t \dt+C\int_0^T\RevX{(1+\normmm{\lambdas}^2)}\Erel{t}\dt.
	\label{estC1}
\end{align}
Again, we need to find a nonlocal control for the Lagrange multipliers. Namely, by means of \eqref{Lagrange}, using the compatibility condition \eqref{eq:compatibility}, together with Gau\ss\ Theorem, we obtain
\begin{align}
	\nonumber&\int_{0}^{T} \int_{\Omega\times \bS} (\lambdas(t) +\lambda(t))B\cdot p \d\mu_t \dt\nonumber=\int_{0}^{T} \int_{\pA(t)} (\lambdas(t) +\lambda(t))B\cdot \bn \d\H^{d-1} \dt\\&
	\nonumber =-\int_{0}^{T} (\lambdas(t) +\lambda(t)) \int_{\Omega}\xhi_{A(t)}\Div B \dx \dt\nonumber
	=-\int_{0}^{T} \int_{\Omega} (\lambdas(t) +\lambda(t))(\xhi_{A(t)}-\xhi_{\mA(t)})\Div B \dx \dt\\&\nonumber
	\leq  \int_{0}^{T} \int_{\Omega} (\normmm{\lambdas(t)} +\normmm{\lambda(t)})\normmm{\xhi_{A(t)}-\xhi_{\mA(t)}}\normmm{\Div B} \dx \dt\\&
	\stackrel{\eqref{div},\eqref{coerciv4}}{\leq}  \int_{0}^{T} (\normmm{\lambdas(t)} +\normmm{\lambda(t)})\int_{\RevX{\supp B(t)}} \normmm{\xhi_{A(t)}-\xhi_{\mA(t)}}O(\dist(\cdot,\cI(t))) \dx \dt\nonumber\\&\leq  \int_{0}^{T} (\normmm{\lambdas(t)} +\normmm{\lambda(t)}) \Ebulk{t} \dt.\label{estlambda1}
\end{align}

To close the estimates we can now repeat word by word the estimates on the bulk error, which do not rely on property \eqref{normal}, to get \eqref{estF}. Summing this inequality with \eqref{estC1}, and exploiting also \eqref{estlambda1}, we finally end up with 
\begin{align*}
	&\Erel{T}+\Ebulk{T}\\&+\frac{1}{4} \int_{0}^{T} \int_{\Omega} |V+\Div \xi-\lambdas|^2+\frac{1}{2} \int_{0}^{T} \int_{\Omega} \abs{\bH-\lambda\media}^2 (1-\rho^2)\d\omega_t \dt\\&+\frac14	\int_0^T\int_\Omega \normmm{(\bH-\lambda\media)\cdot \media- B\cdot \xi}^2\d\omega_t\dt \d\omega_t \dt\\&\nonumber\leq E[\mu_0,\xhi_0|\mA(0)]+E_{\mathrm{bulk}}[\xhi_0|\mA(0)]	\nonumber\\&+C\int_0^T\RevX{(1+\normmm{\lambdas(t)}^2)}\Erel{t}\dt+C\int_0^T(1+\normmm{\lambda(t)}+\normmm{\lambdas(t)})\Ebulk{t}\dt,
\end{align*}
so that, since $\lambdas,\lambda\in L^2(0,T)$, we can apply Gronwall's Lemma and conclude the proof of \eqref{eq:stabilityglobal} also in this case.
\\

  \textbf{Acknowledgments.}
  \extra{The author thanks the anonymous referees for the appropriate comments and useful remarks. }The author wants to warmly thank Tim Laux for fundamental and inspiring discussions, especially concerning the new gradient-flow calibration of Theorem  \ref{thmcalibration}, and for giving some useful comments on a first version of this manuscript. Part of this contribution was completed while visiting Helmut Abels at the Faculty of Mathematics of the University of Regensburg, whose hospitality is kindly acknowledged.  Moreover, the author thanks Keisuke Takasao for some helpful discussions about the sharp interface limits of nonlocal Allen-Cahn equation during his visit to Regensburg.
  
  The author also gratefully acknowledges support from the Alexander von Humboldt Foundation for the stay in Regensburg. This research was funded by the Austrian Science Fund (FWF) \href{https://doi.org/10.55776/ESP552}{10.55776/ESP552}.
The author is also a member of Gruppo Nazionale per l’Analisi Matematica, la Probabilità e le loro Applicazioni (GNAMPA) of
Istituto Nazionale per l’Alta Matematica (INdAM).  
For open access purposes, the author has applied a CC BY public copyright license to
any author accepted manuscript version arising from this submission.
\\

\textbf{Conflict of interest.} There is no conflict of interests.
\\

\textbf{Data availability.} Since our paper only consists of mathematical proofs, no additional datasets have been
generated or analyzed. All relevant materials can be obtained from the references cited in our paper.

  \bibliographystyle{siam}
	\bibliography{Bib_VPMCF}
	
\end{document}